\input amstex
\documentstyle{amsppt}
\magnification 1200
{\catcode`\@=11\gdef\logo@{}}

\voffset-7mm
\pagewidth{135mm}


\loadbold

\font \Snf = cmss12
\font \snf = cmss10
\font \ssnf = cmss7 

\def\ort{\hbox{\lower0.3ex\hbox{$\scriptstyle\sim$}%
\hskip-0.687ex\lower-.12ex\hbox{\vrule height0.95ex depth0.25pt%
width0.095ex\hskip0.6ex}}}
\def\Ort{\hbox{\lower0.34ex\hbox{$\scriptstyle\sim$}%
\hskip-0.74ex\lower-.15ex\hbox{\vrule height1.25ex depth0.32pt%
width0.110ex\hskip0.75ex}}}


\def \Ee   {\text{\Snf{E}}}

\def \Mm   {\text{\snf{M}}}
\def \Nn   {\text{\snf{N}}}
\def \nn   {\text{\ssnf{N}}}
\def \mm   {\text{\ssnf{M}}}
\def \bN   {\Bbb N}
\def \sbN  {{}^*\Bbb N}
\def \sbNi  {{}^*{\Bbb N}_{\infty}}
\def \bZ   {\Bbb Z}
\def \sbZ  {{}^*\Bbb Z}
\def \bQ   {\Bbb Q}
\def \sbQ  {{}^*\Bbb Q}
\def \bR   {\Bbb R}
\def \sbR  {{}^*\Bbb R}
\def \IR   {{\Bbb I}{}{^*\Bbb R}}
\def \FR   {{\Bbb F}{\,}{^*\Bbb R}}
\def \bC   {\Bbb C}
\def \sbC  {{}^*\Bbb C}
\def \bT   {\Bbb T}
\def \sbT  {{}^*\Bbb T}
\def \IC   {{\Bbb I}{}{^*\Bbb C}}
\def \FC   {{\Bbb F}{\,}{^*\Bbb C}}
\def \IT   {{\Bbb I}{}{^*\Bbb T}}
\def \IiC  {{\Bbb I}_{_\infty}{\!}{^*\Bbb C}}
\def \FiC  {{\Bbb F}_{_\infty}{\!}{^*\Bbb C}}
\def \IpC  {{\Bbb I}_{p}{\!}{^*\Bbb C}}
\def \FpC  {{\Bbb F}_{p}{\!}{^*\Bbb C}}

\def \FjC  {{\Bbb F}_{1}{\!}{^*\Bbb C}}

\def \InC  {{\Bbb I}_{_{\nn}}{\!}{^*\Bbb C}}
\def \FnC  {{\Bbb F}_{_{\nn}}{\!}{^*\Bbb C}}
\def \bV   {\Bbb V}
\def \BB   {\Cal B}
\def \DD   {\Cal D}
\def \PP   {\Cal P}
\def \wPP  {\widetilde{\Cal P}}

\def \II   {\Cal I}
\def \KK   {\Cal K}
\def \LL   {\Cal L}

\def \CCb  {\Cal C_{\operatorname{b}}}
\def \CCbu {\Cal C_{\operatorname{bu}}}

\def \CCc  {\Cal C_{\operatorname{c}}}
\def \CMp  {\Cal C\Cal M^p}

\def \CCo  {\Cal C_0}

\def \FF   {\Cal F}
\def \MM   {\Cal M}

\def \RR   {\Cal R}
\def \QQ   {\Cal Q}
\def \SS   {\Cal S}
\def \UU   {\Cal U}
\def \VV   {\Cal V}
\def \WW   {\Cal W}
\def \ab   {\boldkey a}

\def \fb   {\boldkey f}
\def \gb   {\boldkey g}
\def \ub   {\boldkey u}
\def \vb   {\boldkey v}
\def \xb   {\boldkey x}
\def \yb   {\boldkey y}

\def \mb   {\boldkey{m}}
\def \nb   {\boldkey{n}}
\def \mbd  {\boldkey{m}_d}
\def \Ab   {\bold A}
\def \Bb   {\bold B}
\def \Gb   {\bold G}
\def \sGb  {{}^*{\bold G}}
\def \dGb  {\widehat{\bold G}}
\def \ddGb {\widehat{\widehat{\bold G}}}
\def \sdGb {{}^*\widehat{\bold G}}
\def \IG   {{\Bbb I}{}{^{\,*}\bold G}}
\def \FG   {{\Bbb F}{}{^{\,*}\bold G}}
\def \IdG  {{\Bbb I}{}{^{\,*}\widehat{\bold G}}}
\def \FdG  {{\Bbb F}{}{^{\,*}\widehat{\bold G}}}
\def \Db   {\bold D}
\def \oDb  {\overline{\bold D}}
\def \Eb   {\bold E}
\def \Fb   {\bold F}
\def \Hb   {\bold H}
\def \Kb   {\bold K}
\def \Lbb  {\bold L}
\def \Nb   {\bold N}
\def \Sbb  {\bold S}
\def \sKb  {{}^*\bold K}
\def \Ub   {\bold U}
\def \sUb  {{}^*\bold U}
\def \Vb   {\bold V}
\def \sVb  {{}^*\bold V}
\def \Wb   {\bold W}
\def \Xb   {\bold X}
\def \Yb   {\bold Y}
\def \Id   {\operatorname{Id}}
\def \C    {\operatorname{C}}
\def \Cb   {\operatorname{C}_{\operatorname{b}}}
\def \Cbo  {\operatorname{C}_0}
\def \Cbc  {\operatorname{C}_{\operatorname{c}}}
\def \Cbu  {\operatorname{C}_{\operatorname{bu}}}
\def \Mb   {\operatorname{M}}
\def \Lb   {\operatorname{L}}
\def \Lbp  {\operatorname{L}^p}
\def \Lbq  {\operatorname{L}^q}

\def \Gama  {\varGamma}
\def \Gamb  {\boldsymbol{\Gamma}}

\def \Dela  {\varDelta}
\def \Delb  {\boldsymbol{\Delta}}
\def \Thb   {\boldsymbol{\Theta}}
\def \vFi   {\varPhi}

\def \Omb   {\boldsymbol{\Omega}}
\def \Yps   {\varUpsilon}
\def \Ypb   {\boldsymbol{\Upsilon}}
\def \vfi   {\varphi}
\def \pphi  {\boldsymbol{\phi}}
\def \alfa  {\alpha}
\def \gama  {\gamma}
\def \ggam  {\boldsymbol{\gamma}}
\def \cchi  {\boldsymbol{\chi}}
\def \mmu   {\boldsymbol{\mu}}

\def \oom   {\boldsymbol{\omega}}
\def \lam   {\lambda}
\def \eps   {\varepsilon}
\def \ro    {\varrho}
\def \sig   {\sigma}
\def \eeta  {\boldsymbol{\eta}}
\def \th    {\vartheta}
\def \tth   {\boldsymbol{\theta}}

\def \Go    {G_0}
\def \Gf    {G_{\operatorname{f}}}
\def \Ho    {H_0}
\def \Hf    {H_{\operatorname{f}}}
\def \Xf    {X_{\operatorname{f}}}
\def \Yf    {Y_{\operatorname{f}}}
\def \dG    {\widehat G}
\def \ddG   {\widehat{\!\widehat G}}
\def \dH    {\widehat H}
\def \arg   {\operatorname{arg}}
\def \card  {\operatorname{card}}
\def \dom   {\operatorname{dom}}
\def \supp  {\operatorname{supp}}
\def \st    {\operatorname{st}}
\def \Hom   {\operatorname{Hom}}
\def \Re    {\operatorname{Re}}

\def \dif   {\operatorname{d}\!}
\def \Ns    {\operatorname{Ns}}
\def \Mon   {\operatorname{Mon}}
\def \Bohr  {\operatorname{Bohr}}
\def \Spec  {\operatorname{Spec}}
\def \SO    {\operatorname{SO}}
\def \et    {\ \,\&\,\ }
\def \imp   {\ \Rightarrow\ }
\def \iff   {\ \Leftrightarrow\ }
\def \exs   {\exists\,}
\def \all   {\forall\,}
\def \alli  {\forall^{\,\operatorname{int}}}
\def \lv   {\left|}
\def \rv   {\right|}
\def \LV   {\left\|}
\def \RV   {\right\|}

\def \br   {\langle}
\def \kt   {\rangle}
\def \nin  {\notin}
\def \sbs  {\subseteq}
\def \sps  {\supseteq}
\def \nsg  {\vartriangleleft}
\def \apr  {\approx}
\def \mto  {\mapsto}

\def \iso  {\cong}
\def \sms  {\smallsetminus}
\def \cd   {\cdot}
\def \cx   {\times}
\def \co   {\circ}
\def \rst  {\!\restriction\!}
\def \wh   {\widehat}
\def \wtl  {\widetilde}
\def \ovl  {\overline}
\def \be   {\flat}
\def \kr   {\sharp}
\def \ee   {{\operatorname{e}}}
\def \ii   {{\operatorname{i}}}
\def \dh   {\hat d}
\def \sA  {{}^{*\!}A}
\def \sB  {{}^{*\!}B}
\def \sD  {{}^{*\!}D}
\def \sU  {{}^{*\!}U}
\def \sX  {{}^{*\!}X}
\def \sXb {{}^*{\bold X}}

\topmatter

\title
Gordon's Conjectures: \\
Pontryagin-van\,Kampen duality and \\
the Fourier transform in hyperfinite setting 
\endtitle

\rightheadtext{Gordon's Conjectures}

\author
Pavol Zlato\v{s}
\endauthor

\dedicatory
Dedicated to the memory of Petr Vop\v{e}nka (1935--2015)
\enddedicatory

\address
Faculty of Mathematics Physics and Informatics, Comenius University,
\newline\phantom{xx\,}%
Mlynsk\'{a} dolina, 842\,48~Bratislava, Slovakia
\endaddress

\email
{\tt zlatos{\@}fmph.uniba.sk}
\endemail

\abstract
Using the ideas of E.\,I.~Gordon we present and farther advance
an approach, based on nonstandard analysis, to simultaneous
approximations of locally compact abelian groups and their duals
by (hyper)finite abelian groups, as well as to approximations of
various types of Fourier transforms on them by the discrete Fourier
transform. Combining some methods of nonstandard analysis and
additive combinatorics we prove the three Gordon's Conjectures
which were open since 1991 and are crucial both in the formulations
and proofs of the LCA groups and Fourier transform approximation
theorems.
\endabstract



\thanks
Research supported by the Scientific Grant Agency of Slovak Republic VEGA
\endthanks

\endtopmatter
\vskip-15pt

\document
\centerline{\smc Contents}
\medskip
{\parindent 21pt
\item{0.} Introduction \hfill{2}
\item{1.} Nonstandard analysis \hfill{5}
\itemitem{1.1.} General setting \hfill{5}
\itemitem{1.2.} Bounded monadic spaces \hfill{7}
\itemitem{1.3.} Functional spaces 1: Continuous functions \hfill{12}
\itemitem{1.4.} Functional spaces 2: Loeb measures and
                Lebesgue spaces \hfill{15}
\itemitem{1.5.} Bounded monadic groups \hfill{20}
\item{2.} Pontryagin-van\,Kampen duality in hyperfinite setting \hfill{31}
\itemitem{2.1.} The dual triplet \hfill{31}
\itemitem{2.2.} Fourier transforms, Bohr sets and spectral sets
                in finite abelian groups \hfill{36}
\itemitem{2.3.} The dual triplet continued: Normal multipliers and
\itemitem{}     proofs of Gordon's Conjectures 1 and 2 \hfill{42}
\itemitem{2.4.} Some (mainly) standard equivalents:
                Hrushovski style theorems \hfill{45}
\itemitem{2.5.} Simultaneous approximation of an LCA group
                and its dual \hfill{50}
\item{3.} The Fourier transform in hyperfinite dimensional setting \hfill{64}
\itemitem{3.1.} A characterization of liftings \hfill{64}
\itemitem{3.2.} The Smoothness-and-Decay Principle \hfill{68}
\itemitem{3.3.} Hyperfinite dimensional approximations of the Fourier transforms:
\itemitem{}     Generalized Gordon's Conjecture 3 \hfill{71}
\itemitem{3.4.} Standard consequences: Simultaneous approximation of a function
\itemitem{}     and its Fourier transform \hfill{74}
\item{} References \hfill{84}
\item{}}

{}\phantom{x}
{}\bigskip

\head
0. Introduction
\endhead
{}\bigskip\bigskip

\flushpar
Locally compact abelian groups (briefly, LCA groups) and
the celebrated Pontryagin-van\,Kampen duality theorem provide a
general background on which all the particular instances of the
(commutative) Fourier transform can be treated in a uniform way.
Among the most important cases they include Fourier series of
periodic functions $f\:\bR \to \bC$ with a fixed period
$T > 0$, Fourier transforms of functions $f\:\bR \to \bC$ and
$f\:\bR^n \to \bC$, the semidiscrete Fourier transform of sequences
$f\:\bZ \to \bC$, and, of course, the discrete Fourier transform
of functions ($n$-dimensional vectors) $f\:\bZ_n \to \bC$ or,
more generally, of functions $f\:G \to \bC$ defined on an arbitrary
finite abelian group $G$ (cf\. \cite{HR1}, \cite{HR2}, \cite{Rd2},
\cite{Tr}).

For a finite abelian group $G$ its dual group $\dG = \Hom(G,\bT)$,
where $\bT$ denotes the compact multiplicative group of complex
units, is isomorphic (though not canonically) to $G$, the
$|G|$-dimensional  vector space $\bC^G$ is endowed with the
Hermitian inner product
$$
\br f,\,g\kt_d = d \,\sum_{x \in G} f(x)\,\ovl{g(x)}
$$
where $d > 0$ is some scaling or normalizing coefficient, the
characters $\gama \in \dG$ form an orthogonal basis in $\bC^G$
and the discrete Fourier transform (DFT) $\FF\:\bC^G \to \bC^{\dG}$
is defined as the inner (or scalar) product
$$
\FF(f)(\gama) = \wh f(\gama) = \br f,\,\gama\kt_d \,,
$$
for $f \in \bC^G$, $\gama \in \dG$. Once the inner product
$\br\vfi,\psi\kt_{\dh}$ on $\bC^{\dG}$ is defined using the
adjoint scaling coefficient $\dh = (d\lv G\rv)^{-1}$, we have
the Fourier inversion formula
$$
f = \dh\,\sum_{\gama \in \dG} \wh f(\gama)\,\gama \,,
$$
and the Plancherel identity 
$$
\br f,\,g\kt_d = \bigl\br\wh f,\,\wh g\,\bigr\kt_{\dh} \,,
$$
turning the DFT \,$\FF\:\bC^G \to \bC^{\dG}$ into a linear
isometry of unitary spaces.

For general LCA groups the picture is by far not so simple. The
dual group $\dG$ of \,$G$ consists of all continuous homomorphisms
(characters) $\gama\:G \to \bT$, and the Fourier transform is
primarily defined on the Lebesgue space $\Lb^1(G) = \Lb^1(G,m)$,
where $m = m_G$ is the Haar measure on $G$, as the bounded linear
operator $\FF\:\Lb^1(G) \to \Cbo\bigl(\dG\bigr)$ given by
$$
\FF(f)(\gama) = \wh f(\gama) = \int f\,\ovl\gama \dif m \,,
$$
for $f \in \Lb^1(G)$, $\gama \in \dG$. It can be extended to the
so called the Fourier-Stieltjes transform
$\FF\:\Mb(G) \to \Cbu\bigl(\dG\bigr)$
from the Banach space $\Mb(G) \sps \Lb^1(G)$ of all complex-valued
regular Borel measures on $G$ with finite total variation to the
Banach space $\Cbu\bigl(\dG\bigr)$ of all bounded uniformly
continuous functions $\dG \to \bC$, defined by
$$
\FF(\mu)(\gama) = \wh\mu(\gama) = \int \ovl\gama \dif\mu \,,
$$
for $\mu \in \Mb(G)$.
\pagebreak

Using the density of the intersection $\Lb^1(G) \cap \Lbp(G)$ in
the Lebesgue space $\Lbp(G)$ with respect to its norm
$\LV\cd\RV_{_p}$, the Fourier transform can also be extended to
a bounded linear operator $\FF\:\Lbp(G) \to \Lbq\bigl(\dG\bigr)$
for $p \in (1,2\,]$ and the adjoint exponent
$q = p/(p-1) \in [\,2,\infty)$. Under a proper normalization of
the Haar measure $m_{\dG}$ on the dual group $\dG$ we have the
Fourier inversion formula
$$
f = \int \wh f(\gama)\,\gama \dif m_{\dG}
$$
(both with respect to the supremum norm $\|\!\cd\!\|_{_\infty}$
and the $\Lb^p$-norm $\|\!\cd\!\|_{_p}$) just for
$f \in \Lb^p(G) \cap \FF\bigl[\Mb\bigl(\dG\bigr)\bigr]
\sbs \Lb^p(G) \cap \Cbu(G)$, with $\FF$ denoting the
Fourier-Stieltjes transform $\Mb\bigl(\dG\bigr) \to \Cbu(G)$,
here.

For $p=q=2$  we obtain the isometric linear isomorphism
$\FF\:\Lb^2(G) \to \Lb^2\bigl(\dG\bigr)$ of Hilbert spaces,
called the Fourier-Plancherel transform. Then we have the
Plancherel identity
$$
\br f,\,g\kt = \int f\,\ovl g \dif m_G
= \int \wh f\ \ovl{\wh g} \dif m_{\dG}
= \bigl\br\wh f,\,\wh g\,\bigr\kt
$$
(just) for $f,g \in \Lb^2(G)$. Unfortunately, unless $G$ is
compact, $\dG \cap \Lb^2(G) = \emptyset$ and the inner product
$\br\gama,\chi\kt$ of characters $\gama,\chi \in \dG$ is never
defined, so that one can speak of the orthogonal basis formed by
the characters at most in a metaphorical sense.

Taking additionally into account that the Fourier transform on
finite abelian groups can be computed using the extremely fast and
powerful algorithms of the Fast Fourier Transform, there naturally
arises the following question:

\proclaim\nofrills{}\par
Given any LCA group $G$, isn't there some ``universal extension'',
encompassing all the spaces $\Lbp(G)$ and $\Mb(G)$, and a uniform
scheme defining the Fourier transform on this extension, covering
all the above mentioned particular Fourier transforms, like if
\,$G$ were finite?
\endproclaim

The main goal of this paper is to provide arguments that namely
Nonstandard Analysis offers not only a reasonable and satisfactory
but also a fairly elegant solution to this question, as well
as several additional insights. Our approach is based on the idea
of an infinitesimal approximation of any LCA group $G$ \,by
a hyperfinite abelian group, yielding infinitesimal approximations
of all the above mentioned Fourier transforms by the hyperfinite
dimensional DFT on (the hyperfinite dimensional vector space over)
the approximating hyperfinite group. In standard terms this means
that each LCA group $G$ admits an approximating system consisting
of finite abelian groups, yielding approximations of the particular
Fourier transforms on the $\Lbp(G)$s and $\Mb(G)$ by the DFTs on
(the finite dimensional vector spaces over) the finite approximating
groups.

For the first time the methods of nonstandard analysis were applied
in the study of Fourier series of functions $\bT \to \bC$ by
Luxemburg \cite{Lx3}. The key idea consisted in embedding the
group of integers $\bZ$ into the hyperfinite cyclic group $\bZ_n$,
where $n \in \sbN \sms \bN$, and infinitesimal approximation of the
group of complex units $\bT$ by the hyperfinite subgroup
$\bigl\{\ee^{2\pi\ii k/n}\: k \in \bZ_n\bigr\} \iso \bZ_n$ of its
nonstandard extension $\sbT$. The first treatment of abstract
(commutative) harmonic analysis by nonstandard methods in full
generality is due to Gordon. In a~series of works culminating in
\cite{Go1}, \cite{Go2} he elaborated a~nonstandard approach to
approximations of LCA groups by hyperfinite abelian groups,
formulated a version of Pontryagin-van\,Kampen duality for them
and developed an approach to the approximation of the classical
Fourier-Plancherel transform $\Lb^2(G) \to \Lb^2\bigl(\dG\bigr)$
by the DFT on the approximating hyperfinite group. At the same time,
he formulated three rather fundamental conjectures in \cite{Go1}
which remained open since 1991 until 2012.

In the present paper we will recapitulate Gordon's approach and
results, introducing some conceptual and notational modifications,
based mainly on some results from \cite{ZZ}, and prove the three
Gordon's conjectures, generalizing the last one to the case of the
classical Fourier transforms on all of the $\Lbp(G)$ and $\Mb(G)$
spaces. The methods used in the proofs are mainly combinations of
various methods of nonstandard analysis and harmonic analysis with
the Fourier analytic methods of additive combinatorics by
Green-Ruzsa \cite{GR} and Tao-Vu \cite{TV}.

Finally, we will present classical (standard) equivalents to some of
the obtained nonstandard results. Should we encapsulate their moral
in a single sentence, the best we can do seems to be to phrase it as
a response to the question formulated in the title of the paper
by Epstein \cite{Ep}:
\medskip
\centerline{\it How well does the finite Fourier transform
approximate the Fourier transform?}
\medskip
The response there in the Abstract is {\it ``very well indeed'',}
and farther in the Conclusion also {\it ``as well as it possibly
could''}. We hope to convince the reader to agree finally with the
following:
\medskip
\centerline{\it Even better than one could ever hope.%
\footnote{In fact, in \cite{Ep} that question is asked for the
Fourier transform of periodic functions $\bR \to \bC$,
only. In our response we have in mind Fourier transforms on
arbitrary LCA groups.}}
\medskip
Preliminarily, we can unfold the above slogan in the following
imprecise and intuitive way:
{\parindent 21pt
\item{1.}
for every LCA group $\Gb$ one can find an ``arbitrarily good pair of
adjoint approximations'' of \,$\Gb$ by a finite abelian group $G$ and
of its dual group $\dGb$ by the dual $\dG$ of the finite group $G$;
\item{2.}
any of the Fourier transforms $\Mb(\Gb) \to \Cbu\bigl(\dGb\bigr)$,
\hbox{$\Lb^1(\Gb) \to \Cbo\bigl(\dGb\bigr)$,}
$\Lbp(\Gb) \to \Lbq\bigl(\dGb\bigr)$, for adjoint exponents
$1 < p \le 2 \le q < \infty$, can be ``arbitrarily well''
approximated by the discrete Fourier transform $\bC^G \to \bC^{\dG}$,
based on some adjoint approximations of \,$\Gb$, $\dGb$ by finite
abelian groups $G$, $\dG$, respectively.
\item{}
\vskip-12 pt}
\smallskip
\flushpar
Precise formulation of \,1 is the Strongly Adjoint Finite LCA Group
Approximation Theorem~2.5.8, while precise formulations of \,2 are
the three Finite Fourier Transform Approximation Theorems~3.4.4,
3.4.5 and 3.4.6. They are derived from their nonstandard
counterparts: the Adjoint Hyperfnite LCA Group Approximation
Theorem (Corollary~2.5.2), and the three Hyperfinite Dimensional
Fourier Transform Approximation Theorems~3.3.1, 3.3.2 and 3.3.4,
respectively. In their formulations as well as in their proofs
Gordon's Conjectures~1 and~2 are crucial. These are stated and
proved in Sections~2.1--2.3. Gordon's Conjecture~3 on hyperfinite
dimensional approximation of the Fourier-Plancherel transform
$\Lb^2(\Gb) \to \Lb^2\bigl(\dGb\bigr)$ is in fact a special
case of Theorem~3.3.2 (Corollary~3.3.3).

\definition{Acknowledgement}
In the end of this introductory part I would like to express
my deep indebtedness and gratitude to Zhenya Gordon for having
introduced me to the topic of nonstandard  approach to
Pontryagin-van\,Kampen duality and the Fourier transform, as well as
for many valuable discussions and permanent encouragement and
support when my long year wrestling with his Conjectures seemed
desperately hopeless.
\enddefinition
\newpage

{}\phantom{x}
{}\bigskip

\head
1. Nonstandard Analysis
\endhead
{}\bigskip\bigskip

\flushpar
The reader is assumed to have some basic acquaintance with nonstandard
analysis, including the nonstandard approach to topology and continuity
in terms of monads and equivalence relations of infinitesimal nearness,
the Loeb measure construction and internal Banach spaces and their
nonstandard hulls. Besides the original Robinson's book \cite{Rb},
the standard general references are, e.g., the monographs \cite{AFHL},
\cite{ACH} (mainly the parts \cite{Hn2} and \cite{Lb2}), \cite{Dv},
and \cite{Gb}. For Loeb measure also the survey \cite{Ct} can be
consulted. The canonical reference for nonstandard Banach space
theory is the paper \cite{HM}. Additionally, we refer to \cite{CK}
for the ultraproduct construction and its use in constructing
nonstandard extensions of first-order structures and superstructures,
as well as for the connections between certain properties of
ultrafilters and the nonstandard universes obtained via the
ultraproduct construction with respect to them.

\specialhead
1.1.~General setting
\endspecialhead
\flushpar
Our exposition takes place in a~nonstandard universe ${}^*\bV$
which is an elementary extension of a superstructure $\bV$ over
some set of individuals containing at least all (classical)
complex numbers and the elements of the topological space or
topological group, as well as index sets, etc., dealt with.
In particular, every standard mathematical (first-order)
structure $A \in \bV$ is embedded into its nonstandard extension
$\sA \in {}^*\bV$ via the mapping \hbox{$a \mto {}^*a\: A \to \sA$}
such that, for any formula $\Phi(x_1,\dots,x_n)$ in the language
of $A$ and elements $a_1,\dots,a_n \in A$, $\Phi(a_1,\dots,a_n)$
is satisfied in $A$ if and only if $\Phi({}^*a_1,\dots,{}^*a_n)$
is satisfied in $\sA$ ({\it transfer principle\/}). Whenever
there threatens no confusion we tend to identify $a \in A$ with
${}^*a \in \sA$, and to denote the corresponding operations and
relations in $A$ and $\sA$ by the same sign, dropping the $^*$ in
the latter. Similarly, to a function $f\: A \to B$ in $\bV$ there
(functorially) corresponds a function ${}^{*\!}f\: \sA \to \sB$ in
${}^*\bV$, etc.

In particular, we have the structures of {\it hypernatural numbers\/}
$\sbN$, {\it hyperintegers\/} $\sbZ$, {\it hyperrational numbers\/}
$\sbQ$, {\it hyperreal numbers\/} $\sbR$, and {\it hypercomplex numbers\/}
$\sbC$ with the usual (and, possibly, some additional) operations and
relations, extending the structures of natural numbers $\bN$, integers
$\bZ$, rational numbers $\bQ$, real numbers $\bR$, and complex numbers
$\bC$, respectively.

As the superstructure $\bV$ is {\it transitive\/}, i.e.,
$X \sbs \bV$ for any set $X \in \bV$, the same is true for
${}^*\bV$. The sets belonging to ${}^*\bV$ are called
{\it internal\/}; other subsets of \,${}^*\bV$ are called
{\it external\/}. Additionally, we assume that the nonstandard
universe ${}^*\bV$ is either $\kappa$-saturated for some
uncountable cardinal $\kappa$ or even {\it polysaturated}, i.e.,
$\kappa$-saturated for some $\kappa$ bigger than the cardinality
of any set in the original (standard) universe $\bV$. However,
for the sake of generality, we do not specify the saturation
degree $\kappa$ explicitly. Instead we use the term a {\it set\/}
or {\it system of admissible size\/} referring to (external)
subsets of the nonstandard universe with the (external)
cardinality $< \kappa$, and assume that the universe ${}^*\bV$
is {\it sufficiently saturated\/}, meaning that
$\bigcap \SS \ne \emptyset$ for any system of internal
sets $\SS \sbs {}^*\bV$ of admissible size with the
finite intersection property. For most applications an
$\aleph_1$-saturated nonstandard universe (i.e.,
$\kappa = \aleph_1$) would be sufficient; in that case
a~system of admissible size is simply a~countable one.

Internal sets $A$ which can be put into a one-to-one correspondence
via an internal bijection with sets of the form $\{1,\dots,n\}$
for some $n \in \sbN$ are called {\it hyperfinite\/}; in that
case $n = |A|$ is referred to as the {\it number of elements\/} of
$A$. Hyperfinite sets (briefly, HF sets) behave within the internal
context much like finite sets though, for $n \in \sbN \sms \bN$,
they are (externally) infinite.

Internal and, particularly, hyperfinite sets are the simplest objects
in the descriptive hierarchy of the at least in principle describable
external sets. Intuitively, they serve as mathematical models of
unambiguous, sharply demarcated groupings of objects. Next to them in
this hierarchy there are the {\it galactic} or $\Sigma^0_1(\kappa)$-sets
and the {\it monadic\/} or $\Pi^0_1(\kappa)$-sets, defined as unions and
intersections, respectively, of admissible size systems of internal sets.
From the formal point of view, it is namely these sets which make work the
powerful techniques of Nonstandard Analysis based on {\it saturation}.
Intuitively, galactic sets can be viewed as mathematical models of
groupings demarcated by some observational horizon toward which the
property possibly fades. Then the monadic sets are just complements of
galactic sets with respect to internal sets; they model negations of
properties demarcated by some observational horizon.

An example of a galactic set is the set $\bN$ of all (standard)
natural numbers, in other words, the set of all {\it finite\/}
elements of the internal set of all hypernatural numbers $\sbN$.
Then the {\it infinite\/} hypernatural numbers form the monadic
set $\sbNi = \sbN \sms \bN$. As $n+1 \in \bN$ for each $n \in \bN$,
and, similarly, $n-1 \in \sbNi$ for each $n \in \sbNi$, it is rather
a horizon than a sharp, firm boundary separating the finite and
infinite hypernatural numbers. At the same time it is worthwhile to
notice one important difference between the galactic subset $\bN$
and the monadic subset $\sbNi$ of \,$\sbN$ which is due to
{\it saturation\/}: There are sequences $(a_k)_{k \in \bN}$ of
\,{\it finite\/} natural numbers (e.g., all the strictly increasing
ones) having no continuation within $\bN$, i.e., any internal
sequence $(a_k)_{k \le n}$, where $n \in \sbNi$, extending
$(a_k)_{k \in \bN}$ necessarily contains an initial segment
$(a_k)_{k \le m}$ such that $m \in \sbNi$, $m \le n$ and
$a_k \in \sbNi$ for each infinite $k \le m$. On the other hand,
every sequence $(b_k)_{k \in \bN}$ of {\it infinite\/} natural numbers
has a continuation within $\sbNi$, i.e., there is an internal sequence
$(b_k)_{k \le n}$, where $n \in \sbNi$, extending it such that
$k \in \sbNi$ for each $k \le n$.

Further typical and fairly important examples of monadic sets are the
equivalence relations of indiscernibility or infinitesimal nearness,
arising in nonstandard models of topological spaces. Intuitively, two
objects are {\it indiscernible\/} within some method of comparison or
observation if they are behind the horizon of its {\it discernibility\/}.

Every function $f\:X \to Y$ is considered to be equal to the set
of ordered pairs $\bigl\{(x,f(x))\: x \in X\bigr\}$. If
\,$R \sbs X \cx Y$ is a relation, then a function $f$ is called a
{\it choice function\/} from $R$ on a set $A \sbs \dom R$ if
\,$A \sbs \dom f$ and $f\rst A \sbs R$, i.e., if
$(a,f(a)) \in R$ for each $a \in A$.

\proclaim{1.1.1.~Internal Choice Lemma}
Let $X$, $Y$ be internal sets in a sufficiently saturated
nonstandard universe, and $R \sbs X \cx Y$ be a relation
such that for every internal set $D \sbs \dom R$ the restriction
$R\rst D$ is a monadic set. Then for every galactic set
$A \sbs \dom R$ there exists an internal choice function $f$
from $R$ on $A$.
\endproclaim

\demo{Sketch of proof}
If $A$ is internal, then the monadic relation $R\rst A$ can be
written as the intersection $R\rst A = \bigcap_{i \in I} R_i$ of
admissibly many internal relations $R_i \sbs X \cx Y$ with common
domain $A$. Then, for any nonempty finite set $J \sbs I$, we readily
obtain an internal choice function $f_J$ on $A$ from the relation
$\bigcap_{i \in J} R_i$, by applying the {\it transfer principle\/}
to the axiom of choice. (For hyperfinite $A$ the axiom of choice
is even not needed in this point.) The existence of an internal
choice function $f$ from $R$ on $A$ follows by the virtue of
saturation.

If $A$ is a galactic set, then the desired conclusion follows from
the internal case by applying the saturation argument once again.
\enddemo

\specialhead
1.2.~Bounded monadic spaces
\endspecialhead
\flushpar
In view of \cite{Lx2}, \cite{Hn1} and \cite{Gn2} we accept as
a bare fact that there is no canonical way how to define the
{\it finite elements\/} in the nonstandard extension of a uniform
space. Instead of looking for the adequate definition in terms of
the uniform structure and standard elements we will treat the
{\it infinitesimal nearness\/} or {\it indiscernibility\/} on one
hand, and {\it finiteness\/} or {\it accessibility\/} on the other
hand as related but different phenomena to which there correspond
different primitive concepts. Our basic nonstandard objects, by
means of which we will study (sufficiently regular) topological
spaces, will be ordered triples of the form $(X,E,\Xf)$ where $X$
is an internal set, $E$ is a monadic equivalence relation on $X$
and $\Xf$ is a galactic subset of \,$X$ which is $E$-closed, i.e.,
$x \in \Xf$ and $(x,y) \in E$ imply $y \in \Xf$, for $x,y \in X$.
We will call them alternatively {\it bounded monadic spaces}, or
{\it IMG spaces\/} like in \cite{ZZ}, or {\it IMG triplets},
indicating that we do not consider this terminology as definitive.

Intuitively, $X$ is viewed as the underlying or ambient set of the
triplet, $E$ is the relation of indiscernibility or infinitesimal
nearness on $X$, and $\Xf$ is the set of elements of \,$X$
encompassed by some observational horizon. The elements of \,$\Xf$
will be briefly referred to as the {\it finite\/} or
{\it accessible\/} ones. To stress the role of the equivalence $E$
we will preferably write $x \apr y$ instead of $(x,y) \in E$, for
$x,y \in X$, and call the set
$$
E[x] = \{y \in X\: y \apr x\}
$$
of points indiscernible from the point $x \in X$ the
{\it $E$-monad\/} or just the {\it monad\/} of $x$. The
$E$-closeness of \,$\Xf$ can be now expressed as the
condition $E[x] \sbs \Xf$ for any $x \in \Xf$. The quotient
$$
\Xf/E = \Xf/\!\!\apr\ = \{E[x]\: x \in \Xf\}
$$
is called the {\it observable trace\/} of the triplet $(X,E,\Xf)$.

The restricted quotient mapping $\Xf \to \Xf/E$ reminds of the
standard part mapping $\Ns({}^*\Xb) \to \Xb$ in nonstandard
extensions of Hausdorff uniform spaces, sending every point
$\xb \in \Ns({}^*\Xb)$ to its standard part, i.e., the unique
element ${}^\co\xb = \st\xb \in \Xb \iso \Ns({}^*\Xb)/\!\!\apr$
infinitesimally close to $\xb$. In order to underline this analogy
(especially when viewing the monads as individual points and
forgetting about their ``sethood'') we introduce the notation
$E[x] = x^\be$ for the monad of $x \in X$, and
$$
A^\be = \bigl\{a^\be\: a \in A \cap \Xf\bigr\}
$$
for the {\it observable trace\/} of any set $A \sbs X$. In
particular, the observable trace $X^\be = \Xf^\be = \Xf/E$ of
the triplet $(X,E,\Xf)$ should not be confused with the full
quotient $X/E \sps X^\be$. Conversely, for any $\Yb \sbs X^\be$,
we call the following set the {\it pretrace\/} of \,$\Yb$:
$$
\Yb^\kr = \bigl\{x \in \Xf\: x^\be \in \Yb\bigr\}\,.
$$

Given an IMG triplet $(X,E,\Xf)$, it is an easy exercise in
{\it saturation\/} to show that for each internal relation $R \sps E$
on $X$ there is a symmetric internal relation $S \sps E$ on $X$ such
that $S \co S \sbs R$. Similarly, for any internal set $A \sbs \Xf$
there is an internal set $B \sbs \Xf$ and a symmetric internal
relation $S \sps E$ on $X$ such that $S[A] \sbs B$. It follows that
there is a~downward directed system $\RR$ of reflexive and symmetric
internal relations on $X$, and an upward directed system $\BB$ of
internal subsets of~$X$, both of admissible size, satisfying the
following conditions:
$$\xalignat{3}
&(\all\! R \in \RR)(\exs S \in \RR)(S \co S \sbs R)\,,
& &\text{and} & E &= \bigcap\ \RR\,,\\
(\all\! A &\in \BB)(\exs B \in \BB)(\exs S \in \RR)(S[A] \sbs B)\,,
& &\text{and} & \Xf &= \bigcup\ \BB\,.
\endxalignat$$
Then $\RR$ becomes a~base of a~uniformity $\UU_E$ on~$X$
(non-Hausdorff, unless $E = \Id_X$). Another base for this
uniformity (though not necessarily of admissible size) is formed by
all the internal relations $R$ on $X$ such that $E \sbs R$. A~set
$Y \sbs X$ is open in the induced topology if and only if for each
$y \in Y$ there is an internal set $A$ such that
$E[y] \sbs A \sbs Y$. In particular, $\Xf$ is an open subset of~$X$.
The closure of any set $Y \sbs X$ is $\bigcap_{R\in\RR} R[Y]$; for
internal $Y$ this is equal to $E[Y]$.

The observable traces
$$
R^\be = \bigl\{\bigl(x^\be,y^\be\bigr)\:
               (x,y) \in R \cap (\Xf \cx \Xf)\bigr\}
$$
of internal relations $R \in \UU_E$ (or just $R \in \RR$) form a
uniformity base on the observable trace $X^\be = \Xf/E$, inducing
a Hausdorff  completely regular topology on it.

We are particularly interested in representing Hausdorff locally
compact spaces as observable traces of IMG triplets $(X,E,\Xf)$
with hyperfinite ambient set $X$. To this end we introduce some
types of indices of internal sets $A \sbs X$ with respect to
reflexive and symmetric internal relations $S \sbs X \cx X$:
{\parindent 21pt
\item{(a)}
the {\it covering index\/} or {\it entropy\/} of $A$ with respect to
$S$, denoted by $\lfloor A:S\rfloor$, is the least $n\in\sbN$, such
that $A \sbs S[F]$ for some hyperfinite sets $F \sbs X$ with $n$
elements, or the symbol $\infty$ if there is no such~$n$;
\item{(b)}
the {\it inner covering index\/} of $A$ with respect to $S$, denoted
by $\lfloor A:S\rfloor_\ii$, is the least $n\in\sbN$, such that
$A \sbs S[F]$ for some hyperfinite sets $F \sbs A$ with $n$ elements,
or the symbol $\infty$ if there is no such~$n$;
\item{(c)}
the {\it independence index\/} or the {\it capacity\/} of $A$
with respect to $S$ denoted by $\lceil A:S\rceil$, is the biggest
$n\in\sbN$ such that there is an $n$-element set $F \sbs A$
satisfying $(x,y) \nin S$ for any distinct $x,y \in F$, or the
symbol $\infty$ if there is no biggest $n$ with that property.
\item{}}
\vskip-12pt
\flushpar
Then we have the following obvious inequalities (cf.~\cite{Ro}):
$$
\bigl\lceil A : (S \co S)\bigr\rceil \le \lfloor A : S\rfloor
\le \lfloor A : S\rfloor_\ii \le \lceil A : S\rceil\,.
$$

If $G$ is an internal group then, instead of the symmetric and
reflexive internal relation $S$ on $G$, we can take a symmetric
internal subset $S \sbs G$ containing the unit element $1 \in G$.
The obvious modification of the above definitions and the last
inequalities to this situation is left to the reader.

In the following Proposition the expression $[A : S]$ denotes
any of the indices $\lfloor A:S\rfloor$, $\lfloor A:S\rfloor_\ii$
or $\lceil A:S\rceil$. The above inequalities guarantee that the
formulations of (ii) with any particular choice for $[A : S]$ are
all equivalent.

\proclaim{1.2.1.~Proposition}
Let $(X,E,\Xf)$ be a bounded monadic space. Then the following
conditions are equivalent:
{\parindent 21pt
\item{\sl(i)}
all the internal subsets of \,$\Xf$ are relatively compact;
\item{\sl(ii)}
for any internal set $A \sbs \Xf$ and every symmetric internal
relation $S \sps E$ on $X$ the index $[A : S]$ is finite;
\item{\sl(iii)}
there is an external set $P \sbs \Xf$ of admissible size such that
$x \not\apr y$ for any distinct $x,y \in P$, and $\Xf \sbs S[P]$ for
every internal relation $S \sps E$;
\item{\sl(iv)}
for each $n \in \sbNi$ there is a hyperfinite set $H \sbs X$ with
at most $n$ elements such that $\Xf \sbs E[H]$;
\item{\sl(v)}
for every infinite hyperfinite set $H \sbs \Xf$ there are at least
two distinct elements $x,y \in H$ such that $x \apr y$.
\item{}}
\vskip-12pt
\endproclaim

\demo{Proof}
(i)$\iff$(ii)
is plain, especially for the index $\lfloor A : S\rfloor_\ii$.

(ii)$\imp$(iii)
Let $\RR$ and $\BB$ be systems of admissible size consisting of
internal symmetric relations $R \sbs X \cx X$ and internal sets
$B \sbs X$, respectively, such that $E = \bigcap \RR$ and
$\Xf = \bigcup \BB$. Condition (ii) with the index
$\lfloor A : S\rfloor_\ii$ implies that there is a system of finite
sets $\{F_{RB}\: R \in \RR,\,B \in \BB\}$ such that
$$
F_{RB} \sbs B \sbs R[F_{RB}]
$$
for all $R \in \RR$, $B \in \BB$. Then the set
$Q = \bigcup_{R,B} F_{RB} \sbs \Xf$ is of admissible size and,
obviously, $\Xf \sbs S[Q]$ for any internal relation $S \sps E$ on $X$.
Now, it suffices to take for $P$ any subset of \,$Q$ containing exactly
one point from every $E$-monad intersecting $Q$.

(iii)$\imp$(iv)
Let $P$ be the set guaranteed by (ii) and $n$ be an arbitrary
infinite hypernatural number. Then, for each finite set $F \sbs P$,
we have $|F| \le n$. By {\it saturation}, there is a hyperfinite set
$H$ such that $P \sbs H \sbs X$ and $|H| \le n$. It can be easily
verified that $\Xf \sbs E[H]$.

(iv)$\imp$(v)
Let $H \sbs \Xf$ be any infinite hyperfinite set and $n$ be an
infinite hypernatural number such that $n < |H|$. By (iv) there is a
hyperfinite set $J \sbs X$ such that $|J| \le n$ and $\Xf \sbs E[J]$.
In particular, $H \sbs E[J]$. As $|J| < |H|$, there is at least one
element $z \in J$ such that the set $H \cap E[z]$ contains at least
two distinct elements $x$, $y$. Then $x \apr y$.

(v)$\imp$(ii)
Assume that the independence index $\lceil A : S\rceil$ $A \sbs \Xf$
is not finite for some internal set $A \sbs \Xf$ and an internal
symmetric relation $S \sps E$. Then, by {\it saturation}, there is an
infinite hyperfinite set $H \sbs A$ such that $x \not\apr y$ for any
distinct $x,y \in H$, contradicting (v).
\enddemo

\remark{Remark}
Let us turn the reader's attention to the following two details.
First, in order to prove the implication (iii)$\imp$(iv) the set
$P$ just has to be of admissible size and satisfy the condition
$\Xf \sbs S[P]$ for every internal relation $S \sps E$; it is not
necessary that it consists of pairwise discernible points. Second,
(v) could be strengthened to the following Ramsey type statement:
{\it every infinite hyperfinite set $H \sbs \Xf$ contains an
infinite hyperfinite subset $H_0$ such that $x \apr y$ for any
$x,y \in H_0$}. Then the implication (iv)$\imp$(v) could still be
proved by picking an $n \in \sbNi$ such that, say, $n^2 \le |H|$.
\endremark
\medskip

The last condition (v) suggests to call the IMG triplets
satisfying any (hence all) of the above conditions
{\it condensing\/} (cf\. \cite{ZZ}). Obviously, the observable
trace of any condensing IMG space is locally compact, and the
compact subsets of \,$X^\be$ are exactly the observable traces
$A^\be$ of internal subsets $A \sbs \Xf$. However, it should be kept
in mind that this is a considerably stronger condition than just
the local compactness of \,$X^\be$. Nevertheless, we still have the
following representation theorem.

\proclaim{1.2.2.~Proposition}
Let $\Xb$ be a Hausdorff locally compact topological space. Then,
in every sufficiently saturated nonstandard universe, there is a
condensing IMG triplet $(X,E,\Xf)$ such that $\Xb$ is homeomorphic
to the observable trace $X^\be$. If desirable, one can additionally
arrange that the ambient space $X$ be hyperfinite.
\endproclaim

\demo{Proof}
Let $\UU$ be some uniformity inducing the topology of \,$\Xb$
and $\kappa$ be an uncountable cardinal bigger than the minimal
cardinality of some base of \,$\UU$ as well as of some open cover
of \,$\Xb$ by relatively compact sets. Let us embed $\Xb$ into its
nonstandard extension $\sXb$ in some $\kappa$-saturated nonstandard
universe. Put $E_\UU = \bigcap_{\Ub \in \UU}\sUb$. It can be
easily verified that $\bigl(\sXb, E_\UU, \Ns(\sXb)\bigr)$
is a condensing IMG space whose observable trace (nonstandard
hull) $\Ns(\sXb)/E_\UU$ is homeomorphic to $\Xb$.

Let $n$ be an arbitrary infinite hypernatural number and
$X = H \sbs \sXb$ be the $n$-element hyperfinite set guaranteed by
(iv) of 1.2.1. Now, it suffices to put $\Xf = \Ns(\sXb) \cap H$,
\,$E = E_\UU \cap (H \cx H)$, and we get another condensing
IMG triplet  $(X,E,\Xf)$ with the observable trace
$\Xf/E \iso \Ns(\sXb)/E_\UU \iso \Xb$ and hyperfinite ambient
space $X$.
\enddemo

The crucial property of the internal inclusion mapping
$\Id_X: X \to \sXb$ (under the identification ${}^*\xb = \xb$ for
$\xb \in \Xb$) is namely that for each $\xb \in \Xb$ there is some
$x \in X$ such that $x \apr \xb$. More generally, given any hyperfinite
set $X$, an internal mapping $\eta\: X \to \sXb$ is called a
{\it hyperfinite infinitesimal approximation}, briefly
{\it HFI approximation}, of the Hausdorff uniform space $\Xb$ if
$$
(\all \xb \in \Xb)(\exs x \in X)(\eta(x) \apr \xb)
$$
It is called an {\it injective HFI approximation\/} if, additionally,
$\eta$ is an injective mapping.

The standard counterpart of this notion can be formulated in terms of
approximating systems by finite sets.
Let $\Xb$ be any set, $\Kb \sbs \Xb$ be nonempty set and
$\Ub \sbs \Xb \cx \Xb$ be a reflexive relation. A~mapping
$\eta\: X \to \Xb$ is called a
{\it finite $(\Kb,\Ub)$ approximation\/} of \,$\Xb$ if \,$X$
is a finite set and
$$
(\all \xb \in \Kb)(\exs x \in X)\bigl((\eta(x),\xb) \in \Ub\bigr)\,.
$$
The meaning of the term {\it injective $(\Kb,\Ub)$ approximation\/} is
obvious.

The term {\it poset\/} $(I,\le)$ is used as a shorthand for a partially
ordered set, i.e. a set $I$ equipped with a reflexive, antisymmetric and
transitive binary relation $\le$; if the relation $\le$ is upward directed,
then $(I,\le)$ is referred to as a {\it directed poset}.

Let $(\Xb,\UU)$ be a Hausdorff uniform space and $(I,\le)$ be
a directed poset. Then a system of mappings
$(\eta_i\:X_i \to \Xb)_{i\in I}$ is called an {\it approximating
system\/} of the space $(\Xb,\UU)$ provided each $X_i$ is a finite
set, and for any $\Ub \in \UU$ and any compact set $\Kb \sbs \Xb$
there is an $i \in I$ such that $\eta_j\:X_j \to \Xb$ is a $(\Kb,\Ub)$
approximation of \,$\Xb$ for each $j \in I$, $j \ge i$.

Now, it is almost obvious that every Hausdorff locally compact
uniform space $(\Xb,\UU)$ has some approximating system
$(\eta_i\:X_i \to \Xb)_{i\in I}$ such that each $X_i$ is a finite
subset of \,$\Xb$ and $\eta_i\: X_i \to \Xb$ is the inclusion mapping.
Assuming $\kappa$-saturation for some sufficiently big $\kappa$,
we can take a ${}^*$compact set $\Kb_0 \sbs \sXb$ and a
$\Ub_0 \in {}^*\UU$, such that $\sKb \sbs \Kb_0$ and
$\Ub_0 \sbs \sUb$ for all compact $\Kb \sbs \Xb$ and
$\Ub \in \UU$. Then there is an $i \in {}^{*\!}I$ such that the
hyperfinite set $X_i \sbs \sXb$ (together with the inclusion mapping
$\eta_i\:X_i \to \sXb$) is an internal $(\Kb_0,\Ub_0)$ approximation,
hence an HFI approximation, of \,$\Xb$. Putting $X = X_i$,
$E = E_\UU \cap (X \cx X)$, and $\Xf = \Ns(\sXb) \cap X$, we get
a condensing IMG triplet $(X,E,\Xf)$ with hyperfinite ambient
space $X$ and observable trace $X^\be = \Xf/E \iso \Xb$. This gives
another proof of Proposition~1.2.2.

Given two IMG spaces $(X,E,\Xf)$, $(Y,F,\Yf)$, an internal
mapping $f\:D \to Y$ is called a {\it triplet morphism\/} if
$\Xf \sbs D \sbs X$, it {\it preserves finiteness}, i.e.,
$f(x) \in \Yf$ for any $x \in \Xf$, and it is {\it $S$-continuous\/}
on $\Xf$, i.e.,
$$
x \apr y \imp f(x) \apr f(y)\,,
$$
for any $x,y \in \Xf$. In such a case we write
$f\:(X,E,\Xf) \to (Y,F,\Yf)$. Note that every triplet morphism $f$
with domain $D$ can be formally extended to an everywhere defined
triplet morphism $\wtl f\: X \to Y$ in an arbitrary way.

Every triplet morphism $f\:(X,E,\Xf) \to (Y,F,\Yf)$ induces a
continuous mapping $f^\be\:X^\be \to Y^\be$, called the {\it
observable trace\/} of $f$, (correctly) defined by
$$
f^\be\bigl(x^\be\bigr) = f(x)^\be\,,
$$
for $x \in \Xf$. Two triplet morphisms
$f,g\:(X,E,\Xf) \to (Y,F,\Yf)$ are called {\it equivalent\/} if
they have the same observable trace $f^\be = g^\be$. However, we
can put
$$
f \apr_{\Xf} g \iff (\all x \in \Xf)(f(x) \apr g(x))
$$
for any internal functions $f\: D_1 \to Y$, $g\:D_2 \to Y$ such that
$\Xf \sbs D_1 \cap D_2$. This {\it relation of infinitesimal nearness
on finite elements\/} is a monadic equivalence on any set $Y^D$ of
all internal functions $D \to Y$, where $\Xf \sbs D$. For triplet
morphisms $f,g\:(X,E,\Xf) \to (Y,F,\Yf)$ we have
$$
f \apr_{\Xf} g \iff f^\be = g^\be\,
$$
indicating a fundamental role of this indiscernibility equivalence.
Let us remark, without being precise, that for a condensing IMG
space $(X,E,\Xf)$ the equivalence $\apr_{\Xf}$ corresponds to the
compact-open topology on the space $\C\bigl(X^\be,Y^\be\bigr)$ of
all continuous functions $X^\be \to Y^\be$.

Another important feature of condensing IMG spaces is the lifting
property. Given a mapping $\fb\: X^\be \to Y^\be$ between the
observable traces, an internal mapping $f\: D \to Y$ is called an
{\it $S$-continuous lifting\/} of \,$\fb$ if \,$\Xf \sbs D \sbs X$
and
$$
\fb\bigl(x^\be\bigr) = f(x)^\be\,,
$$
for each $x \in \Xf$. Notice that an internal mapping $f\: D \to Y$,
satisfying the last equality necessarily is a triplet morphism,
hence $S$-continuous on $\Xf$. In other words, a triplet morphism
$f\:(X,E,\Xf) \to (Y,F,\Yf)$ is a lifting of \,$\fb$ if and only if
\,$\fb = f^\be$ is the observable trace of $f$. Then $\fb$
necessarily is continuous, as well. Thus only continuous mappings
between observable traces of \,IMG triplets have liftings
$S$-continuous on $\Xf$. The point is that for a {\it condensing\/}
IMG triplet $(X,E,\Xf)$ this necessary continuity condition is
already sufficient for the existence of liftings.

\proclaim{1.2.3.~Proposition}
Let $(X,E,\Xf)$, $(Y,F,\Yf)$ be two IMG spaces. If \,$(X,E,\Xf)$
is con\-dens\-ing, then a mapping $\fb\: X^\be \to Y^\be$ has
an internal lifting $S$-continuous on $\Xf$ if and only if \,$\fb$
is continuous.
\endproclaim

\demo{Sketch of proof}
Let's focus on the nontrivial implication, only. Assume that $\fb$
is continuous and denote by
$$
\fb^\kr = \bigl\{(x,y) \in X \cx Y\:
   \fb\bigl(x^\be\bigr) = y^\be\bigr\}
$$
its pretrace considered as a set $\fb^\kr \sbs \Xf \cx \Yf$ in the
IMG space $(X \cx Y, E \cx F, \Xf \cx \Yf)$. It suffices to show
that, for every internal set $D \sbs \Xf$, the restriction
$\fb^\kr\rst D$ is monadic. Then, by the Internal Choice
Lemma~1.1.1, there is an internal choice function $f$ from $\fb^\kr$
on its domain $\Xf$. Obviously, this $f$ is an $S$-continuous
lifting of $\fb$.

Preliminarily we can only assure that that there is some (not
necessarily internal) function $\vfi\:\Xf \to Y$, such that
$\fb\bigl(x^\be\bigr) = \vfi(x)^\be$ for each $x \in \Xf$.

Let us fix some nonempty internal set $D \sbs \Xf$. Let further
$\RR$, $\SS$ be some systems of admissible size consisting of
symmetric internal relations on $X$, $Y$ respectively, such that
$E = \bigcap \RR$, $F =\bigcap \SS$, and $P \sbs \Xf$ be a set of
admissible size, such that $R[P] = \Xf$ for every $R \in \RR$,
whose existence is guaranteed in Proposition~1.2.1(iii). Denote
by $\II_D$ the set of all ordered triples $(S,R,A)$ such that
$S \in \SS$, $R \in \RR$, $A$ is a finite subset of $P$ subject
to $D \sbs R[A]$, and
$$
\bigl(x^\be,y^\be\bigr) \in R^\be \imp
\bigl(\fb\bigl(x^\be\bigr),\fb\bigl(y^\be\bigr)\bigr) \in S^\be
   $$
for all $x,y \in D$. Obviously, $\II_D$ is of admissible size.
Moreover, $D^\be \sbs X^\be$ is compact, hence $\fb$ is uniformly
continuous on $D^\be$. Therefore, for any $S \in \SS$ there is an
$R \in \RR$ and a finite $A \sbs P$ such that $(S,R,A) \in \II_D$.
For any triple $i = (S,R,A) \in \II_D$ we
put
$$
\vFi_i = \bigcup_{a \in A}
   \bigl(D \cap R[a]\bigr) \cx S\bigl[\vfi(a)\bigr]\,.
$$
As $A$ is finite, each $\vFi_i$ is an internal relation. Using
the continuity of $\fb$, the equality
$$
\fb^\kr\rst D = \bigcap_{i\in\II_D} \vFi_i
$$
can be checked in a routine way, resembling the standard proof of
the closed graph theorem for a continuous function into a Hausdorff
space.
\enddemo

\proclaim{1.2.4.~Corollary}
Let $(X,E,\Xf)$, $(Y,F,\Yf)$ be condensing bounded monadic spaces
with homeo\-morphic observable traces \,$X^\be \iso Y^\be$. Then there
exist triplet morphisms $f\:(X,E,\Xf) \to (Y,F,\Yf)$,
$g\:(Y,F,\Yf) \to (X,E,\Xf)$ such that
$$
g(f(x)) \apr x \qquad\text{and}\qquad f(g(y)) \apr y\,,
$$
for all $x \in \Xf$, $y \in \Yf$.
\endproclaim

Naturally, a triplet morphism $f\:(X,E,\Xf) \to (Y,F,\Yf)$ to which
there is a triplet morphism $g\:(Y,F,\Yf) \to (X,E,\Xf)$ satisfying
the above condition will be called a {\it triplet isomorphism}. Now,
the last Corollary can be restated as follows: {\it Condensing\/}
IMG triplets are isomorphic if and only if they have homeomorphic
observable traces.

\proclaim{1.2.5.~Corollary}
Let \,$\Xb$ be a Hausdorff locally compact uniform space and
$(X,E,\Xf)$ be a condensing IMG space with hyperfinite ambient
set and the observable trace $X^\be$ homeomorphic to $\Xb$.
Then there is an HFI approximation $\eta\:X \to \sXb$ which
is a triplet isomorphism
$\eta\:(X,E,\Xf) \to \bigl(\sXb, E_\UU,\Ns(\sXb)\bigr)$.
\endproclaim

\specialhead
1.3.~Functional spaces 1: Continuous functions
\endspecialhead
\flushpar
In this and the next following section $\Xb$ is a Hausdorff locally
compact topological space, whose topology is induced by a uniformity
$\UU$, represented as the observable trace $\Xb \iso X^\be$ of a
condensing IMG triplet $(X,E,\Xf)$ with a {\it hyperfinite\/}
ambient set $X$ by means of a (not necessarily injective)
HFI  approximation $\eta\: X \to \sXb$. Then
$\Xf = \eta^{-1}\bigl[\Ns({}^*\Xb)\bigr]$ and we can assume,
without  loss of generality, that
$x \apr y \iff \eta(x) \apr \eta(y)$ for all $x$, $y$ in $X$
and not just in $\Xf$. Identifying the observable trace
$X^\be = \Xf/E$ with $\Xb \iso \Ns({}^*\Xb)/\!\!\apr$ via the
homeomorphism $\eta^\be$, we regard each point
$\eta^\be\bigl(x^\be\bigr) = {}^\co\eta(x) \in \Xb$ as
observable trace $x^\be$ of \,$x \in \Xf$.

The hyperfiniteness of \,$X$ enables to represent various Banach
spaces of functions $\Xb \to \bC$ by means of the {\it hyperfinite
dimensional\/} linear space $\sbC^X$ of all internal functions
$X \to \sbC$. We will systematically employ the advantage of such
an approach.

Let's start with the observation that the nonstandard extensions
$\sbR$, $\sbC$ of real and complex numbers, respectively, give
rise to the condensing IMG triplets $(\sbR,\IR,\FR)$ and
$(\sbC,\IC,\FC)$. Here $\FR \sbs \sbR$ or $\FC \sbs \sbC$ denote
the galactic subrings of finite (bounded) hyperreal or hypercomplex
numbers, and the monadic ideals  $\IR \sbs \FR$ or $\IC \sbs \FC$ of
infinitesimal hyperreal or hypercomplex numbers are used instead of
the equivalence relations $\apr$ of infinitesimal nearness on $\sbR$
or $\sbC$ in the notation of the triplets. Then, obviously,
$\bR \iso \FR/\IR$ and $\bC \iso \FC/\IC$ as topological fields.
For $x \in \sbR$, $x \ge 0$, we write $x < \infty$ instead of
$x \in \FR$, and $x \sim \infty$ instead of $x \in \sbR \sms \FR$.

From now on we will focus on spaces of complex functions, leaving
the reader the formulation of the real version.

The hyperfinite dimensional (HFD) linear space $\sbC^X$ admits
several internal norms. For any internal norm $\Nn$ on $\sbC^X$
we denote by
$$\align
\InC^X &= \bigl\{f \in \sbC^X\: \Nn(f) \apr 0\bigr\}\,, \\
\FnC^X &= \bigl\{f \in \sbC^X\: \Nn(f) < \infty\bigr\}\,,
\endalign$$
the $\FC$-linear subspaces (more precisely, $\FC$-submodules)
of \,$\sbC^X$, consisting of functions which are infinitesimal
or finite, respectively, with respect to the norm $\Nn$.
Then the arising IMG triplet $\bigl(\sbC^X,\InC^X,\FnC^X\bigr)$
(with $\InC$ standing in place of the indiscernibility equivalence
relation $f \apr_{_{\nn}} g \iff f - g \in \InC^X$) which (at least
for a ``reasonable'' norm $\Nn$) is condensing if and only if $X$
is finite. Its observable trace (nonstandard hull)
$$
\bigl(\sbC^X\bigr)_{\nn}^\be = \FnC^X\!/\,\InC^X
$$
becomes a (standard) Banach space under the norm $\Nn^\be$, given
by
$$
\Nn^\be\bigl(f_{_{\nn}}^\be\bigr) = {}^\co\Nn(f) = \st\Nn(f)\,,
$$
where $f_{_{\nn}}^\be \in \bigl(\sbC^X\bigr)_{\nn}^\be$ denotes
the observable trace of the function $f \in \FnC^X$ with respect
to the norm $\Nn$. Typically, $\bigl(\sbC^X\bigr)_{\nn}^\be$ is
nonseparable unless $X$ is finite.

An internal function $f\:D \to \sbC$, such that $\Xf \sbs D \sbs X$,
is a triplet morphism $(X,E,\Xf) \to (\sbC,\IC,\FC)$ if and only if
$f$ is $S$-continuous on $\Xf$ and $f[\Xf] \sbs \FC$. Its observable
trace is the function $f^\be\: X^\be \iso \Xb \to \bC$ given by
$$
f^\be\bigl(x^\be\bigr) = {}^{\co\!}f(x) = \st f(x)
$$
for $x \in \Xf$. However, unless $\Xf = X$, the monadic equivalence
relation $\apr_{\Xf}$ on $\sbC^X$, corresponding to the compact-open
topology on the space $\C(\Xb)$ of all continuous functions
$\Xb \to \bC$, is not of the form $\apr_{_{\nn}}$ for any internal
norm $\Nn$ on $\sbC^X$.

When dealing with $S$-continuous functions, the {\it maximum norm}
$$
\|f\|_{_\infty} = \max_{x\in X} |f(x)|\,,
$$
where $f \in \sbC^X$, becomes rather important. Denoting by
$$\align
\IiC^X &= \bigl\{f \in \sbC^X\:
          \|f\|_{_\infty} \apr 0\bigr\}\,, \\
\FiC^X &= \bigl\{f \in \sbC^X\:
          \|f\|_{_\infty} < \infty\bigr\}\,,
\endalign$$
the $\FC$-linear subspaces of \,$\sbC^X$, consisting of internal
functions which are infinitesimal or finite, respectively,
with respect to the max-norm we get the IMG triplet
$\bigl(\sbC^X,\IiC^X,\FiC^X\bigr)$.

We are particularly interested in the Banach spaces $\Cb(\Xb)$,
$\Cbu(\Xb)$, and $\Cbo(\Xb)$ of all bounded continuous, bounded
uniformly continuous, and continuous vanishing at infinity functions
$\Xb \to \bC$, respectively, and also in the (non-Banach) normed
linear space $\Cbc(\Xb)$ of all continuous functions $\Xb \to \bC$
with compact support, all with the supremum norm denoted by
$\|\!\cd\!\|_{_\infty}$, as well. Let us denote by
$$\align
\CCb(X,E,\Xf) &= \bigl\{f \in \FiC^X\:
(\all x,y \in \Xf)(x \apr y \imp f(x) \apr f(y))\bigr\}\,,\\
\CCbu(X,E) &= \bigl\{f \in \FiC^X\:
(\all x,y \in X)(x \apr y \imp f(x) \apr f(y))\bigr\}\,,\\
\CCo(X,E,\Xf) &= \bigl\{f \in \CCb(X,E,\Xf)\:
(\all x \in X \sms \Xf)(f(x) \apr 0)\bigr\}\,,\\
\CCc(X,E,\Xf) &= \bigl\{f \in \CCb(X,E,\Xf)\:
(\all x \in X \sms \Xf)(f(x) = 0)\bigr\}
\endalign$$
their intended nonstandard counterparts. Each $f \in \CCb(X,E,\Xf)$
is an everywhere defined triplet morphism
$(X,E,\Xf) \to (\sbC,\IC,\FC)$; (however there can be also such
triplet morphisms not belonging to $\CCb(X,E,\Xf)$). Moreover,
we have the obvious inclusions of \,$\FC$-linear subspaces of
\,$\sbC^X$:
$$
\CCc(X,E,\Xf) + \IiC^X \sbs \CCo(X,E,\Xf) \sbs \CCbu(X,E)
\sbs \CCb(X,E,\Xf)\,.
$$

\proclaim{1.3.1.~Proposition}
Let $\fb\:\Xb \to \bC$ be any function. Then
{\parindent 21pt
\item{\sl(a)}\,
$\fb \in \Cb(\Xb)$ if and only if \,$\fb$ has a lifting
$f \in \CCb(X,E,\Xf)$;
\item{\sl(b)}\,
$\fb \in \Cbu(\Xb)$ if and only if \,$\fb$ has a lifting
$f \in \CCbu(X,E)$;
\item{\sl(c)}\,
$\fb \in \Cbo(\Xb)$ if and only if \,$\fb$ has a lifting
$f \in \CCo(X,E,\Xf)$;
\item{\sl(d)}\,
$\fb \in \Cbc(\Xb)$ if and only if \,$\fb$ has a lifting
$f \in \CCc(X,E,\Xf)$.
\item{}}
\vskip-12pt
\flushpar
Conversely, every internal function $f \in \CCb(X,E,\Xf)$ is
lifting of a function $\fb \in \Cb(\Xb)$, every internal function
$f \in \CCbu(X,E)$ is lifting of a function $\fb \in \Cbu(\Xb)$,
every internal function $f \in \CCo(X,E,\Xf)$ is lifting of a function
$\fb \in \Cbo(\Xb)$, and every internal function $f \in \CCc(X,E,\Xf)$
is lifting of a function $\fb \in \Cbc(\Xb)$.
\endproclaim

\demo{Proof}
Let us start with the observation that, under the identification
$x^\be = {}^\co\eta(x)$ for $x \in \Xf$, we have
$({}^{*\!}\fb \co \eta)^\be = \fb$ for any continuous function
$\fb\: \Xb \to \bC$. Next, we leave the reader to verify the
following easy facts:
$$\align
\fb \in \Cb(\Xb) &\imp {}^{*\!}\fb \co \eta \in \CCb(X,E,\Xf)\,,\\
\fb \in \Cbu(\Xb) &\imp {}^{*\!}\fb \co \eta \in \CCbu(X,E)\,,\\
\fb \in \Cbo(\Xb) &\imp {}^{*\!}\fb \co \eta \in \CCo(X,E,\Xf)\,, \\
\fb \in \Cbc(\Xb) &\imp {}^{*\!}\fb \co \eta \in \CCc(X,E,\Xf)\,, \\
f \in \CCb(X,E,\Xf) &\imp f^\be \in \Cb(\Xb)\,,\\
f \in \CCbu(X,E) &\imp f^\be \in \Cbu(\Xb)\,,\\
f \in \CCo(X,E,\Xf) &\imp f^\be \in \Cbo(\Xb)\,,\\
f \in \CCc(X,E,\Xf) &\imp f^\be \in \Cbc(\Xb)\,,
\endalign$$
for  $\fb \in \bC^{\Xb}$, $f \in \sbC^X$. Now, (a), (b), (c) and (d)
follow from the first quadruple of implications (the additional
property $E = \{(x,y) \in X \cx X\: \eta(x) \apr \eta(y)\}$ is
needed for the second implication.
The ``conversely part'' follows from the second quadruple of
implications.
\enddemo

However, in general the observable traces $f^\be$ and
$f^\be_{\!_\infty}$ of \,$S$-continuous functions should not be
confused. While for $f, g \in \CCo(X,E,\Xf)$ (and the more for
$f,g \in \CCc(X,E,\Xf)$) we have the equivalence
$$
f^\be = g^\be \iff \|f - g\|_{_\infty} \apr 0\,,
$$
for $f,g \in \CCbu(X,E)$, and the more for $f,g \in \CCb(X,E,\Xf)$,
we have just the implication
$$
\|f - g\|_{_\infty} \apr 0 \imp f^\be = g^\be\,,
$$
and, unless $\Xf$ is dense in $X$, there are functions
$f,g \in \CCbu(X,E)$ such that $f^\be = g^\be$ but
$\|f-g\|_{_\infty} \not\apr 0$. Summing up we have

\proclaim{1.3.2~Corollary}
The observable trace map $f \mto f^\be$ induces
{\parindent 21pt
\item{\sl(a)}
a bounded surjective linear mapping of the
subspace $\CCb(X,E,\Xf)/\,\IiC^X$ of the nonstandard hull
\,$\FiC^X\!/\,\IiC^X$ onto the Banach space $\Cb(\Xb)$;
\item{\sl(b)}
a bounded surjective linear mapping of the
subspace $\CCbu(X,E)/\,\IiC^X$ of the nonstandard hull
\,$\FiC^X\!/\,\IiC^X$ onto the Banach space $\Cbu(\Xb)$;
\item{\sl(c)}
a Banach space isomorphism of the subspace $\CCo(X,E,\Xf)/\,\IiC^X$
of the nonstandard hull \,$\FiC^X\!/\,\IiC^X$ onto the Banach space
$\Cbo(\Xb)$;
\item{\sl(d)}
a normed space isomorphism of the subspace $\CCc(X,E,\Xf)/\,\IiC^X$
of the nonstandard hull \,$\FiC^X\!/\,\IiC^X$ onto the normed space
$\Cbc(\Xb)$.
\item{}}
\vskip-12pt
\endproclaim

\specialhead
1.4.~Functional spaces 2: Loeb measures and Lebesgue spaces
\endspecialhead
\flushpar
Let us recall that $\Xb$ denotes a Hausdorff locally compact space,
represented as the observable trace $\Xb \iso X^\be$ of a condensing
IMG triplet $(X,E,\Xf)$ with a {\it hyperfinite\/} ambient set $X$
by means of an HFI approximation $\eta\: X \to \sXb$.

The natural way of getting functions $\fb\:\Xb \to \bC$ as
observable traces $\fb = f^\be$ of internal functions
$f\:X \to \sbC$ works just under the assumption that $f$ is
finite and $S$-continuous on $\Xf$. In this way, however, just
{\it continuous} functions $\fb\:\Xb \to \bC$ can be obtained.
Therefore, by far not all internal functions $f\:X \to \sbC$
represent classical functions $\Xb \to \bC$. On the other hand,
they can be used to represent various objects of different nature:
measures, distributions, etc. The class of that way representable
functions $\fb$ on $\Xb = X^\be$ can be extended to encompass the
Lebesgue spaces $\Lbp(\Xb)$ by relaxing the equality
$\fb\bigl(x^\be\bigr) = {}^{\circ\!}f(x)$ on $\Xf$ to the equality
{\it almost everywhere\/} on $\Xf$ with respect to some measure.
Strictly speaking, the elements of $\Lbp(\Xb)$ themselves are not
genuine functions but certain equivalence classes of functions.
We are going to make this point more precise.

Let $d\:X \to \sbR$ be an internal function such that $d(x) \ge 0$
for each $x \in X$. Intuitively, $d(x)$ is viewed as the ``weight''
of the point $x$. Then $d$ induces the internal (hyper)finitely
additive measure $\nu_d\:\PP(X) \to \sbR$ on the internal Boolean
algebra $\PP(X)$ of all internal subsets of \,$X$, given by
$\nu_d(A) = \sum_{a \in A} d(a)$ for internal $A \sbs X$. Putting
$$
({}^{\co}\nu_d)(A) = \cases
{}^{\co}(\nu_d(A)) &\text{if \,$\nu_d(A) < \infty$,} \\
\infty &\text{if \,$\nu_d(A) \sim \infty$,}
\endcases
$$
we get a finitely additive (non-negative) measure
${}^{\co}\nu_d\:\PP(X) \to \bR \cup \{\infty\}$ which has a unique
extension to a $\sig$-additive measure
$\lam_d\:\wPP(X) \to \bR \cup \{\infty\}$ defined on the
$\sig$-algebra $\wPP(X)$ of all subsets of \,$X$ generated by all
monadic (or, equivalently, by all galactic) subsets of \,$X$. (If
$\kappa = \aleph_1$, then $\wPP(X)$ is simply the $\sig$-algebra
generated by all internal subsets of \,$X$.) Then $\lam_d$ is the
{\it Loeb measure\/} induced by the internal function $d$ (cf\.
\cite{Ct}, \cite{Lb1}, \cite{LR}).

Assume that the function $d$ satisfies additionally the condition
\,$\nu_d(A) \in \FR$ for each internal set $A \sbs \Xf$.
Then the set $\Xf$ belongs to the algebra $\wPP(X)$ and the system
$$
\wPP(\Xf) = \bigl\{A \in \wPP(X)\: A \sbs \Xf\bigr\}
$$
is again a $\sig$-algebra of subsets of \,$\Xf$. Moreover, for each
Borel set $\Yb \sbs \Xb$ its pretrace $\Yb^\kr$ belongs to
$\wPP(\Xf)$, hence the observable trace map
\,${}^\be\:\Xf \to \Xb = X^\be$ is measurable. Pushing down the
Loeb measure along this map, i.e., putting
$$
\mbd(\Yb) = \lam_d\bigl(\Yb^\kr\bigr)
$$
for Borel $\Yb \sbs \Xb$, we get a regular Borel measure $\mbd$ on
$\Xf$ (which, if desirable, can be extended to a complete measure
by the Carath\'eodory construction). Equally important is the converse.

\proclaim{1.4.1.~Proposition}
Every nonnegative regular Borel measure $\mb$ on $\Xb$ has the form
$\mb = \mbd$ for some internal function $d\:X \to \sbR$, such that
\,$d(x) \ge 0$ for each \hbox{$x \in X$,} and \,$\nu_d(A) \in \FR$
for each internal set \,$A \sbs \Xf$. Additionally, if \,$\mb$ is
not the identically zero measure, then it can be arranged that
\,$d(x) > 0$ for each $x \in X$.
\endproclaim

\demo{Sketch of proof}
There is a symmetric entourage $\Ub_0 \in {}^*\UU$, ${}^*$open in
$\sXb \cx \sXb$, and a ${}^*$compact set $\Kb_0 \sbs \sXb$ such
that $\Ub_0 \sbs E_{\UU}$, $\Ns(\sXb) \sbs \Kb_0$, and $\eta$ is
a $(\Ub_0,\Kb_0)$ approximation of \,$\sXb$. Since $\eta$ has
hyperfinite range, there is an internal mapping $\sig\:\sXb \to X$
such that $(\eta \co \sig \co \eta)(x) = \eta(x)$ and
$\bigl((\eta \co \sig)(\xb),\xb\bigr) \in \Ub_0$ for $x \in X$,
$\xb \in \Kb_0$, and each of the sets
$\{\xb \in \Kb_0\: (\eta\co\sig)(\xb)=\eta(x)\}$
is ${}^*$Borel in $\sXb$. We put
$$
d(x) = \cases
\dsize\frac{{}^*\mb\bigl(\{\xb \in \Kb_0\:
(\eta\co\sig)(\xb)=\eta(x)\}\bigr)}
{\bigl|\{y \in X\: \eta(y)=\eta(x)\}\bigr|} \,,
     & \text{if \,$\eta(x) \in \Ub_0[\Kb_0]$,} \\
0\,, &\text{if \,$\eta(x) \nin \Ub_0[\Kb_0]$,}
\endcases
$$
for $x \in X$.%
\footnote{In fact, the values of $d(x)$ for
$\eta(x) \nin \Ub_0[\Kb_0]$ can be chosen arbitrarily without
affecting the resulting measure $\mb_d$.}
The verification that the internal mapping $d$ has
all the required properties is left to the reader. Replacing each
value $d(x)$ by $d(x) + \eps$, where $\eps$ is a positive
infinitesimal such that $\eps\lv X\rv \apr 0$, we can satisfy also
the additional requirement.
\enddemo

Similarly, every internal function $g\:X \to \sbC$, such that
$\sum_{x \in X} |g(x)|$ is finite, gives rise to the finite
complex Loeb measure $\lam_g\:\wPP(X) \to \bC$, such that
$$
\lam_g(A) = {\vphantom{\Bigl)}}^\co
             \biggl(\,\sum_{x \in A} g(x)\biggr)
$$
for internal $A \sbs X$, and to a complex regular Borel measure
$\tth_g$ on $\Xb$, given by
$$
\tth_g(\Yb) = \lam_g\bigl(\Yb^\kr\bigr)
$$
for Borel $\Yb \sbs \Xb$. Moreover, $\tth_g$ has finite variation
$$
\|\tth_g\| \le {\vphantom{\Bigl)}}^\co
             \biggl(\,\sum_{x \in X} |g(x)|\biggr)
= \lam_{|g|}(X)\,.
$$

\proclaim{1.4.2.~Proposition}
Every complex regular Borel measure $\mmu$ on $\Xb$ with finite
variation has the form $\mmu = \tth_g$ for some internal function
$g\:X \to \sbC$, such that
$$
\|\mmu\| = {\vphantom{\Bigl)}}^\co\biggl(\,\sum_{x \in X} |g(x)|\biggr)\,
\quad\text{and}\qquad
\sum_{x \in Z} |g(x)| \apr 0
$$
for each internal set $Z \sbs X \sms \Xf$.
\endproclaim

\demo{Sketch of proof}
Essentially the same construction as used to obtain the function $d$
from the nonnegative measure $\mb$ in Proposition~1.4.1 works for
every complex regular Borel measure $\mmu$ with finite variation to give
the function $g$. One just has to take care that the  hyperfinite ${}^*$Borel
partition of \,$\Kb_0$ formed by the sets
$\{\xb \in \Kb_0\: (\eta\co\sig)(\xb)=\eta(x)\}$, where $x$ runs over
some maximal internal set $X_0 \sbs \eta^{-1}[\Ub_0[\Kb_0]]$ such that
the restriction $\eta\rst X_0$ is injective, satisfies additionally
$$
\sum_{x \in X_0}
\bigl|{}^{*\!}\mmu\bigl(\{\xb \in \Kb_0\:
(\eta\co\sig)(\xb)=\eta(x)\}\bigr)\bigr|
\le \|\mmu\|\,,
$$
which can be achieved by the virtue of {\it saturation}. A complete proof
of Proposition~1.4.2, following a slightly different chain of arguments,
can be found in \cite{ZZ, Section~4}.
\enddemo

In the rest of this section $d\:X \to \sbR$ denotes a fixed
internal function such that $d(x) > 0$ for each $x \in X$,
\,$\nu_d(A) \not\apr 0$ at least for one, and \,$\nu_d(A) \in \FR$
for each internal set $A \sbs \Xf$. The strict positivity of \,$d$
entails that the formula
$$
\|f\|_{_{p,d}} = \|f\|_{_p} =
\biggl(\,\sum_{x \in X} \lv f(x)\rv^p d(x)\biggr)^{1/p}
$$
defines an internal norm on the linear space $\sbC^X$ for each real
number $p \ge 1$. Particularly, for each internal set $A \sbs X$,
we have
$$
\nu_d(A) = \sum_{a \in A} d(a) =
\|1_A\|_{_1} = \|1_A\|_{_p}^{\,p}\,,
$$
where $1_A\: X \to \{0,1\}$ is the {\it indicator\/} or
{\it characteristic function\/} of the subset $A$ in $X$.

Suppressing $d$ in our notation, we denote by
$$\align
\IpC^X &= \bigl\{f \in \sbC^X\:
          \|f\|_{_p} \apr 0\bigr\}\,, \\
\FpC^X &= \bigl\{f \in \sbC^X\:
          \|f\|_{_p} < \infty\bigr\}\,,
\endalign$$
the $\FC$-linear subspaces of \,$\sbC^X$, consisting of internal
functions which are infinitesimal or finite, respectively, with
respect to the $p$-norm $\|\!\cd\!\|_{_p}$. Thus we get the IMG
triplet $\bigl(\sbC^X,\IpC^X,\FpC^X\bigr)$.

We also fix the notation $\mb = \mbd$ for the nonnegative
(and not identically 0) regular Borel measure on $\Xb = X^\be$
induced by $d$, and $\Lbp(\Xb) = \Lbp(\Xb,\mb)$ for the
corresponding Lebesgue spaces with the norms
$$
\|\fb\|_{_p} = \biggl(\,\int |\,\fb\,|^p \dif\mb\biggr)^{1/p}\,.
$$
We will relate them to some subspaces of \,$\FpC^X$ and
of the observable trace (nonstandard hull)
$\bigl(\sbC^X\bigr)^\be_p = \FpC^X\!/\,\IpC^X$.

Let $\Mb(\Xb)$ denote the Banach space of all complex regular Borel
measures $\mmu$ on $\Xb$ with finite total variation $\|\mmu\|$.
According to the Riesz representation theorem, $\Mb(\Xb)$ is
isomorphic to the dual $\Cbo(\Xb)^\star$ of the Banach space
$\Cbo(\Xb)$. By the Radon-Nikodym theorem, $\Lb^1(\Xb)$ can be
identified with the closed subspace of all measures
$\mmu \in \Mb(\Xb)$ absolutely continuous with respect to $\mb$.

This, together with the representation of functions
$\fb \in \Cbo(\Xb)$ by their liftings, which are internal functions
$f \in \CCo(X,E,\Xf)$ (see Proposition~1.3.1), justifies the
following notion. An internal function $g \in \sbC^X$ is called
a {\it weak lifting\/} of the measure $\mmu \in \Mb(\Xb)$, if
$$
\int f^\be \dif\mmu = {\vphantom{\Bigl)}}^\co
\biggl(\,\sum_{x \in X} f(x)\,g(x)\,d(x)\biggr)\,,
$$
for every function $f \in \CCo(X,E,\Xf)$. This is obviously
equivalent to the condition
$$
\mmu = \tth_{gd}\,.
$$
If $\mmu$ is absolutely continuous with respect to $\mb$ and
$\dif\mmu = \gb\dif\mb$, where $\gb \in \Lb^1(\Xb)$, then
$g \in \sbC^X$ is called a {\it weak lifting\/} of \,$\gb$ if \,$g$
is a weak lifting of the measure $\mmu$, i.e., if and only if
$$
\int f^\be \gb \dif\mb = {\vphantom{\Bigl)}}^\co
\biggl(\,\sum_{x \in X} f(x)\,g(x)\,d(x)\biggr) \,,
$$
for every function $f \in \CCo(X,E,\Xf)$.

Before formulating what we have just proved, let us introduce some
notation and terminology. Showing explicitly the weight function
$d$ we denote by $\MM(X,\Xf,d)$ the $\FC$-linear subspace of
\,$\sbC^X$ consisting of all internal functions $g\:X \to \sbC$
satisfying
$$
\|g\|_{_1} < \infty\,
\qquad\text{and}\qquad
\|g \cd 1_Z\|_{_1} \apr 0
$$
for each internal set $Z  \sbs X \sms \Xf$. The last condition
simply says that the Loeb measure $\lam_{|g|d}$ is concentrated on
the galaxy of accessible elements $\Xf$. Therefore, if
\,$g \in \MM(X,\Xf,d)$, then
$$
\int f^\be \gb \dif\mb = {\vphantom{\Bigl)}}^\co
\biggl(\,\sum_{x \in X} f(x)\,g(x)\,d(x)\biggr)\,,
$$
holds even for all $f \in \CCb(X,E,\Xf)$. Now, the $\FC$-linear
subspace $\SS(X,\Xf,d)$ of \,$\sbC^X$ consists of all functions
$g \in \MM(X,\Xf,d)$, satisfying additionally
$$
\nu_d(A) \apr 0 \imp \|g \cd 1_A\|_{_1} \apr 0
$$
for each internal set $A \sbs X$; functions $g \in \SS(X,\Xf,d)$
are called {\it $S$-integrable}. Obviously, the last condition is
equivalent to absolute continuity of the Loeb measure $\lam_{|g|d}$
with respect to the Loeb measure $\lam_d$, as well as to absolute
continuity of \,$\tth_{|g|d}$ with respect to $\mb$. Summing up,
we have

\proclaim{1.4.3.~Proposition}
{\sl(a)} Every measure $\mmu \in \Mb(\Xb)$ has a weak lifting
$g \in \MM(X,\Xf,d)$ such that $\|\mmu\| = {}^\co\|g\|_{_1}$.
Conversely, every function $g \in \FjC^X$, in particular, every
$g \in \MM(X,\Xf,d)$, is a weak lifting of the measure
$\tth_{gd} \in \MM(X,\Xf,d)$.

{\sl(b)} A measure $\mmu \in \Mb(\Xb)$ has a weak lifting
$g \in \SS(X,\Xf,d)$ if and only if $\mmu$ is absolutely continuous
with respect to the measure $\mb$. Conversely, every function
$g \in \SS(X,\Xf,d)$ is a weak lifting of the measure
$\tth_{gd} \in \MM(X,\Xf,d)$ (which is absolutely continuous with
respect to $\mb$).
\endproclaim

Now, there arises a natural question, namely what's the relation
between the weak lifting $g \in \SS(X,\Xf,d)$ of an absolutely
continuous measure $\mmu \in \Mb(\Xb)$ and the (unique) function
$\gb \in \Lb^1(\Xb)$ such that $\dif\mmu = \gb\dif\mb$, i.e.,
between $\gb$ and its weak lifting $g$. Unless $\gb$ is continuous,
we cannot have $\gb = g^\be$, and, unless $g$ is $S$-continuous on
$\Xf$, the formula for $g^\be$ doesn't make sense. Nevertheless,
we can still generalize the original notion of lifting of
continuous functions in the following sense. An internal function
$g\:X \to \sbC$ is called a {\it lifting\/} of a function
$\gb\:\Xb \to \bC$ (with respect to the weight function $d$)
if the equality
$$
\gb\bigl(x^\be\bigr) = {}^{\co\!}g(x)
$$
holds for {\it almost all\,} $x \in \Xf$ with respect to the Loeb
measure $\lam_d$. As the function $\gb \in \Lb^1(\Xb)$ is determined
up to the equality almost everywhere with respect to the measure
$\mb = \mbd$, only, this is the maximum one can expect.

\proclaim{1.4.4.~Proposition}
{\sl(a)} Let \,$\gb \in \Lb^1(\Xb)$ and \,$g \in \FjC^X$. Then $g$
is a weak lifting of \,$\gb$ if and only if \,$g$ is a lifting of
\,$\gb$.

{\sl(b)} Let \,$\gb\:\Xb \to \bC$. Then the following conditions
are equivalent:
{\parindent 25pt
\itemitem{\sl(i)}
\,$\gb \in \Lb^1(\Xb)$;
\itemitem{\sl(ii)}
\,$\gb$ has a weak lifting $g \in \SS(X,\Xf,d)$;
\itemitem{\sl(iii)}
\,$\gb$ has a lifting $g \in \SS(X,\Xf,d)$.
}
\endproclaim

\demo{Sketch of proof}
(a) If \,$g \in \FjC^X$ is a lifting of \,$\gb$, then, obviously,
it is a weak lifting of \,$\gb$. The reversed implication follows
from Proposition~1.3.1(c) and the uniqueness part of the
Radon-Nikodym  theorem.

(b) As the implications (iii)$\imp$(ii)$\imp$(i) are obvious,
it suffices to prove (i)$\imp$(iii). To this end denote by
$\mmu \in \Mb(\Xb)$ the measure satisfying $\dif\mmu = \gb\dif\mb$,
and by $\wtl g \in \sbC^X$ the internal function guaranteed to
$\mmu$ in Proposition~1.4.2. Then the function $g = \wtl g/d$
has all the required properties, and
${}^\co g(x) = \gb\bigl(x^\be\bigr)$ for $\lam_d$-almost all
$x \in \Xf$, due to the uniqueness part of the Radon-Nikodym
theorem, again.
\enddemo

\remark{Remark}
It is worthwhile to notice that, for a ``typical'' $x \in \Xf$,
$$
g(x) \apr \frac{{}^*\mmu\bigl(\{\xb \in \Kb_0\:
(\eta\co\sig)(\xb)=\eta(x)\}\bigr)}
{{}^*\mb\bigl(\{\xb \in \Kb_0\:
(\eta\co\sig)(\xb)=\eta(x)\}\bigr)}\,,
$$
where the right hand term is the mean value of the function
${}^*\gb = {}^*(\dif\mmu/\!\dif\mb)$ on the set
$\{\xb \in \Kb_0\: (\eta\co\sig)(\xb)=\eta(x)\}
\sbs \Ub_0[\eta(x)]$.
This is in accord with the intuition that the Radon-Nikodym
derivative $\gb(\xb) = (\dif\mmu/\!\dif\mb)(\xb)$ is the ratio
$\mmu(\Vb)/\mb(\Vb)$ of measures of some ``infinitesimal
neighborhood'' $\Vb$ of the point $\xb \in \Xb$.
\endremark
\medskip

Unfortunately, not every $S$-integrable function is lifting of
some function $\gb \in \Lb^1(\Xb)$. For instance, every function
$g \in \FiC^X$ with internal support
$$
\supp g = \{x \in X\: g(x) \ne 0\}
$$
contained in $\Xf$ is $S$-integrable, however, unless $E$ is
internal, such a function need not be lifting of any function
$\gb \in \Lb^1(\Xb)$. One can naturally expect that, in order to
lift a function $\gb \in \Lb^1(\Xb)$, the function $g \in \sbC^X$
has to display some ``reasonable amount'' of continuity, which is
not clear for the moment. This leads us to define the external
subspace $\LL^1(X,E,\Xf) \sbs \sbC^X$ as the space of all
internal functions $g \in \MM(X,\Xf,d)$ which are liftings of
functions $\gb \in \Lb^1(\Xb)$. Further we put
$$\align
\MM^p(X,\Xf,d) &=
\bigl\{f \in \sbC^X\: \lv f\rv^p \in \MM(X,\Xf,d)\bigr\}\,,\\
\SS^p(X,\Xf,d) &=
\bigl\{f \in \sbC^X\: \lv f\rv^p \in \SS(X,\Xf,d)\bigr\}\,,\\
\LL^p(X,E,\Xf) &=
\bigl\{f \in \sbC^X\: \lv f\rv^p \in \LL^1(X,E,\Xf)\bigr\}\,,
\endalign$$
for $1 \le p < \infty$. Obviously, all the functions in
$\LL^1(X,E,\Xf)$ are $S$-integrable, hence
$\LL^p(X,E,\Xf) \sbs \SS^p(X,\Xf,d)$, and the subspaces
$\LL^p(X,E,\Xf)$ are formed by the liftings of functions
$\gb \in \Lbp(\Xb)$.

Presently we do not dispose of a more explicit description
of the spaces $\LL^p(X,E,\Xf)$. However, in case that the triplet
$(X,E,\Xf)$ corresponds to a locally compact abelian group $\Gb$ in
a sense to be made precise in the next section and $\mb_d$ is the
Haar measure on $\Gb$, we will give a characterization of functions
in $\LL^p(X,E,\Xf)$ as those belonging to $\MM^p(X,\Xf,d)$ and
satisfying certain natural continuity condition, indeed (see
Theorem~3.1.4).

From the definition of \,$\LL^p(X,E,\Xf)$ and the last Proposition
we readily obtain the following result, justifying our notation.

\proclaim{1.4.5.~Proposition}
Let \,$1 \le p < \infty$. Then the Lebesgue space $\Lbp(\Xb)$ is
isomorphic to the closed subspace $\LL^p(X,E,\Xf)/\IpC^X$
of the nonstandard hull \,$\FpC^X\!/\,\IpC^X$.
\endproclaim

\remark{Remark}
Though, in general, $\LL^p(X,E,\Xf)$ is a proper subspace
of \,$\SS^p(X,\Xf,d)$, from 1.4.3 and 1.4.4 it follows that
$\LL^p(X,E,\Xf)$ is dense in $\SS^p(X,\Xf,d)$ with respect to
some rather natural weak topology which we need not to describe
precisely here.
\endremark

\specialhead
1.5.~Bounded monadic groups
\endspecialhead
\flushpar
A {\it bounded monadic group\/} is an ordered triple $(G,\Go,\Gf)$
consisting of an internal group $G$ (which means that $G$ is
an internal set, endowed with internal operations of group
multiplication and taking inverses), a monadic subgroup
$\Go \sbs G$ and a galactic subgroup $\Gf \sbs G$, such that
$\Go \nsg \Gf$ (i.e., $\Go$ is a normal subgroup of \,$\Gf$).
Intuitively, $\Go$ is viewed as the subgroup of infinitesimals
and $\Gf$ is viewed as the subgroup of finite elements in $G$.
Bounded monadic groups will be alternatively referred to as
{\it IMG group triplets\/} (cf\. \cite{ZZ}).

Every IMG group triplet gives rise to two IMG~spaces:
$(G,E_l,\Gf)$ and $(G,E_r,\Gf)$, where $E_l$, $E_r$ denote the left
and the right equivalence relation on $G$ corresponding to~$G_0$,
respectively. Though they may differ and induce different
uniformities on the ambient group $G$, they still induce the same
uniformity on $\Gf$. In other words, the bounded monadic spaces
$(G,E_l,\Gf)$, $(G,E_r,\Gf)$ are isomorphic via the identity mapping
$\Id_G\: G \to G$.

The group $\Gf$, as well as the observable trace $G^\be = \Gf/\Go$,
endowed with the topologies described in Section~1.2, become
topological groups, and the observable trace map $x \mto x^\be$ is
a~continuous surjective homomorphism of topological groups
$\Gf \to G^\be$. (On the other hand, unless $\Go \nsg G$, that
topology on $G$ does not turn it into a topological group.)

Now, the systems $\RR$ and $\BB$ from 1.2 can be replaced by a
single system $\QQ$ of admissible size, directed both upward and downward,
consisting of symmetric internal subsets of \,$G$, such that
$$\gather
(\all Q \in \QQ)(\exs R, S \in \QQ)
\bigl(R^2  \sbs Q \et Q^2 \sbs S\bigr)\,,\\
(\all Q, R \in \QQ)(\exs S \in \QQ)
   \biggl(\,\bigcup_{x \in R} xSx^{-1} \sbs Q\biggr)\,,\\
\Go = \bigcap \QQ\,,
\qquad\qquad
\Gf = \bigcup \QQ\,.
\endgather$$
If $\Go$ is the intersection and $\Gf$ is the union of {\it countably
many\/} internal sets, then we can assume that $\QQ =\{Q_n\: n \in \bZ\}$
and the symmetric internal sets $Q_n \sbs G$ satisfy
$$\gather
(\all n \in \bZ)\bigl(Q_n^2 \sbs Q_{n+1}\bigr)\,, \\
(\all n \in \bN)\biggl(
   \,\bigcup_{x \in Q_n} x\,Q_{-n-1}x^{-1} \sbs Q_{-n}\biggr)\,, \\
\Go = \bigcap_{n \in \bZ} Q_n = \bigcap_{n \in \bN} Q_{-n}\,,
\qquad\qquad
\Gf = \bigcup_{n \in \bZ} Q_n = \bigcup_{n \in \bN} Q_n\,.
\endgather$$

An internal function $\ro\:G \to \sbR$ is called a {\it valuation\/}
on the internal group $G$ if
$$\gather
\ro(x) = 0 \iff x = 1\,, \\
\ro(x) = \rho\bigl(x^{-1}\bigr)\,, \\
\ro(xy) \le \rho(x) + \ro(y)\,,
\endgather$$
for all $x,y, \in G$ ($\ro(x) \ge 0$ already follows). If $\ro$ is
a valuation on $G$ then, of course, $\ro\bigl(x^{-1}y\bigr)$ is a
left invariant and $\ro\bigl(xy^{-1}\bigr)$ is a right invariant
metric on $G$. The set of all valuations on $G$ is ordered
by the relation
$$
\ro \le \sig \iff (\all x \in G)(\ro(x) \le \sig(x))\,.
$$
A set $\VV$ of valuations on $G$ is called {\it bidirected\/} if
it is both downward and upward directed with respect to the
just introduced order relation $\le$.

Using a slightly modified Birkhoff-Kakutani style argument, one can
prove the following version of metrization theorem for bounded monadic
group triplets.

\proclaim{1.5.1.~Proposition}
Let $(G,\Go,\Gf)$ be an IMG group triplet. If \,$\Go$ is
intersection and \,$\Gf$ is union of countably many internal sets,
then there is a valuation $\ro\:G \to \sbR$ such that
$$\align
\Go &= \{x \in G\: \ro(x) \apr 0\}\,,\\
\Gf &= \{x \in G\: \ro(x) < \infty\}\,.
\endalign$$
In general, there is a bidirected set \,$\VV$ of admissible
size of valuations on $G$, such that
$$\align
\Go &= \bigl\{x \in G\:
   (\all \ro \in \VV)(\ro(x) \apr 0)\bigr\}\,,\\
\Gf &= \bigl\{x \in G\:
   (\exs \ro \in \VV)(\ro(x) < \infty)\bigr\}\,.
\endalign$$
\endproclaim

\demo{Sketch of proof}
In the countable case we give just the formula for $\ro$; the
verification that it has all the required properties is just
a matter of skill.

There is a sequence $(A_n)_{n\in\bZ}$ of symmetric internal subsets
of $G$ such that
$$
\Go = \bigcap_{n\in\bZ} A_n\,, \qquad
\Gf = \bigcup_{n\in\bZ} A_n\,,
$$
and
$$
A_n \cd A_n \cd A_n \sbs A_{n+1}
$$
for each $n \in \bZ$. By saturation, it can be extended to
an internal sequence $(A_n)_{-m \le n \le m}$ for some
$m \in \sbN_\infty$, such that $A_{-m} = \{0\}$, $A_m = G$,
and the above inclusions hold for all $n$ such that $-m \le n < m$.
Let us put
$$
\mu(x) = \min\,\{n\: -m \le n \le m\ \&\ x \in A_n\}\,,
$$
and denote by
$$
\bold P(G) = \bigcup_{n\in\sbN} G^n
$$
the internal set of all hyperfinite internal progressions in $G$.
For any progression $\bold x = (x_1,\dots,x_n) \in \bold P(G)$
we denote by $\lv \bold x\rv = n$  its length and put
$$
\varPi(\bold x) = \prod_{i=1}^n x_i\,,
\qquad\text{and}\qquad
w(\bold x) = \sum_{i=1}^n 2^{\mu(x_i)}\,;
$$
for the empty progression $() = \emptyset$ (i.e., in case $n=0$)
this means that $\varPi(\emptyset) = 1$ and $w(\emptyset) = 0$.
Then, finally
$$
\ro(x) = {}^{*\!}\inf\,\bigl\{w(\bold x)\:
\bold x \in \bold P(G)\ \&\ \varPi(\bold x) = x\bigr\}
$$
is the desired internal valuation.

In the general case let us invoke the system $\QQ$ introduced in
the beginning of this section. For each $Q \in \QQ$ there is a
sequence $(A^Q_n)_{n \in \bZ}$ of sets from $\QQ$ such that
$A^Q_0 = Q$,
$$
\Go \sbs \bigcap_{n\in\bZ} A^Q_n\,, \qquad
\bigcup_{n\in\bZ} A^Q_n \sbs \Gf\,,
$$
and
$$
A^Q_n \cd A^Q_n \cd A^Q_n \sbs A^Q_{n+1}
$$
for each $n$. Let $\ro^Q$ be the valuation constructed from
(some hyperfinite extension of) this sequence. Let us denote
$$
\VV_0 = \bigl\{\rho^Q\: Q \in \QQ\bigr\}\,.
$$
Then $\VV_0$ obviously has all the properties required for $\VV$,
except perhaps for the bidirectedness. To fix this issue, we define
$$\align
(\rho_1 \land \rho_2)(x) &=
{}^{*\!}\inf\,\Biggl\{\sum_{i=1}^{|\bold x|}
   \min\{\ro_1(x_i), \ro_2(x_i)\}\:
\bold x \in \bold P(G) \et \varPi(\bold x) = x\Biggr\}\,,\\
(\rho_1 \lor \rho_2)(x) &= \max\{\ro_1(x),\ro_2(x)\}\,,
\endalign$$
for any valuations $\ro_1$, $\ro_2$ and $x \in G$. As easily seen,
both $\ro_1 \land \ro_2$, $\ro_1 \lor \ro_2$ are valuations and
$$
\ro_1 \land \ro_2 \le \ro_1,\ro_2 \le \ro_1 \lor \ro_2\,.
$$
Taking for $\VV$ the closure of \,$\VV_0$ with respect to the
operations $\land$ and $\lor$, we are done.
\enddemo

\remark{Remark}
The last metrization theorem can be proved even for more general
triplets $(G,\Go,\Gf)$, consisting of an internal group $G$,
a monadic subgroup $\Go$ and a galactic subgroup $\Gf$, such that
$\Go \sbs \Gf \sbs G$, without the assumption that $\Go$ is normal
in~$\Gf$. Conversely, any single valuation $\ro$, as well as any
downward directed set $\VV$ of admissible size of valuations on
an internal group $G$ gives rise to a monadic subgroup $\Go$ and
a galactic subgroup $\Gf$, {\it defined\/} by the formulas from
1.5.1, such that $\Go \sbs \Gf \sbs G$. However, one cannot prove
$\Go \nsg \Gf$, in general.
\endremark
\medskip

Topological groups $\Gb$ embeddable into observable traces
$\Gb \iso G^\be$ of IMG group triplets can be easily characterized.

Let $\Gb$ be a topological group and $\Ub$ be a symmetric
neighborhood of the unit element $1 \in \Gb$. Then $\Gb$ is called
{\it $\Ub$-locally uniform\/} if the group multiplication in $\Gb$
restricted to the set $\Ub \cx \Ub$ is uniformly continuous in
the left (or, equivalently, in the right) uniformity on $\Gb$
(cf\. \cite{Gn2}). $\Gb$ is {\it locally uniform\/} if it is
$\Ub$-locally uniform for some $\Ub$.

Obviously, any subgroup of a locally uniform topological group
is itself locally uniform (in the subgroup topology).

The easy proof of the following nonstandard formulation of
\,$\Ub$-local uniformity is left to the reader.

\proclaim{1.5.2.~Lemma}
Let \,$\Gb$ be a topological group and \,$\Ub$ be a symmetric
neighborhood of the unit element \,$1 \in \Gb$. Denote by $\Nb$
the normalizer of the monad \,$\Mon(1)$ in $\sGb$. Then $\Gb$ is
$\Ub$-locally uniform if and only if \,$\sUb \sbs \Nb$.
\endproclaim

\proclaim{1.5.3.~Proposition}
Let \,$\Gb$ be a Hausdorff topological group. Then $\Gb$ can
be embedded into the observable trace $G^\be = \Gf/\Go$ of some
bounded monadic group $(G,\Go,\Gf)$ if and only if \,$\Gb$ is
locally uniform.
\endproclaim

\demo{Proof}
Assume that $\Gb$ is isomorphic to the observable trace $G^\be$
of some IMG group triplet $(G,\Go,\Gf)$. Then, as $\Go \nsg \Gf$,
we have
$$
(\all x_1,x_2,y_2,y_2 \in \Gf)
(x_1 \apr y_1 \et x_2 \apr y_2 \imp x_1 x_2 \apr y_1 y_2)\,,
$$
hence the multiplication in $G$ is $S$-continuous on $\Gf \cx \Gf$,
yielding uniform continuity of the multiplication in $\Gb$ on every
set of the form $U^\be \cx U^\be$, where $U$ is internal and
$\Go \sbs U \sbs \Gf$. Since such sets $U^\be$ form a neighborhood
base of \,$1 \in \Gb$, \,$\Gb$ is locally uniform. If $\Gb$ is just
embedded into the observable trace $G^\be$, then $\Gb$ is isomorphic
to a subgroup of a locally uniform group, hence it is locally
uniform, as well.

Now, assume that $\Gb$ is locally uniform. Let $\kappa$ be the least
uncountable cardinal, such that the topology of \,$\Gb$ has a base
$\BB$ of cardinality $< \kappa$, and $\sGb$ be a nonstandard extension
of \,$\Gb$ in a $\kappa$-saturated nonstandard universe. Let
$\Ub \in \BB$ be a symmetric neighborhood of \,$1 \in \Gb$ such
that the group multiplication is uniformly continuous on
$\Ub \cx \Ub$. Denote by $\IG$ the monad of the unit element
$1 \in \sGb$ and by $\Nb$ its normalizer in $\sGb$.
Finally we put
$$
\Sbb = \bigcup\bigl\{{}^*\Bb\:
        \Bb \in \BB \et {}^*\Bb \sbs \Nb\bigr\}\,,
$$
and denote by $\FG = \br\Sbb\kt$ the subgroup of \,$\sGb$
generated by $\Sbb$. Then, obviously, $\FG$ is a union of admissibly
many internal sets, and $\sUb \sbs \FG \sbs \Nb$, which means
that $\IG \nsg \FG$. The inclusion $\Gb \sbs \FG$ follows from the
fact that for each $\ab \in \Gb$ the group operation is uniformly
continuous on $\Ub\ab \cx \Ub\ab$, hence also on $\Bb \cx \Bb$ where
$\Bb \in \BB$ is a neighborhood of $\ab$ such that $\Bb \sbs \Ub\ab$.

Thus $(\sGb,\IG,\FG)$ is an IMG group triplet and the star
map $\ab \mto {}^*\ab$ induces an embedding of \,$\Gb$ into
its observable trace $\sGb^\be = \FG/\,\IG$.
\enddemo

For locally compact groups even more can be proved.

\proclaim{1.5.4.~Proposition}
Let \,$\Gb$ be a Hausdorff locally compact topological group.
Then, in a sufficiently saturated nonstandard universe, $\Gb$
is isomorphic to the observable trace $\Ns(\sGb)/\Mon(1)$
of the condensing IMG group triplet
$\bigl(\sGb,\Mon(1),\Ns(\sGb)\bigr)$.
\endproclaim

\demo{Sketch of proof}
It suffices to have the nonstandard universe $\kappa$-saturated
where $\kappa$ is the least uncountable cardinal bigger than the
cardinality of some neighborhood base of \,$1 \in \Gb$ and such
that $\Gb$ can be covered by the interiors of less than $\kappa$
compact sets.
\enddemo

We are mainly interested in {\it condensing\/} IMG group triplets
$(G,\Go,\Gf)$ with a {\it hyperfinite\/} ambient group $G$.
Condensing IMG triplets can be characterized using Proposition~1.2.1.
Below $[A : B]$ denotes any of the indices $\lfloor A:B\rfloor$,
$\lfloor A:B\rfloor_\ii$ or $\lceil A:B\rceil$. The simple proof
of the following facts is left to the reader.

\proclaim{1.5.5.~Proposition}
{\sl (a)} An IMG group triplet $(G,\Go,\Gf)$ is condensing if and only if for
any symmetric internal sets $A$, $B$ between $\Go$ and $\Gf$ the
index $[A : B]$ is finite. If \,$G$ is hyperfinite, then this is
equivalent to
$$
0 \not\apr \frac{\lv A\rv}{\lv B\rv} < \infty
$$
for any internal sets $A$, $B$ between $\Go$ and $\Gf$.

{\sl (b)} The observable trace $G^{\be}$ of an condensing group triplet
$(G,\Go,\Gf)$ is discrete if and only if \,$\Go$ is an internal group;
$G^{\be}$ is compact if and only  if \,$\Gf$ internal.
\endproclaim

Given a condensing IMG group triplet $(G,\Go\Gf)$ with a
hyperfinite ambient group $G$, a positive number $d \in \sbR$, such
that $d\lv A\rv \in \FR \sms \IR$ for some (or, equivalently, for
each) internal set $A$ between $\Go$ and $\Gf$, is called a
{\it normalizing multiplier\/} or {\it normalizing coefficient\/}
for $(G,\Go,\Gf)$. According to 1.5.5, if $d$ is a normalizing
multiplier for $(G,\Go,\Gf)$, and $0 < d' \in \sbR$, then $d'$ is
a normalizing multiplier if and only if $d/d' \in \FR \sms \IR$.
In particular, for any internal set $A$ between $\Go$ and $\Gf$,
$d = 1/\lv A\rv$ is a normalizing coefficient for $(G,\Go,\Gf)$.

From the results of Section~1.4 it follows directly

\proclaim{1.5.6.~Proposition}
Let \,$(G,\Go,\Gf)$ be a condensing IMG group triplet with a
hyperfinite ambient group $G$ and a normalizing multiplier $d$.
Let $\lam_d$ denote the Loeb measure induced by the constant
function $d(x) = d$ on $G$. Then the measure $\mbd$ obtained by
pushing down the Loeb measure $\lam_d$ is both left and right
invariant Haar measure on the observable trace $G^\be = \Gf/\Go$.
\endproclaim

\remark{Remark}
A reasonable characterization of locally compact topological
groups isomorphic to observable traces of condensing IMG group
triplets with hyperfinite ambient group is still missing. According
to Proposition~1.5.6, every locally compact group representable as an
observable trace of such a group triplet is necessarily unimodular.
On the other hand, there are even compact topological groups, as,
e.g., $\SO(3)$, which do not admit any such representation (see
\cite{GG}, \cite{GGR}). Similarly, not even all finitely generated
discrete groups admit such a representation (cf\.~\cite{AGG},
\cite{GV}). Nevertheless, as we shall see latter on, locally
compact {\it abelian\/} groups behave well.
\endremark
\medskip

Let $\Gb$ be a group, $\Ub \sbs \Gb$ be any set containing the
unit $1 \in \Gb$ and $\Kb$ be a nonempty subset of \,$\Gb$.
A~mapping $\eta\: G \to \Gb$ is called a {\it finite
$(\Kb,\Ub)$ approximation\/} of \,$\Gb$ if $G$ is a finite group,
and $\eta$ satisfies the following two conditions:
$$\gather
(\all \xb \in \Kb)(\exs x \in G)\bigl(\eta(x) \in \Ub\xb\bigr)\,,\\
(\all x,y \in G)\bigl(\eta(x),\eta(y) \in \Kb
   \imp \eta(x)\eta(y) \in \Ub\eta(xy)\bigr)\,.
\endgather$$
A $(\Kb,\Ub)$ approximation is called {\it injective\/} if \,$\eta$
is an injective mapping; it is called {\it strict\/} if
$\eta(1) = 1$, $\eta\bigl(x^{-1}\bigr) = \eta(x)^{-1}$ for all
$x \in G$, and the second of the above conditions can be
strengthened to
$$
(\all x,y \in G)\bigl(\eta(x),\eta(y) \in \Kb
   \imp \eta(xy) = \eta(x)\eta(y)\bigr)\,.
$$
Notice that any $(\Kb,\Ub)$ approximation $\eta\:G \to \Gb$
satisfies $\eta(1) \in \Ub$ if \,$1 \in \eta^{-1}[\Kb]$,
and $\eta\bigl(x^{-1}\bigr) \in \Ub\eta(x)^{-1}$ if
\,$1,x,x^{-1} \in \eta^{-1}[\Kb]$. Similarly, the strengthened
second condition implies $\eta(1) = 1$ if \,$1 \in \eta^{-1}[\Kb]$,
and $\eta\bigl(x^{-1}\bigr) = \eta(x)^{-1}$ if
\,$1,x,x^{-1} \in \eta^{-1}[\Kb]$. However, the convenient
conditions $\eta(1) = 1$, $\eta\bigl(x^{-1}\bigr) = \eta(x)^{-1}$
alone can always be assumed without loss of generality.

If \,$\Gb$ is a Hausdorff topological group and $(I,\le)$ is
a directed poset, then a system of mappings
$(\eta_i\:G_i \to \Gb)_{i \in I}$, where each $G_i$ is a finite
group, is called an {\it approximating system\/} of \,$\Gb$ if
for every compact set $\Kb \sbs \Gb$ and every neighborhood
$\Ub$ of the unit in $\Gb$ there is an $i \in I$ such that
$\eta_j\:G_j \to \Gb$ is a $(\Kb,\Ub)$ approximation of \,$\Gb$
for each $j \in I$, $j \ge i$.

A system $\bigl((\Kb_i,\Ub_i)\bigr)_{i \in I}$ of pairs of subsets
of \,$\Gb$ is called a {\it directed double base\/} of \,$\Gb$,
briefly a {\it DD~base}, if the sets $\Kb_i$ are compact and their
interiors cover $\Gb$, the sets $\Ub_i$ form a neighborhood base of
\,$1 \in \Gb$, and, for all $i,j \in I$, $i \le j$ implies
$$
\Ub_j \sbs \Ub_i \sbs \Kb_i \sbs \Kb_j\,.
$$
An approximating system $(\eta_i\:G_i \to \Gb)_{i \in I}$ is
called {\it well based\/} if there is a directed double base
$\bigl((\Kb_i,\Ub_i)\bigr)_{i \in I}$ such that each $\eta_i$
is a $(\Kb_i,\Ub_i)$ approximation of \,$\Gb$. Is it the case,
then we say that the approximating system $(\eta_i)_{i \in I}$
{\it is well based with respect to the DD base\/}
$\bigl((\Kb_i,\Ub_i)\bigr)_{i \in I}$.

For the sake of convenience it is useful to realize that, given a
DD base $\bigl((\Kb_i,\Ub_i)\bigr)_{i \in I}$ of \,$\Gb$ over a
directed poset $(I,\le)$, any system of mappings
$(\eta_i\:G_i \to \Gb)_{i \in I}$ such that each $\eta_i$ is a
$(\Kb_i,\Ub_i)$ approximation of \,$\Gb$ by some finite group
$G_i$ is already an approximating system of \,$\Gb$ (which is,
of course, well based with respect to that DD base).

If $\Gb$ is discrete, then it is enough to deal with its DD bases
consisting just of pairs $(\Kb_i,\{1\})$, where $\Kb_i$ are finite
sets whose union  is $\Gb$. Similarly, if $\Gb$ is compact, then it
is enough to consider its DD bases with members of form $(\Gb,\Ub_i)$
where $(\Ub_i)_{i \in I}$ is a neighborhood base of \,$1 \in \Gb$.

An internal mapping $\eta\:G \to \sGb$ is a {\it hyperfinite
infinitesimal approximation}, briefly {\it HFI approximation},
of \,$\Gb$ if \,$G$ is a hyperfinite group and
$$
\gather
(\all \xb \in \Gb)(\exs x \in G)\bigl(\eta(x) \apr \xb\bigr)\,,\\
(\all x,y \in G)\bigl(\eta(x),\eta(y) \in \Ns(\sGb)
\imp \eta(xy) \apr \eta(x)\eta(y)\bigr)\,.
\endgather
$$
The notions of {\it injective\/} and {\it strict\/} HFI
approximation, respectively, are defined in the obvious way.

Every hyperfinite infinitesimal approximation $\eta\:G \to \sGb$
of a Hausdorff locally compact group $\Gb$ in a sufficiently
saturated nonstandard universe gives rise to an IMG group triplet
$(G,\Go,\Gf)$ such that
$$\align
\Go &= \eta^{-1}\bigl[\Mon(1)\bigr] =
       \{x \in G\: \eta(x) \apr 1\}\,, \\
\Gf &= \eta^{-1}\bigl[\Ns(\sGb)\bigr] =
       \{x \in G\: \eta(x) \in \Ns(\sGb)\}\,,
\endalign$$
and $\Gb$ is isomorphic to its observable trace $G^\be = \Gf/\Go$.

\proclaim{1.5.7.~Proposition}
Let $\Gb$ be a Hausdorff locally compact group. Then the following
conditions are equivalent:
{\parindent 21pt
\item{\sl(i)}
$\Gb$ is isomorphic to the observable trace $G^\be$ of
a hyperfinite condensing IMG group triplet $(G,\Go,\Gf)$;
\item{\sl(ii)}
there is an HFI approximation $\eta\: G \to \sGb$ of \,$\Gb$;
\item{\sl(iii)}
for every compact set \,$\Kb \sbs \Gb$ and every neighborhood
\,$\Ub$ of \,$1 \in \Gb$  there is a finite
$(\Kb,\Ub)$ approximation $\eta\: G \to \Gb$;
\item{\sl(iv)}
there is a well based approximating system
$(\eta_i\:G_i \to \Gb)_{i \in I}$ of \,$\Gb$ by finite groups $G_i$;
\item{\sl(v)}
there is an approximating system $(\eta_i\:G_i \to \Gb)_{i \in I}$
of \,$\Gb$ by finite groups $G_i$.
\item{}}
\vskip-12pt
\endproclaim

\demo{Proof}
(ii)$\imp$(i) was proved in advance, right before formulating the
Proposition.

(i)$\imp$(ii)
Assume that $(G,\Go,\Gf)$ is a hyperfinite condensing IMG group
triplet and $\eeta\:\Gf/\Go \to \Gb$ is an isomorphism of
topological groups. By Proposition 1.2.3, $\eeta$ has a lifting
$\eta\:G \to \sGb$. Now, it is a routine to check that $\eta$ is
an HFI approximation of \,$\Gb$.

(ii)$\imp$(iii)
Let $\eta\: G \to \sGb$ be an HFI approximation of \,$\Gb$.
For any neighborhood $\Ub$ of \,$1 \in \Gb$ and a compact set
$\Kb \sbs \Gb$ denote by $\Phi(\Kb,\Ub)$ the set of all finite
$(\Kb,\Ub)$ approximations $\phi\:F \to \Gb$ of \,$\Gb$. Then
$\eta \in {}^*\Phi(\Kb,\Ub)$. By the {\it transfer principle},
$\Phi(\Kb,\Ub) \ne \emptyset$.

(iii)$\imp$(iv)
Let $\KK$ be an upward directed set of compacts in $\Gb$ such that
the interiors of members of \,$\KK$ cover $\Gb$ and $\UU$ be a
neighborhood base of \,$1 \in \Gb$. Denote by $I$ the set of
all pairs $(\Kb,\Ub) \in \KK \cx \UU$ such that $\Ub \sbs \Kb$,
and order it by the relation
$$
(\Kb,\Ub) \le (\Kb',\Ub') \iff \Kb \sbs \Kb' \et \Ub' \sbs \Ub\,;
$$
obviously, $(I,\le)$ is a directed poset. Let $\eta_i\:G_i \to \Gb$
be a finite $(\Kb,\Ub)$ approximation of \,$\Gb$ for
$i = (\Kb,\Ub) \in I$. Clearly, the mappings $\eta_i$ form
a well based approximating system of \,$\Gb$.

(iv)$\imp$(v) is trivial.

(v)$\imp$(ii)
Let $(\eta_i\:G_i \to \Gb)_{i \in I}$ be an approximating system
of \,$\Gb$ over $(I,\le)$. Let us embed the situation into a
$\kappa$-saturated nonstandard universe, where $\kappa > \aleph_0$
and $(I,\le)$ has a cofinal subset of cardinality $< \kappa$.
Then there is a $k \in {}^{*\!}I$ such that $i \le k$ for all
$i \in I$; it follows that $\eta_k\:G_k \to \sGb$ is an HFI
approximation of \,$\Gb$.
\enddemo

\remark{Remark}
From the above proof it is clear that one can strengthen
conditions (ii)--(v) of Proposition~1.5.7 by requiring injectivity
of the corresponding approximation(s), and the modified conditions
will be equivalent, again. What is not clear is the way how one
should modify condition (i) to make it equivalent with them (maybe
even no modification is needed at all).
\endremark

\medskip
One could expect that a tighter connection than the merely
existential one stated in the last Proposition can be established
by means of the ultraproduct construction. To this end consider a
directed poset $(I,\le)$.  An {\it ultrafilter\/} $\DD$ on $I$
is called {\it directed\/} if it contains all the sets
$[i) = \{j \in I\: i \le j\}$ for $i \in I$. It can be easily seen
that any set $J \sbs I$ belonging to a directed ultrafilter $\DD$ on
$I$ is cofinal in $(I,\le)$.

We denote by $\sGb = \Gb^I\!/\DD$ the ultrapower of a group $\Gb$
and, for any set $\Xb \sbs \Gb$, we denote by
$\sXb = \Xb^I\!/\DD$ its ultrapower viewed as a subset of \,$\sGb$.
Given a system of finite groups $(G_i)_{i\in I}$, their ultraproduct
$G = \prod_{i \in I} G_i/\DD$ is a hyperfinite group. Finally, if
$(\eta_i\:G_i\to\Gb)_{i\in I}$ is a system of mappings, then
$\eta = (\eta_i)_{i \in I}/\DD\: G \to \sGb$ denotes the internal
mapping given by $\eta(a) = \bigl(\eta_i(a_i)\bigr)_{i \in I}/\DD$
for $a = (a_i)_{i \in I}/\DD \in G$.

The following statement, making use of the above ultraproduct
notation, is a more detailed version of the equivalence
(ii)$\iff$(iii) from Proposition~1.5.7.

\proclaim{1.5.8.~Proposition}
Let \,$\Gb$ be a Hausdorff locally compact group, $(I,\le)$
be a directed poset, $\bigl((\Kb_i,\Ub_i)\bigr)_{i \in I}$ be
a directed double base of \,$\Gb$ and $\DD$ be a directed ultrafilter
on~$I$. Assume that $(G_i)_{i \in I}$ is a system of finite groups
endowed with mappings $\eta_i\: G_i \to \Gb$. Then the following hold
true:
{\parindent 21pt
\item{\sl(a)}
If \,$(\eta_i\:G_i \to \Gb)_{i\in I}$ is a well based approximating
system of \,$\Gb$ with respect to the DD~base
$\bigl((\Kb_i,\Ub_i)\bigr)_{i \in I}$, then $\eta\:G \to \sGb$ is
an HFI approximation of \,$\Gb$. In such a case, $\eta$ is an
injective HFI approximation if and only if
\,$\{i \in I\: \eta_i\ \text{\rm is injective}\} \in \DD$.
\item{\sl(b)}
Let \,$\eta\:G \to \sGb$ be an HFI approximation of \,$\Gb$. Then
there is a function \hbox{$\tau\:I \to I$} such that $i \le \tau(i)$
for each $i \in I$, and
$\bigl(\eta_{\tau(i)}\: G_{\tau(i)} \to \Gb\bigr)_{i \in I}$ is
a well based approximating system of \,$\Gb$ with respect to the
DD~base $\bigl((\Kb_i,\Ub_i)\bigr)_{i \in I}$. If \,$\eta$ is
injective, then the system $\bigl(\eta_{\tau(i)}\bigr)_{i \in I}$
can be chosen to consist of injective mappings, as well.
\item{}}
\vskip-12pt
\endproclaim

\demo{Proof}
(a) Let $(\eta_i\:G_i \to \Gb)_{i \in I}$ be a well based
approximating system of \,$\Gb$ with respect to the DD~base
$\bigl((\Kb_i,\Ub_i))_{i \in I}$.

First we show that, for each $\xb \in \Gb$, there is an
$a = (a_i)_{i\in I}/\DD \in G$ such that $\eta(a) \apr \xb$,
i.e., $\eta(a) \in \sUb\xb$ for each neighborhood $\Ub$ of
\,$1 \in \Gb$. Taking arbitrary $\xb$ and $\Ub$, there is an
$i \in I$ such that $\Ub_i\xb \sbs \Ub\xb \cap \Kb_i$. As
$\eta_j$ is a $(\Kb_j,\Ub_j)$ approximation of $\Gb$ for any
$j \in I$, for $j \ge i$ there is an $a_j \in G_j$ such that
$\eta_j(a_j) \in \Ub_j\xb \sbs \Ub_i\xb$. Picking arbitrary
elements $a_j \in G_j$ for $j \not\le i$, we have
\,$[i) \sbs \{j \in I\: \eta_j(a_j) \in \Ub\xb\} \in \DD$.
Then $\eta(a) \in \sUb\xb$ by {\L}os' theorem.

The almost homomorphy condition can be expressed in form of
the implication
$$
(\all a,b \in G)(\all \xb,\yb \in \Gb)
   \bigl(\eta(a) \apr \xb \et \eta(b) \apr \yb
      \imp \eta(ab) \apr \xb\yb\bigr)\,.
$$
Its straightforward verification does not even rely on the well
baseness of the approximating system, and is left to the reader.
The characterization of injectivity of \,$\eta$ is obvious, too.

(b) Let $\eta\:G \to \sGb$ be an HFI approximation of \,$\Gb$. Then
$\eta$ is a $(\sKb_i,\sUb_i)$ approximation of \,$\sGb$, for each
$i \in I$. By {\L}os' theorem there is a set $J \in \DD$ such that
$\eta_j$ is a $(\Kb_i,\Ub_i)$ approximation of \,$\Gb$ for each
$j \in J$. As $J$ is cofinal, there is a $\tau(i) \in J \cap [i)$.
Now, it is clear that the mappings
$\eta_{\tau(i)}\: G_{\tau(i)} \to \Gb$ form an approximating system
of \,$\Gb$ over $(I,\le)$, well based with respect to its DD base
$\bigl((\Kb_i,\Ub_i))_{i \in I}$. The supplement on injectivity is
plain.
\enddemo

\remark{Remark}
It is not true that an arbitrary approximating system
$(\eta_i\:G_i \to \Gb)_{i\in I}$ of a Hausdorff locally
compact group $\Gb$ by finite groups $G_i$ gives rise to
some HFI approximation $\eta\:G \to \sGb$ for an ultraproduct
$G = \prod_{i \in I} G_i/\DD$ and an ultrapower
$\sGb = \Gb^I\!/\DD$ with respect to some directed ultrafilter
$\DD$ on $(I,\le)$. The reason is that, roughly speaking, the
existence of an $I$-indexed approximating system of \,$\Gb$ does
not imply the existence of an $I$-indexed neighborhood base of
the unit in $\Gb$. Similarly, such an ultraproduct HFI
approximation $\eta\:G \to \sGb$ of \,$\Gb$ does not automatically
yield an approximating system $(\eta_i\:G_i \to \Gb)_{i \in I}$.
Thus it is not possible to avoid mentioning the $I$-indexed
DD~base $\bigl((\Kb_i,\Ub_i))_{i \in I}$ and the function
$\tau\: I \to I$ in the formulation of the last proposition.
A~more detailed analysis of this issue will be published
separately.
\endremark
\medskip

The following examples (a), (b), (c), as well as some special cases
of (d) are essentially taken from \cite{Go2}; (b) in fact (without
using the present terminology) can be found already in \cite{Lx3}.

\example{1.5.9.~Example}
Let $1 \le n \in \sbN$, and
$\bZ_n = \bigl\{-\lfloor\frac{n-1}{2}\rfloor,\dots,-1,0, 1,
\dots,\lceil\frac{n-1}{2}\rceil\bigr\}$
be the (hyper)finite cyclic group of order~$n$, represented as the
set of absolutely smallest remainders modulo~$n$.
\smallskip
(a) If \,$n \in \bN$ and $0 \le k < n/4$, then the identity
mapping \,$\bZ_n \to \bZ$ is a strict injective
$\bigl(K,\{0\}\bigr)$ approximation of the group $\bZ$,
where $K =  \{0,\pm 1, \dots, \pm k\}$.

If \,$n \in \sbNi$, then the identity mapping \,$\bZ_n \to \sbZ$ is
a strict HFI approximation of \,$\bZ$. The IMG group triplet arising
form it has the form $\bigl(\bZ_n, \{0\}, \bZ\bigr)$.
\smallskip
(b) If \,$n \in \bN$ and $\eps > \pi/n$, then the homomorphism
$a \mto \ee^{2\pi\ii a/n}\: \bZ_n \to \bT$ is a strict
injective $(\bT,\Ub)$ approximation of the group $\bT$, where
$\Ub = \{u \in \bT\: \lv\arg u\rv \le \eps\}$.

If \,$n \in \sbNi$, then the internal homomorphism
$a \mto \ee^{2\pi\ii a/n}\:\bZ_n \to \sbT$ is a strict
injective HFI approximation of \,$\bT$. The corresponding IMG
group triplet is $(\bZ_n,\Go,\bZ_n)$, where
$$
\Go = \{a \in \bZ_n\: a/n \apr 0\}\,.
$$

(c) If \,$n \in \bN$, $0 \le k < n/4 \in \bN$, and $d$ and
$\eps > d/2$ are positive real numbers, then the mapping
$a \mto ad\:\bZ_n \to \bR$ is a strict injective $(\Kb,\Ub)$
approximation of the group $\bR$, where \,$\Kb = [-kd,kd]$,
\,$\Ub = [-\eps,\eps]$.

If $n \in \sbNi$, and $d$ is a positive infinitesimal
such that $nd \sim \infty$, then the internal mapping
$a \mto ad\:\bZ_n \to \sbR$ is a strict injective HFI
approximation of~\,$\bR$, inducing the IMG group triplet
$(\bZ_n,\Go,\Gf)$, where
$$
\Go = \{a \in \bZ_n\: ad \apr 0\}\,,
\qquad\qquad
\Gf = \{a \in \bZ_n\: \lv a d\rv < \infty\}\,.
$$

(d) Let $\Gb$ be any Hausdorff locally compact group (written
multiplicatively) and $\Ub \sbs \Kb \sbs \Gb$ be its compact open
subgroups such that $\Ub$ is normal in $\Kb$. Then the quotient
$G = \Kb/\Ub$ is a finite group and any (necessarily injective)
mapping $\eta\:G \to \Gb$ such that $\Ub\eta(x) = x$ for each
$x \in G$ is an injective $(\Kb,\Ub)$ approximation of \,$\Gb$.

If \,$\Gb$ has a DD base $\bigl((\Kb_i,\Ub_i)\bigr)_{i \in I}$
consisting of compact open subgroups $\Ub_i \nsg \Kb_i$ of \,$\Gb$,
then there are $^*$compact $^*$open subgroups $\Ub \nsg \Kb$ of
\,$\sGb$, such that $\Ub \sbs \Mon(1)$, $\Ns(\sGb) \sbs \Kb$ and
the quotient $G = \Kb/\Ub$ is a hyperfinite group. Any internal mapping
$\eta\:G \to \sGb$ such that $\Ub\eta(x) = x$ for each $x \in G$ is
an injective HFI approximation of \,$\Gb$. The IMG group triplet
$(G,\Go,\Gf)$ obtained from $\eta$ satisfies
$$\align
\Go &= \bigcap_{i \in I} \eta^{-1}[\sUb_i]
= \{\Ub\xb\: (\all i \in I)(\xb \in \Ub_i)\}\,,\\
\Gf &= \bigcup_{i \in I} \eta^{-1}[\sKb_i]
= \{\Ub\xb\: (\exs i \in I)(\xb \in \Kb_i)\}\,.
\endalign$$

Both in the finite and the hyperfinite case, $\eta$ is not a strict
approximation, unless it is a genuine homomorphism.
\endexample

Obviously, if \,$\eta_1\:G_1 \to \Gb_1$ is a finite $(\Kb_1,\Ub_1)$
approximation of \,$\Gb_1$ and \,$\eta_2\:G_2 \to \Gb_2$ is a
finite $(\Kb_2,\Ub_2)$ approximation of $\Gb_2$, then the mapping
$\eta_1 \cx \eta_2\: G_1 \cx G_2 \to \Gb_1 \cx \Gb_2$ is a finite
$(\Kb_1 \cx \Kb_2, \Ub_1 \cx \Ub_2)$ approximation of
\,$\Gb_1 \cx \Gb_2$. Moreover, it is clear, that $\eta_1 \cx \eta _2$
is strict or injective, respectively, once both the approximations
$\eta_1$, $\eta_2$ have the respective property. Analogous observations
apply to HFI approximations, as well. Thus the items (a), (b), (c) of
the last example enable to construct strict injective approximating
systems and strict injective HFI approximations for all {\it elementary
LCA groups}, i.e., LCA groups of the form $F \cx \bZ^k \cx \bT^m \cx \bR^n$
where $F$ is a finite abelian group and $k,m,n \in \bN$.

Example~1.5.9(d) gives a direct hint how to construct both
approximating systems and HFI approximations of any LCA group
$\Gb$ with a DD base $\bigl((\Kb_i,\Ub_i)\bigr)_{i \in I}$
consisting of pairs of subgroups  $\Ub_i \sbs \Kb_i$ of \,$\Gb$.
The following LCA groups are included as special cases:
\smallskip
{\parindent 21pt
\item{1.}
the compact additive groups of $\tau$-adic integers $\Delta_\tau$,
where $\tau = (\tau_n)_{n\in\bN}$ is any increasing sequence of
positive integers such that $\tau_n\,|\,\tau_{n+1}$ for each $n$;
\item{2.}
the additive LCA groups of $\tau$-adic numbers $\bQ_\tau$, where
$\tau = (\tau_n)_{n\in\bZ}$ is any sequence of positive integers,
such that $\tau_n < \tau_{n+1}$, $\tau_n\,|\,\tau_{n+1}$ for
$n \ge 0$, and $\tau_n < \tau_{n-1}$, $\tau_n\,|\,\tau_{n-1}$
for $n \le 0$;
\item{3.}
torsion (discrete abelian) groups, as they are direct (inductive)
limits of finite abelian groups;
\item{4.}
profinite (compact abelian) groups, i.e., inverse (projective)
limits of finite abelian groups.%
\footnote{Obviously, 1.~is a special case of 4.}
\item{}
\vskip-12pt}
\smallskip
In \cite{Go2} also HFI approximations of \,$\tau$-adic solenoids
$$
\Sigma_\tau = (\bR \cx \Delta_{\tau})/\{(a,a)\: a \in \bZ\}\,,
$$
not falling within the scope of 1.5.9(d), are described.

Actually, any LCA group admits arbitrarily good finite approximations
and, henceforth, HFI approximations, too. The key ingredient of the
proof of this fact is the following structure theorem for compactly
generated LCA groups (i.e., LCA groups generated by some compact subset).
Its proof is relatively elementary (though certainly not trivial), in
particular, it does not rely on the Pontraygin-van\,Kampen duality
theorem (cf.~\cite{Pn, Theorem~50}, \cite{HR1, Theorem~9.6}).

\proclaim{1.5.10.~Proposition}
Let \,$\Gb$ be a compactly generated Hausdorff LCA group. Then every
compact neighborhood \,$\Vb$ of \,$0$ in $\Gb$ contains a closed subgroup
$\Hb$ of \,$\Gb$ such that the quotient \,$\Gb/\Hb$ is an elementary
LCA group.
\endproclaim

\proclaim{1.5.11.~Finite LCA Group Approximation Theorem}
Let \,$\Gb$ be a Hausdorff LCA group. Then, for any compact set
\,$\Kb \sbs \Gb$ and any neighborhood \,$\Ub$ of \,$0 \in \Gb$,
there is a finite {\rm abelian\/} group $G$ and an injective
$(\Kb,\Ub)$ approximation $\eta\:G \to \Gb$ of \,$\Gb$, such that
$\eta(0) = 0$, and $\eta(-a) = -\eta(a)$ for each $a \in G$.
Equivalently, $\Gb$ admits some well based approximating system
$(\eta_i\:G_i \to \Gb)_{i \in I}$ by finite groups, such that all
the mappings $\eta_i$ are injective and preserve $0$ and inverses.
\endproclaim

\demo{Proof}
As $\Gb$ is locally compact, we can assume, without loss of
generality, that $\Ub \sbs \Kb$, and pick a compact symmetric
neighborhood $\Vb$ of \,0 such that $\Vb + \Vb \sbs \Ub$. According to
the above Prroposition~1.5.10, $\Vb$ contains a closed subgroup $\Hb$
such that the quotient $\Eb = \br\Kb\kt/\Hb$ of the subgroup $\br\Kb\kt$
of \,$\Gb$ generated by the compact set $\Kb$ is an elementary LCA group.
Denote by $\Vb'$ and $\Kb'$ the images of the sets $\Vb$, $\Kb$,
respectively, under the canonical projection $\psi\:\br\Kb\kt \to \Eb$.
Then $\Vb'$ is a symmetric neighborhood of \,0 and $\Kb'$ is a compact
set in $\Eb$. According to Example~1.5.9 and the comment following it,
there is a finite group $G$ and a strict injective $(\Kb',\Vb')$
approximation $\zeta\:G \to \Eb$. Let $\sig\:\Eb \to \br\Kb\kt$ be
any (necessarily injective) mapping such that $\sig(0) = 0$,
$\sig(-x) = -\sig(x)$ and $\psi(\sig(x))= x$ for $x \in \Eb$.
Then $\eta = \sig \co \zeta\:G \to \Gb$ is an injective mapping,
satisfying $\eta(0) = 0$, and $\eta(-a) = -\eta(a)$ for $a \in G$.
Straightforward arguments show that
\medskip
\centerline{$
\Kb \sbs \Kb + \Hb \sbs \eta[G] + \Vb + \Hb \sbs \eta[G] + \Ub\,,$}
\flushpar
and
\flushpar
\centerline{$
\eta(a) + \eta(b) - \eta(a+b) \in \Hb \sbs \Vb \sbs \Ub\,,$}
\medskip
\flushpar
for $a,b \in G$ whenever $\eta(a),\eta(b),\eta(a+b) \in \Kb +\Hb$.
Hence, $\eta$ is a $(\Kb,\Ub)$ approximation of \,$\Gb$. The
equivalence of the first and the second formulation is obvious
in view of Proposition~1.5.7.
\enddemo

\proclaim{1.5.12.~Corollary}%
{\bf [Hyperfinite LCA Group Approximation Theorem]}
Let \,$\Gb$ be a Hausdorff LCA group. Then, in a sufficiently saturated
nonstandard universe, there is an internal hyperfinite abelian group $G$
and an injective HFI approximation $\eta\: G \to \sGb$ of \,$\Gb$, such
that $\eta(0) = 0$ and $\eta(-a) = -\eta(a)$ for $a \in G$. It follows
that \,$\Gb$ is isomorphic to the observable trace $G^\be = \Gf/\Go$ of
the $IMG$ group triplet $(G,\Go,\Gf)$ with hyperfinite abelian ambient
group $G$, where \,$\Go = \eta^{-1}\bigl[\Mon(0)\bigr]$, \,and
\,$\Gf = \eta^{-1}\bigl[\Ns(\sGb)\bigr]$.
\endproclaim
\newpage

{}\phantom{x}
{}\bigskip

\head
2. Pontryagin-van\,Kampen Duality in Hyperfinite Ambience
\endhead
{}\bigskip\bigskip

\flushpar
In this chapter we finally come to one of the central topics of
our paper which is the study of condensing IMG group triplets with
hyperfinite abelian ambient group with regard to the celebrated
Pontryagin-van\,Kampen duality theorem. We formulate and prove the
first and second of the three Gordon's conjectures concerning them.

\specialhead
2.1.~The dual triplet
\endspecialhead
\flushpar
From now on (unless we explicitly say something else), $(G,\Go,\Gf)$
denotes a fixed but arbitrary condensing IMG group triplet with
hyperfinite abelian ambient group $G$ (which is tacitly assumed to
be externally infinite).%
\footnote{However, the reader should keep in mind that
some of our accounts remain valid for general internal abelian
ambient group $G$, as well, or require just some minor modification.}
We also fix a system $\QQ$ of admissible size, directed both downward
and upward, consisting of symmetric internal sets such that
$\Go = \bigcap \QQ$ and $\Gf = \bigcup \QQ$.

Let us denote by $\Gb = G^\be = \Gf/\Go$ the observable trace of the
triplet. The LCA group $\Gb$ gives rise to the dual group $\dGb$
of all continuous homomorphisms (characters) $\ggam\:\Gb \to \bT$,
endowed with the compact-open topology. Then $\dGb$, as an LCA
group, can itself be represented as the observable trace of some
condensing IMG group triplet $(H,\Ho,\Hf)$ with hyperfinite
abelian ambient group~$H$. One can naturally expect that the triplet
representing the dual group $\dGb$ can be constructed from the
original triplet $(G,\Go,\Gf)$ in some canonical way.

Let $\dG = \Hom(G,\sbT)$ denote the set of all internal
homomorphisms (characters) $\gama\:G \to \sbT$; then $\dG$ with the
pointwise multiplication is a hyperfinite abelian group internally
isomorphic  to~$G$ (though not in a canonical way). For any sets
$A \sbs G$, $\Gama \sbs \dG$ we define their {\it infinitesimal
annihilators\/} by
$$\align
A^{\ort} &= \bigl\{\gama \in \dG\:
   (\all a \in A)(\gama(a) \apr 1)\bigr\}\,,\\
\Gama^{\ort} &= \bigl\{a \in G\:
   (\all \gama \in \Gama)(\gama(a) \apr 1)\bigr\}\,.
\endalign$$
Obviously, $A^{\ort}$ is a subgroup of \,$\dG$ and $\Gama^{\ort}$
is a subgroup of \,$G$. For any $A,B \sbs G$,
\,$\Gama,\Dela \sbs \dG$ we have
$$
A \sbs \Gama^{\ort} \iff \Gama \sbs A^{\ort}\,,
$$
as well as
$$\xalignat{2}
A \sbs B &\imp B^{\ort} \sbs A^{\ort}\,, &
A &\sbs A^{\ort\ort}\,, \\
\Gama \sbs \Dela &\imp \Dela^{\ort} \sbs \Gama^{\ort}\,, &
\Gama &\sbs \Gama^{\ort\ort}\,,
\endxalignat$$
showing that the assignments $A \mto A^{\ort}$,
\,$\Gama \mto \Gama^{\ort}$ form a Galois connection.

We are particularly interested in the subgroups
$$\align
\Go^{\ort} &= \bigl\{\gama \in \dG\:
   (\all x \in \Go)(\gama(x) \apr 1)\bigr\}\,,\\
\Gf^{\ort} &= \bigl\{\gama \in \dG\:
   (\all x \in \Gf)(\gama(x) \apr 1)\bigr\}
\endalign$$
of the dual group $\dG$. As every $\gama \in \dG$ is a group
homomorphism, it belongs to $\Go^{\ort}$ if and only if it is
$S$-continuous as a mapping $\gama\:G \to \sbT$ with respect to
the monadic equivalence $E_l = E_r$ on $G$ and the
usual equivalence of infinitesimal nearness $\apr$ on $\sbT$,
inherited from the hypercomplex plane $\sbC$. The elements of
\,$\Go^{\ort}$ play the role of {\it finite\/} or {\it accessible
characters}. The characters $\gama \in \Gf^{\ort}$ are
infinitesimally close to 1 on the whole subgroup $\Gf$ of
finite elements of \,$G$, i.e., in front of the horizon they
are indistinguishable from the trivial character $1_G \in \dG$.
They play the role of {\it infinitesimal characters}. This
intuition, however, calls for some justification.

To this end we introduce the {\it Bohr sets\/} (cf.~\cite{TV})
$$\align
\Bohr_\alfa(A) &= \bigl\{\gama \in \dG\:
(\all a \in A)\bigl(\lv\arg\gama(a)\rv \le \alfa\bigr)\bigr\}\,,\\
\Bohr_\alfa(\Gama) &= \bigl\{a \in G\:
(\all \gama \in \Gama)
   \bigl(\lv\arg\gama(a)\rv \le \alfa\bigr)\bigr\}\,,
\endalign$$
defined for any $\alfa \in \sbR$, $0 \le \alfa \le \pi$, and
$A \sbs G$, \,$\Gama \sbs \dG$. Obviously,
$$\align
\Bohr_\alfa(A) &= \bigl\{\gama \in \dG\: (\all a \in A)
   \bigl(\lv\gama(a) - 1\rv \le 2\sin(\alfa/2)\bigr)\bigr\}\\
&= \bigl\{\gama \in \dG\:
   (\all a \in A)\bigl(\Re\gama(a) \ge \cos\alfa\bigr)\bigr\}\,,
\endalign$$
and similarly for $\Bohr_\alfa(\Gama)$. We also have
$$
\Bohr_\beta(A) \sbs \Bohr_\alfa(A)\,,
\qquad\text{and}\qquad
\Bohr_\beta(\Gama) \sbs \Bohr_\alfa(\Gama)\,,
$$
for $0 \le \beta \le \alfa \le \pi$. For any fixed $\alfa$, the
assignments $A \mto \Bohr_\alfa(A)$, $\Gama \mto \Bohr_\alfa(\Gama)$
satisfy analogous Galois type relations like the infinitesimal
annihilators $A^{\ort}$, $\Gama^{\ort}$.

For $A \sbs G$ and any subset $T$ of the interval
$(0,\pi] \sbs \bR$, such that $\inf T = 0$, we obviously have
$$
A^{\ort} = \bigcap_{\alfa \in T} \Bohr_\alfa(A)\,.
$$
The point is that for a {\it subgroup\/} $A$ of \,$G$ a single
$\alfa$ suffices.

\proclaim{2.1.1.~Lemma}
Let \,$\alfa \in (0,2\pi/3)$. If \,$A$ is a subgroup of \,$G$,
then $A^{\ort} = \Bohr_\alfa(A)$.
\endproclaim

\demo{Proof}
As each $\gama \in \dG$ is a homomorphism, for a subgroup $A \sbs G$
and $\gama \in \Bohr_\alfa(A)$, the image $\gama[A]$ must be a
subgroup of \,$\sbT$ contained in the arc
$\{c \in \sbT\: \lv\arg c\rv \le \alfa\}$. For $\alfa < 2\pi/3$,
however, the biggest subgroup of \,$\sbT$ contained there is namely
the monad \,$\IT = \{c \in \sbT\: c \apr 1\}$ of \,$1 \in \sbT$.
Thus $\gama[A] \sbs \IT$, hence $\Bohr_\alfa(A) \sbs A^{\ort}$, while
the reversed inclusion is trivial.
\enddemo

If \,$A$ is an {\it internal\/} subgroup of \,$G$ and
\,$0 < \alfa < 2\pi/3$, then a similar argument shows that both
the sets coincide with the (strict) annihilator $A^\perp$ of \,$A$,
more precisely,
$$
A^{\ort} = \Bohr_\alfa(A) = \Bohr_0(A) =
\bigl\{\gama \in \dG\: (\all a \in A)(\gama(a) = 1)\bigr\}
= A^\perp\,.
$$

\proclaim{2.1.2.~Proposition}
$\bigl(\dG,\Gf^{\ort},\Go^{\ort}\bigr)$ is a condensing IMG group
triplet with hyperfinite abelian ambient group $\dG$.
\endproclaim

\demo{Proof}
Let us pick any (standard) $\alfa \in (0,2\pi/3)$. As \,$\Go$, $\Gf$ are
subgroups of \,$G$, according to the last Lemma we have
$$\align
\Go^{\ort} &= \Bohr_\alfa(\Go)
= \Bohr_\alfa\Bigl(\,\bigcap \QQ\Bigr)
= \bigcup_{Q \in \QQ} \Bohr_\alfa(Q)\,,\\
\Gf^{\ort} &= \Bohr_\alfa(\Gf)
= \Bohr_\alfa\Bigl(\,\bigcup \QQ\Bigr)
= \bigcap_{Q \in \QQ} \Bohr_\alfa(Q)\,.
\endalign$$
While the last equation in the second line is trivial, the
last equation in the first one follows by a straightforward
{\it saturation\/} argument. As all the sets $\Bohr_\alfa(Q)$ are
internal, $\Go^{\ort}$ is a galactic set and $\Gf^{\ort}$ is a
monadic one. It remains to show that the IMG group triplet
$\bigl(\dG,\Gf^{\ort},\Go^{\ort}\bigr)$ is condensing.

To this end consider any symmetric internal sets $\Gama$, $\Dela$
such that $\Gf^{\ort} \sbs \Gama \sbs \Dela \sbs \Go^{\ort}$. Then
there are $P, Q \in \QQ$ such that $\Go \sbs P \sbs Q \sbs \Gf$
and
$$
\Gf^{\ort} \sbs \Bohr_\alfa(Q) \sbs \Gama \sbs \Dela
\sbs \Bohr_{\alfa/4}(P) \sbs \Go^{\ort}\,.
$$
Using an elementary combinatorial argument we will show that
$$
\lceil\Dela:\Gama\rceil \le
\bigl\lceil\Bohr_{\alfa/4}(P):\Bohr_\alfa(Q)\bigr\rceil \le
\biggl\lceil\frac{4\pi}{\alfa}\biggr\rceil^{\lfloor Q : P\rfloor}
< \infty\,,
$$
where the upper integer part $\lceil4\pi/\alfa\rceil$ equals
the covering index $\lfloor \bT : S\rfloor$ of the
circle $\bT$ with respect to the arc
$S = \{c \in \bT\: \lv\arg c\rv \le \alfa/4\}$. Indeed, let
$k = \lfloor \bT : S\rfloor$ and $q = \lfloor Q : P\rfloor$.
Then there are $k$ points $c_1,\dots,c_k \in \bT$ and $q$ points
$x_1,\dots x_q \in G$ such that $\bT \sbs \bigcup_{i=1}^k Sc_i$ and
$Q \sbs \bigcup_{j=1}^q P+x_i$. Let $h\:\bT \to \{c_1,\dots,c_k\}$
be any function such that $c \in S h(c)$ for $c \in \bT$. As there
are only $k^q$ functions
$\{x_1,\dots,x_q\} \to \{c_1,\dots,c_k\}$, given more then $k^q$
functions in $\Bohr_{\alfa/4}(P)$, there are at least two, $\gama$
and $\chi$, say, such that $(h \co \gama)(x_j) = (h \co \chi)(x_j)$
for $1 \le j \le q$. Let's choose any $x \in Q$ and a $j \le q$ such
that $x \in P + x_j$. Then
$\gama(x)\gama(x_j)^{-1}, \chi(x)\chi(x_j)^{-1} \in S$ and, as
$h(\gama(x_j)) = h(\chi(x_j))$, we have
$$
\gama(x_j)\chi(x_j)^{-1} =
\gama(x_j)h(\gama(x_j))^{-1}h(\chi(x_j))\chi(x_j)^{-1} \in S^2\,.
$$
Consequently,
$$
\gama(x)\chi(x)^{-1} =
\gama(x)\gama(x_j)^{-1}\gama(x_j)\chi(x_j)^{-1}\chi(x_j)\chi(x)^{-1}
\in S^4\,,
$$
i.e., $\gama\chi^{-1} \in \Bohr_\alfa(Q)$. It follows that
$\bigl\lceil\Bohr_{\alfa/4}(P):\Bohr_\alfa(Q)\bigr\rceil \le k^q$.
\enddemo

The condensing IMG group triplet
$\bigl(\dG,\Gf^{\ort},\Go^{\ort}\bigr)$ will be called the
{\it dual triplet\/} of the IMG group triplet $(G,\Go,\Gf)$.

For an internal character $\gama \in \Go^{\ort}$ there are
potentially two interpretations of its observable trace $\gama^\be$.
First, it is simply the element $\Gf^{\ort}\gama$ of the quotient
$\Go^{\ort}\!/\Gf^{\ort} = \dG^\be$. Second, $\gama^\be$ is the
observable trace of the $S$-continuous mapping $\gama\: G \to \sbT$,
i.e.,
$$
\gama^\be\bigl(x^\be\bigr) = {}^\co\gama(x)
$$
for $x \in \Gf$. That way $\gama^\be\:G^\be \to \bT$ is a continuous
character of the LCA group $G^\be$ (cf\. Section~1.3). The assignment
$\gama \mto \gama^\be$, depicted in the commutative diagram
$$\CD
G @<{\Id_{\Gf}}<< \Gf @>{{}^\be}>> G^\be \\
@V{\gama}VV @VV{\gama\,\rst\,\Gf}V @VV{\gama^{\be}}V \\ 
\sbT @>>{\Id_{\sbT}}> \sbT @>>{{}^\co}> \bT
\endCD$$
is a group homomorphism $\Go^{\ort} \to \wh{G^\be}$. Its kernel is
the subgroup $\Gf^{\ort} \sbs \dG$ of all infinitesimal characters
in $\dG$. Thus the assignment $\gama \mto \gama^\be$ induces an
injective group homomorphism $\dG^\be \to \wh{G^\be}$ from the
observable trace $\dG^\be = \Go^{\ort}\!/\Gf^{\ort}$ of the dual
triplet $\bigl(\dG,\Gf^{\ort},\Go^{\ort}\bigr)$ into the dual
group $\wh{G^\be} = \wh{\Gf/\Go}$ of the observable trace
$G^\be = \Gf/\Go$ of the original triplet $(G,\Go,\Gf)$. The canonical
injective homomorphism $\Go^{\ort}\!/\Gf^{\ort} \to \wh{\Gf/\Go}$
justifies the identification of the ``two observable traces''
$\Gf^{\ort}\gama$ and $\gama^\be$. As proved by Gordon in \cite{Go1}
(see also \cite{Go2}), even more is true.

\proclaim{2.1.3.~Proposition}
The canonical mapping \,$\Gf^{\ort}\gama \mto \gama^\be$
is an isomorphism of the topological group
$\dG^\be = \Go^{\ort}\!/\Gf^{\ort}$ onto a closed subgroup
of the topological group $\wh{G^\be} = \wh{\Gf/\Go}$.
\endproclaim

\demo{Proof}
Denote $\Gb = G^\be$ and pick some $0 < \alfa < 2\pi/3$. On one hand,
the images of the Bohr sets $\Bohr_\alfa(Q)$, where $Q \in \QQ$,
under the quotient mapping $\Go^{\ort} \to \Go^{\ort}\!/\Gf^{\ort}$
form a neighborhood base of the unit character in $\dG^\be$.
On the other hand, the Bohr sets
$$
\Bohr_\alfa\bigl(Q^\be\bigr) =
   \bigl\{\ggam \in \dGb\: \bigl(\all \xb \in Q^\be\bigr)
      \bigl(\lv\arg\ggam(\xb)\rv \le \alfa\bigr)\bigr\}\,,
$$
where $Q \in \QQ$, form a neighborhood base of the unit
character in $\dGb$. It follows that the canonical injective
group homomorphism $\Go^{\ort}\!/\Gf^{\ort} \to \dGb$ is also
a homeomorphism of \,$\Go^{\ort}\!/\Gf^{\ort}$ onto the subgroup
$\bigl\{\gama^\be\: \gama \in \Go^{\ort}\bigr\}$ of \,$\dGb$.
As a continuous image of an LCA (hence complete) topological
group, it is necessarily closed.
\enddemo

It is both natural and tempting to conjecture that the canonical
mapping $\dG^\be \to \wh{G^\be}$ is also surjective, i.e., that
it is an isomorphism of topological groups. This is indeed the
first of Gordon's Conjectures from \cite{Go1} (see also
\cite{Go2, page~132}).

\proclaim{2.1.4.~Theorem}{\bf [Gordon's Conjecture~1]}
Let $(G,\Go,\Gf)$ be a condensing IMG group triplet with
hyperfinite abelian ambient group $G$. Then the canonical
homomorphism $\Go^{\ort}\!/\Gf^{\ort} \to \wh{\Gf/\Go}$
is an isomorphism of topological groups.
\endproclaim

The proof of Gordon's Conjecture~1 is the first main result of the
present paper. In view of 2.1.3, this amounts just to show that
every continuous character $\ggam$ of the LCA group
$G^\be = \Gf/\Go$ is indeed of the form $\ggam = \gama^\be$
for some internal $S$-continuous character
$\gama \in \Go^{\ort}$. However, we will approach the proof of the
above Theorem indirectly, by investigating the dual triplet of the
dual triplet of the original IMG group triplet $(G,\Go,\Gf)$.

The second dual \,$\ddG$ of the hyperfinite abelian group $G$ can be
naturally identified with the original group $G$. Then the second
dual of the original triplet $(G,\Go,\Gf)$ is defined as the
condensing IMG triplet
$\bigl(G,\Go^{\ort\ort},\Gf^{\ort\ort}\bigr)$.

\proclaim{2.1.5.~Triplet Duality Theorem}
Let $(G,\Go,\Gf)$ be a condensing IMG group triplet with
hyperfinite abelian ambient group $G$. Then
$$
\Go^{\ort\ort} = \Go\,,
\qquad\text{and}\qquad
\Gf^{\ort\ort} = \Gf\,;
$$
in other words, the dual triplet
$\bigl(G,\Go^{\ort\ort},\Gf^{\ort\ort}\bigr)$ of the
dual triplet $\bigl(\dG,\Gf^{\ort},\Go^{\ort}\bigr)$ equals
the original group triplet $(G,\Go,\Gf)$.
\endproclaim

The proof of Theorem~2.1.5 is postponed into the next two sections.
At this place we will just show how Gordon's Conjecture~1, i.e.,
Theorem~2.1.4, can be derived from Theorem~2.1.5. The proof of this
implication relies on the following consequence of the
Pontryagin-van\,Kampen duality theorem.

\proclaim{2.1.6.~Proposition}
Let \,$\Gb$ be a Hausdorff LCA group and \,$\Db$ be a subgroup of its
dual group $\dGb$. Then $\Db$ separates points in $\Gb$ if and only
if it is dense in $\dGb$.
\endproclaim

\demo{Sketch of proof}
$\dGb$ obviously separates points in its dual group $\ddGb$ which, by
the virtue of the duality theorem, can be identified with $\Gb$. Now,
it can be easily seen that any dense subgroup $\Db \sbs \dGb$ separates
points in $\Gb$, as well.

Conversely, if \,$\Db$ is not dense in $\dGb$, then its closure $\oDb$
is a proper closed subgroup of \,$\dGb$. As a consequence of the duality
theorem (cf.~\cite{Mo, Corollary~1 to Theorem~21}), there is a nontrivial
character of \,$\dGb$, i.e., an element $\xb \in \Gb \sms \{0\}$, such
that $\ggam(\xb) = 1$ for all $\ggam \in \oDb$. Therefore, neither $\oDb$
nor $\Db$ separate points in $\Gb$.
\enddemo

\demo{Proof of \,{\rm 2.1.5$\imp$2.1.4}}
By Proposition 2.1.3, the observable trace
$\dG^\be = \Go^{\ort}\!/\Gf^{\ort}$ of the dual triplet
$\bigl(\dG,\Gf^{\ort},\Go^{\ort}\bigr)$ can be identified with the
closed subgroup of the dual group $\wh{G^\be}$ formed by the
observable traces $\gama^\be$ of the internal characters
$\gama \in \Go^{\ort}$. Therefore, in order to show that
$\dG^\be = \wh{G^\be}$, it suffices to prove that $\dG^\be$ is dense
in $\wh{G^\be}$. By Proposition~2.1.6, to this end it is enough to
show that $\dG^\be$ separates points in $G^\be$. But otherwise there
would be an $x \in \Gf \sms \Go$, such that $\gama(x) \apr 1$ for
every $\gama \in \Go^{\ort}$. Then
$$
x \in \Gf \cap \Go^{\ort\ort} = \Gf \cap \Go = \Go\,,
$$
which is a contradiction.
\enddemo

\remark{Remark}
Notice that, in order to derive 2.1.4, just a weaker version
of the first equality in 2.1.5 would be sufficient, namely
$\Go = \Go^{\ort\ort} \cap \Gf$, which, of course, trivially
follows from $\Go^{\ort\ort} = \Go$ and $\Go \sbs \Gf$. The second
equality $\Gf^{\ort\ort} = \Gf$ is not needed to this end. Moreover,
it is already a consequence of the first one.
\endremark

\proclaim{2.1.7.~Lemma}
Let $(G,\Go,\Gf)$ be as above. Then
$$
\Gf^{\ort\ort} = \Gf + \Go^{\ort\ort}\,.
$$
Therefore, $\Go^{\ort\ort} = \Go$ implies $\Gf^{\ort\ort} = \Gf$.
\endproclaim

\demo{Proof}
Let us denote $(H,\Ho,\Hf) = \bigl(\dG,\Gf^{\ort},\Go^{\ort}\bigr)$
the dual triplet of $(G,\Go,\Gf)$. According to the properties of
Galois connections, three $\Ort$'s reduce to one, hence the
group triplet $(H,\Ho,\Hf)$ satisfies the conditions
$$
\Ho^{\ort\ort} = \Ho\,,
\qquad\text{and}\qquad
\Hf^{\ort\ort} = \Hf\,,
$$
so that the dual triplet of $\bigl(\dH,\Hf^{\ort},\Ho^{\ort}\bigr)$
is indeed $(H,\Ho,\Hf)$ and the canonical mapping
$\Ho^{\ort}\!/\Hf^{\ort} \to \wh{\Hf/\Ho}$
is an isomorphism of topological groups.

Now, consider the canonical embedding
$\Go^{\ort}\!/\Gf^{\ort} \to \wh{\Gf/\Go}$ and apply the (external)
duality functor to it. Identifying the dual of the LCA group
$\wh{\Gf/\Go}$ with $\Gf/\Go$, and the dual of the LCA group
$\Hf/\Ho = \Go^{\ort}\!/\Gf^{\ort}$ with
$\Ho^{\ort}\!/\Hf^{\ort} = \Gf^{\ort\ort}\!/\Go^{\ort\ort}$,
we obtain a surjective continuous homomorphism
$$
\Gf/\Go \to \Gf^{\ort\ort}\!/\Go^{\ort\ort}\,,
$$
sending $\gama + \Go \in \Gf/\Go$ to
$\gama + \Go^{\ort\ort} \in \Gf^{\ort\ort}\!/\Go^{\ort\ort}$.
Its surjectivity simply means that $\Gf$, $\Go^{\ort\ort}$,
$\Gf^{\ort\ort}$, as subgroups of \,$G$, satisfy
$\Gf^{\ort\ort} = \Gf + \Go^{\ort\ort}$.
\enddemo

The following Example shows that (a) the condition $X^{\ort\ort} = X$
is not automatically satisfied for monadic or galactic subgroups $X$
of hyperfinite abelian groups, (b) there are plenty of subgroups of
hyperfinite abelian groups satisfying this condition even without
being monadic or galactic.

\example{2.1.8.~Example}
An {\it additive cut\/} on $\sbN$ is any nonempty subset $C \sbs \sbN$,
such that
$$\gather
(\all a \in \sbN)(\all b \in C)(a \le b \imp a \in C)\,,\\
(\all a,b \in C)(a+b \in C)\,.
\endgather$$
Assume that $C$ is an additive cut on $\sbN$, $\bN \sbs C$,
and $n \in \sbN \sms C$.
\smallskip
(a) Let $Z$ be any nontrivial finite abelian group, e.g.,
$Z = \bZ_d$ for some $2 \le d \in \bN$. Let $G = Z^n$ be
the hyperfinite abelian group of all internal sequences
$z =(z_1,\dots,z_n)$ of elements from $Z$. Denote
$$
\supp z = \{k\: 1 \le k \le n \et z_k \ne 0\}\,,
$$
for $z \in G$, and put
$$
S(C) = \{z \in G\: \lv\supp z\rv \in C\}\,.
$$
As $C \ne \{0\}$ and $n \nin C$, as well as $\supp(-z) = \supp z$
and $\supp(y + z) \sbs \supp y \cup \supp z$ for any $y,z \in G$,
$S(C)$ is a nontrivial proper subgroup of \,$G$. Moreover, if \,$C$
is a monadic (galactic) set, then so is $S(C)$. On the other hand,
all the internal subgroups
$$
G_k = \bigl\{z \in G\: \supp z \sbs \{k\}\bigr\}
= \bigl\{z \in G\: (\all i \ne k)(z_i = 0)\bigr\}\iso Z\,,
$$
where $1 \le k \le n$, satisfy $G_k \sbs S(C)$. Identifying
the internal dual $\dG$ with the hyperfinite abelian group
$\wh{Z}^n$ in the obvious way, each character $\gama \in \dG$
is represented as the ordered $n$-tuple
$\gama = (\gama_1,\dots,\gama_n)$, where $\gama_k \in \wh{Z}$.
Then
$$
G_k^{\ort} = \bigl\{\gama \in \dG\: \gama_k = 1_Z\bigr\}\,,
$$
hence,
$$
S(C)^{\ort} \sbs \bigcap_{k=1}^n G_k^{\ort} = \{1_G\}\,,
$$
and, finally, $S(C)^{\ort\ort} = G$.

Now, if we take a monadic additive cut $A$ and a galactic additive
cut $B$, such that $\bN \sbs A \sbs B \sbs \sbN$ and $n \nin B$,
then $S(A)$ is a proper subgroup of \,$S(B)$ and
$\bigl(G,S(A),S(B)\bigr)$ is an (of course, non-condensing) IMG
group triplet with hyperfinite abelian ambient group $G = Z^n$ and
nontrivial observable trace $S(B)/S(A)$. However, its first and
second dual triplets are $\bigl(\dG,\{1_G\},\{1_G\}\bigr)$ and
$(G,G,G)$, respectively; both have trivial observable traces.

\smallskip
(b) Denote $A = C \cup (-C) = \{a \in \sbZ\: \lv a\rv \in C\}$.
Then $A$ is a subgroup of the hyperfinite cyclic group
$\bZ_{2n+1} = \{0,\pm 1,\dots, \pm n\}$ modulo $2n+1$. We leave
the reader as an exercise to verify that $A^{\ort\ort} = A$,
regardless of any further properties of \,$C$.
\endexample

\specialhead
2.2.~Fourier transforms, Bohr sets and spectral sets
     in finite abelian groups
\endspecialhead
\flushpar
In this section we make a necessary digression in order to establish
some inclusions between certain types of (internal) subsets of
(hyper)finite abelian groups and some estimates of their size.
However, in view of the {\it transfer principle}, it is sufficient
to deal with the finite case, only. To this end we will make use
of the discrete Fourier transform. The results thus obtained will
not depend on the scaling coefficients occurring in it. Just as a
matter of convenience we choose $d = 1/|G|$ and $\hat d = 1$.

In what follows $G$ denotes a finite abelian group with the addition
operation. The set $\bC^G$ of all functions $G \to \bC$ becomes
a unitary space endowed with the Hermitian inner product, given as
the expectation
$$
\br f, g\kt_G = \Ee(f \cd \ovl g)
= \Ee_{x \in G}\,f(x)\,\ovl g(x)
= {1 \over |G|} \sum_{x \in G} f(x)\,\ovl{g(x)}\,,
$$
for $f,g \in \bC^G$, and the corresponding $\Lb^2$-norm
$$
\|f\|_{_2} = \sqrt{\br f, f\kt_G} = \sqrt{\Ee(f \cd \ovl f)}\,.
$$
The dual group $\dG = \Hom(G,\bT)$, considered as a subset of
\,$\bC^G$, forms an orthonormal basis of \,$\bC^G$. Thus for the
Fourier transform \ $\FF\:\bC^G \to \bC^{\dG}$ 
$$
\FF(f)(\gama) = \wh f(\gama) = \br f, \gama\kt_G = \Ee(f \cd \ovl\gama)
$$
the Fourier inversion formula takes the form
$$
f = \sum_{\gama \in \dG} \wh f(\gama)\,\gama\,.
$$
The Hermitian inner product on $\bC^{\dG}$
$$
\br \vfi,\psi\kt_{\dG} =
\sum_{\gama \in \dG} \vfi(\gama)\,\ovl\psi(\gama)
$$
ensures the Plancherel identity
$$
\br f,g\kt_G = \bigl\br\wh f,\wh g\,\bigr\kt_{\dG}\,.
$$
The convolution on $\bC^G$ is also defined as the expectation
$$
(f * g)(x) = \Ee_{y \in G}\,f(x-y)\,g(y)\,;
$$
for its Fourier transform we have
$$
\wh{f * g} = \wh f \cd \wh g\,.
$$
Additionally, we will make use of the $\Lbp$-norms on $\bC^G$
and the $\ell^p$-norms on $\bC^{\dG}$
$$\gather
\|f\|_{_p} = \Bigl(\Ee\lv f\rv^p\Bigr)^{1/p}
= \Bigl(\Ee_{x \in G} |f(x)|^p\Bigr)^{1/p}
= \biggl({1\over|G|}\sum_{x\in G}|f(x)|^p\biggr)^{1/ p}\,,\\
\|\vfi\|_{_p} =
\Biggl(\ \sum_{\chi\in\dG} |\vfi(\chi)|^p\Biggr)^{1/p}\,,
\endgather$$
for $1 \le p < \infty$, as well as of the $\Lb^{\infty}$-norm on
$\bC^G$ and the $\ell^{\infty}$-norm on $\bC^{\dG}$
$$
\|f\|_{_\infty} = \max_{x \in G} |f(x)|\,
\qquad\qquad
\|\vfi\|_{_\infty} = \max_{\chi \in \dG} |\vfi(\chi)|\,,
$$
for $f \in \bC^G$, $\vfi \in \bC^{\dG}$.

Let us list some well known and/or obvious relations:
$$
\bigl\|\wh f\,\bigr\|_{_2} = \|f\|_{_2}\,,
\qquad\qquad
\bigl\|\wh f\,\bigr\|_{_\infty} \le \|f\|_{_1}\,,
$$
which are special cases of the Hausdorff-Young inequality
$$
\bigl\|\wh f\,\bigr\|_{_q} \le \|f\|_{_p},
$$
for $1 \le p \le 2$ and $q$ being the dual exponent of $p$,
i.e., $1/p + 1/q = 1$.

We denote by $f_a$ the shift of the function $f\:G \to \bC$ by the
element $a \in G$, i.e.,
$$
f_a(x) = f(x-a)
$$
for $x \in G$. Then
$$
\wh{f_{a\,}}(\gama) = \ovl\gama(a)\, \wh f(\gama)
\qquad\text{and}\qquad
\wh f_{\,\gama} = \wh{\gama f}
$$
for  $\gama \in \dG$, as well as
$$
f * g = \Ee_{a \in G}\,f_a\,g(a) = \Ee_{a \in G}\,f(a)\,g_a
$$
and
$$
(f * g)_a = f * g_a = f_a * g\,,
$$
for $f,g \in \bC^G$.

A norm $\Nn$ on the linear space $\bC^G$ is called {\it translation
invariant} if \,$\Nn(f_a) = \Nn(f)$ for all $f \in \bC^G$, $a \in G$.
For any translation invariant norm $\Nn$ we have
$$
\Nn(f * g) \le \Nn(f)\,\|g\|_{_1}\,.
$$
A norm $\Nn$ on $\bC^G$ is called {\it absolute} if
\,$\Nn(f) \le \Nn(g)$ for any functions $f,g \in \bC^G$ such that
$|f(x)| \le |g(x)|$ for all $x \in G$.
For any absolute norm $\Nn$ we have
$$
\Nn(f\,g) \le \Nn(f)\,\|g\|_{_\infty}\,.
$$
As all the norms $\|\cd\|_{_p}$ are obviously translation invariant
and absolute,
$$
\|f * g\|_{_p} \le \|f\|_{_p}\,\|g\|_{_1}\,,
\qquad\qquad
\|f \,g\|_{_p} \le \|f\|_{_p}\,\|g\|_{_\infty}\,
$$
for $f,g\: G \to \bC$, $1 \le p \le \infty$.
If \,$f,g\: G \to \bR$ are both nonnegative on $G$, then even
$$
\|f * g\|_{_1} = \|f\|_{_1}\,\|g\|_{_1}\,.
$$

For the characteristic function (indicator) $1_A$ of a set $A \sbs G$
and $1 \le p < \infty$ we have
$$
\|1_A\|_{_p}^p = \|1_A\|_{_2}^2 = \|1_A\|_{_1} = {|A| \over |G|}\,.
$$

If $f\: G \to \bC$ is even, i.e., $f(-x) = f(x)$ for $x \in G$, then
so is $\wh f\:\dG \to \bC$ and we have
$$
\wh f(\gama)= \Ee_{x \in G}\,f(x) \Re\gama(x)\,,
\qquad\qquad
f(x) = \sum_{\gama \in \dG} \wh f(\gama) \Re\gama(x)\,.
$$

The {\it support} of a function $g\:G \to \bC$ is defined to
be the set
$$
\supp g = \{x \in G\: g(x) \ne 0\}\,.
$$

The following estimate generalizes an inequality in \cite{GR},
which was part of the proof of Proposition~3.1 there, from
indicators to arbitrary nonnegative functions.

\proclaim{2.2.1~Lemma}
Let $f\:G \to \bR$ be nonnegative, and $D \sbs G$ be a nonempty
set such that $\supp(f * f) \sbs D$. Then
$$
\sum_{\gama \in \dG} \bigl|\wh f(\gama)\bigr|^4
   \ge {\|f\|_{_1}^4 \over \|1_D\|_{_1}}\,.
$$
\endproclaim

\demo{Proof}
By Plancherel identity and the relation between the Fourier
transform and convolution,
$$
\sum_{\gama \in \dG} \bigl|\wh f(\gama)\bigr|^4 =
\bigl\br{\wh f}\,^{\,2},{\wh f}\,^{\,2}\bigr\kt_{\dG}
= \br f*f, f*f\kt_G = \|f * f\|_{_2}^2
$$
for any function $f\:G \to \bC$. Using Cauchy-Schwartz inequality,
the fact that $f * f$ is supported on $D$, and nonnegativity of $f$
we get
$$
\|f * f\|_{_2}^2\,\|1_D\|_{_2}^2 \ge \bigl|\br f*f, 1_D\kt_G\bigr|^2
= \|f * f\|_{_1}^2 = \|f\|_{_1}^4\,.
$$
The claim immediately follows.
\enddemo

The {\it spectral set\/} or {\it spectrum at threshold\/} $t \in \bR$
of a function $f\: G \to \bC$ is the set
$$
\Spec_t(f) = \bigl\{\gama \in \dG\:
   \bigl|\wh f(\gama)\bigr| \ge t\,\|f\|_{_1}\bigr\} \sbs \dG\,.
$$
This is a slight generalization of a definition from \cite{TV, \S 4.6}
were spectral sets of subsets $A \sbs G$ were defined by
$$
\Spec_t(A) = \Spec_t(1_A)
= \bigl\{\gama \in \dG\:
\bigl|\wh{1_A}(\gama)\bigr| \ge t\,\|1_A\|_{_1}\bigr\}\,.
$$
As $\bigl|\wh f(\gama)\bigr| \le
\bigl\|\wh f\,\bigr\|_{_\infty} \le \|f\|_{_1}$,
\,$\Spec_t(f) = \emptyset$ whenever $t > 1$ and $f$ is not
identically~$0$; similarly, $\Spec_t(f) = \dG$ \,for $t \le 0$.
Thus it makes sense to consider the spectral sets just for the
threshold values $0 \le t \le 1$. Also the following implication
is trivial
$$
s \le t \imp \Spec_t(f) \sbs \Spec_s(f)\,.
$$

Let us recall that the {\it Bohr sets\/} $\Bohr_\alfa(A)$,
$\Bohr_\alfa(\Gama)$, for $A \sbs G$, $\Gama \sbs \dG$, were
defined in Section~2.1.

\proclaim{2.2.2.~Lemma}
Let \,$f\:G \to \bR$ be an even nonnegative function,
$D \sbs G$ be a nonempty set such that \,$\supp f \sbs D$,
and \,$0 \le \alfa \le \pi/2$. Then
$$
\Bohr_\alfa(D) \sbs \Spec_t(f)\,,
$$
whenever $0 \le t \le \cos\alfa$.
\endproclaim

\demo{Proof}
Take any $\gama \in \Bohr_\alfa(D)$. According to the assumptions on
$f$, $t$ and $\alfa$ we have
$$\align
\bigl|\wh f(\gama)\bigr| &=
\lv\Ee_{x \in G}\,f(x)\,\ovl\gama(x)\rv
= \lv\Ee_{x \in G}\,1_D(x)\,f(x) \Re\gama(x)\rv \\
&\ge \Ee_{x \in G}\,f(x)\cos\alfa
= \|f\|_{_1} \cos \alfa \ge t\,\|f\|_{_1}\,,
\endalign$$
thus $\gama \in \Spec_t(f)$.
\enddemo

The following result generalizes an inclusion proved in \cite{GR}
during the proof of Proposition~3.1, as well.

\proclaim{2.2.3.~Proposition}
Let \,$f\:G \to \bR$ be a nonnegative function, not identically equal
to~$0$, $D \sbs G$ be a nonempty set such that \,$\supp(f * f) \sbs D$,
and \,$0 < \alfa < \pi/2$, $0 \le t \le  1$. Then
$$
\Bohr_{\alfa}\bigl(\Spec_t(f)\bigr) \sbs D - D\,,
$$
whenever
$$
t \le {\|f\|_{_1} \over \|f\|_{_2} \,\|1_D\|_{_2}}
   \cd \sqrt{\cos\alfa \over 1+\cos\alfa}\,.
$$
\endproclaim

\demo{Proof}
Let $f_{_-}$ denote the function given by $f_{_-}(x) = f(-x)$;
since $f\: G \to \bR$, we have
$\wh{f_{_-}}(\gama) = \ovl{\wh f(\gama)}$ and
$\bigl|\wh f(\gama)\bigr|^2 = \wh{f * f_{_-}}(\gama)$.
Then both the functions $f*f_{_-}$, \,$f*f*f_{_-}*f_{_-}$ are even,
and
$$
\supp(f*f*f_{_-}*f_{_-}) \sbs D - D\,.
$$
Let $x \in \Bohr_\alfa\bigl(\Spec_t(f)\bigr)$. It suffices to prove
that $(f*f*f_{_-}*f_{_-})(x) > 0$. Putting $\Gama = \Spec_t(f)$, we
have $\Re\gama(x) \ge \cos\alfa$ for $\gama \in \Gama$. By the
Fourier inversion formula, Lemma~2.2.1 and Plancherel identity
we get
$$\align
(f * f * f_{_-} * f_{_-})(x) &=
\sum_{\gama \in \dG} \bigl|\wh f(\gama)\bigr|^4\,\gama(x) \\
&= \sum_{\gama \in \Gama} \bigl|\wh f(\gama)\bigr|^4 \Re\gama(x) +
\sum_{\gama\in\dG\sms\Gama} \bigl|\wh f(\gama)\bigr|^4 \Re\gama(x)\\
&> \sum_{\gama \in \Gama} \bigl|\wh f(\gama)\bigr|^4 \cos\alfa
  - \sum_{\gama \in \dG \sms \Gama} \bigl|\wh f(\gama)\bigr|^4 \\
&= \sum_{\gama \in \dG} \bigl|\wh f(\gama)\bigr|^4 \cos\alfa -
\sum_{\gama\in\dG\sms\Gama} \bigl|\wh f(\gama)\bigr|^4 (1+\cos\alfa)\\
&\ge {\|f\|_{_1}^4 \over \|1_D\|_{_1}} \cos\alfa
- \max_{\gama \in \dG \sms \Gama} \bigl|\wh f(\gama)\bigr|^2
  \sum_{\gama \in \dG}\bigl|\wh f(\gama)\bigr|^2 (1 + \cos\alfa) \\
&\ge {\|f\|_{_1}^4 \over \|1_D\|_{_1}} \cos\alfa
   - t^2 \|f\|_{_1}^2\,\|f\|_{_2}^2 (1 + \cos\alfa) \\
&=  \|f\|_{_1}^2\biggl({\|f\|_{_1}^2 \over \|1_D\|_{_2}^2} \cos\alfa
   - t^2 \|f\|_{_2}^2 (1 + \cos\alfa)\biggr)\,.
\endalign$$
The strict inequality in the third line is due to the fact that for the
trivial character $1_G \in \Gama$ we have $\Re 1_G(x) = 1 > \cos\alfa$ as
$\alfa > 0$. According to the assumption on $t$, the expression in the last
line is $\ge 0$.
\enddemo

\remark{Remark}
(a) It would be enough to have the above result for a single fixed value
$\alfa$. The authors in \cite{GR}, \cite{TV} take $\alfa = \pi/3$, which,
in our case, would result in the estimate
$$
t \le {\|f\|_{_1} \over \|f\|_{_2} \,\|1_D\|_{_2} \sqrt{3}}\,.
$$

(b) It is worthwhile to notice that in the special case of $f$ being the
indicator of a nonempty set $A \sbs G$ and $D = A + A$ we have
$$
{\|f\|_{_1} \over \|f\|_{_2} \,\|1_D\|_{_2}} =
   {\|1_A\|_{_1} \over \|1_A\|_{_2} \,\|1_{A+A}\|_{_2}} =
   \sqrt{|A| \over |A + A|} = {1 \over \sqrt{\sig(A)}}\,,
$$
where $\sig(A) = |A+A|/|A|$ is the {\it doubling constant\/} of $A$.
\endremark
\medskip

We will need also some lower and upper bounds of the size of spectral sets
of some functions.

\proclaim{2.2.4.~Proposition}
Let \,$f\:G \to \bR$ be a nonnegative function, not identically
equal to~$0$, $D \sbs G$ be a nonempty set such that
\,$\supp(f * f) \sbs D$, and \,$0 < t \le  1$. Then
$$
{|G| \over |D|} - t^2\,{\|f\|_{_2}^2  \over \|f\|_{_1}^2}
   \le |\Spec_t(f)|
   \le {1 \over t^2}\,{\|f\|_{_2}^2  \over \|f\|_{_1}^2}\,.
$$
\endproclaim

\demo{Proof}
Let us denote $\Gama = \Spec_t(f)$. According to Lemma~2.2.1 we have
$$\align
{\|f\|_{_1}^4 \over \|1_D\|_{_1}}
  \le \sum_{\gama \in \dG} \bigl|\wh f(\gama)\bigr|^4
&= \sum_{\gama\in\Gama}\bigl|\wh f(\gama)\bigr|^4
 + \sum_{\gama\in\dG\sms\Gama}\bigl|\wh f(\gama)\bigr|^4 \\
&\le |\Gama|\,\bigl\|\wh f\,\bigr\|_{_\infty}^4  +
t^2 \|f\|_{_1}^2 \sum_{\gama\in\dG\sms\Gama} \bigl|\wh f(\gama)\bigr|^2\\
&\le |\Gama|\,\|f\|_{_1}^4 + t^2 \|f\|_{_1}^2\,\|f\|_{_2}^2\,,
\endalign$$
using Plancherel formula to pass to the last line. This gives the
lower bound.

The upper bound readily follows from the following computation:
$$
\|f\|_{_2}^2 = \bigl\|\wh f\,\bigr\|_{_2}^2
= \sum_{\gama \in \dG} \bigl|\wh f(\gama)\bigr|^2
\ge \sum_{\gama\in\Gama}
    \bigl|\wh f(\gama)\bigr|^2 \ge t^2 \|f\|_{_1}^2\,|\Gama|\,.
$$
\enddemo

\remark{Remark}
Notice that the lower bound is relevant just in case
$$
t < {\|f\|_{_1} \over \|f\|_{_2} \,\|1_D\|_{_2}}\,,
$$
otherwise it is trivial. One of its consequences can be stated as
$$
\frac{\lv\supp(f * f)\rv \lv\Spec_t(f)\rv}{\lv G\rv}
\ge 1 - t^2\,\frac{\lv\supp(f * f)\rv}{\lv G\rv}\,
              \frac{\LV f\RV_{_2}^2}{\LV f\RV_{_1}^2}
\ge 1 - t^2\,\frac{\LV f\RV_{_2}^2}{\LV f\RV_{_1}^2}\,,
$$
in which form it can be regarded as a kind of the
{\it Uncertainty Principle}. For $f = 1_A$ being the indicator of
a nonempty set $A \sbs G$ the first inequality gives
$$
\frac{\lv A + A\rv\lv\Spec_t(A)\rv}{\lv G\rv}
\ge 1 - t^2 \sig(A)\,.
$$
\endremark
\medskip

Let us close this technical section with a kind of the
{\it Smoothness-and-Decay Principle\/} indicating that spectral sets
of functions continuous in some sense tend to avoid discontinuous
characters. We will return to this topic in Section~3.2.

\proclaim{2.2.5.~Proposition}
Let \,$f\: G \to \bC$ be not identically equal to $0$ and
\,$C, D \sbs G$ satisfy
\,$ \supp f \cup (\supp f + C) \sbs D \ne \emptyset$.
Let \,$0 < t \le 1$, $0 < \alfa < \pi$, and $\eps > 0$
be real numbers. Assume that
$$
\eps \le 2t\,
{\|f\|_{_1} \over \|1_D\|_{_1}}\,\sin{\alfa \over 2}\,,
$$
and \,$\|f_a - f\|_{_\infty} \le \eps$ for $a \in C$. Then
\,$\Spec_t(f) \sbs \Bohr_\alfa(C)$.
\endproclaim

\demo{Proof}
If $a \in C$, then
$\supp(f_a - f) \sbs (\supp f + C) \cup \supp f \sbs D$, hence
$$
\bigl\|\wh{f_{a\,}} - \wh f\,\bigr\|_{_\infty} \le \|f_a - f\|_{_1}
= \Ee_{x \in G} \bigl|f_a(x) - f(x)\bigr| \le \eps\,\|1_D\|_{_1}\,.
$$
Now, take any $\gama \in \Spec_t(f)$ and assume that
$\gama \nin \Bohr_\alfa(C)$. Then there is an $a \in C$ such that
$|\gama(a) - 1| > 2\sin(\alfa/2)$. Thus we have
$$
\eps\,\|1_D\|_{_1} \ge \bigl\|\wh{f_{a\,}} - \wh f\,\bigr\|_{_\infty}
\ge \bigl|\wh{f_{a\,}}(\gama) - \wh f(\gama)\bigr|
= \lv\,\ovl\gama(a) - 1\rv\bigl|\wh f(\gama)\bigr|
> 2t\,\|f\|_{_1} {\tsize\sin{\alfa \over 2}}\,,
$$
contradicting the assumed upper bound for $\eps$.
\enddemo

A brief inspection of the proof yields the following modification
of the last result.

\proclaim{2.2.6.~Corollary}
Let \,$f\: G \to \bC$ be not identically equal to \,$0$, \,$C\sbs G$,
and \,$0 < t \le 1$, $0 < \alfa < \pi$, \,$\eps > 0$ be real
numbers. Assume that \,$\eps \le 2t \sin{\alfa \over 2}$, and
\,$\|f_a - f\|_{_1} \le \eps\LV f\RV_{_1}$ for $a \in C$.
Then \,$\Spec_t(f) \sbs \Bohr_\alfa(C)$.
\endproclaim

\specialhead
2.3.~The dual triplet continued: Normal multipliers and \\
\phantom{2.3.}~proofs of Gordon's Conjectures 1 and 2
\endspecialhead
\flushpar
Now, we are back dealing with a condensing IMG group triplet
$(G,\Go,\Gf)$ with a hyperfinite abelian ambient group~$G$.
Let us recall that a {\it normalizing multiplier\/} of the
triplet is any positive hyperreal number $d$ such that
$0 \not\apr d\lv A\rv < \infty$ for every (or, equivalently, for
some) internal set $A$ between $\Go$, $\Gf$. Then $\mbd$ denotes
the Haar measure on the observable trace $G^\be = \Gf/\Go$, obtained
by pushing down the Loeb measure $\lam_d$ from $G$ to $G^\be$
(see Sections 1.4 and 1.5). The Hermitian inner product, Fourier
transform, $\Lbp$-norms and convolution on the space $\sbC^G$ of all
internal functions $G \to \sbC$ are normalized by the scaling
coefficient~$d$, i.e.,
$$\gather
\br f,g\kt = d \,\sum_{x \in G} f(x)\,\ovl g(x)\,,\\
\FF(f)(\gama) = \wh f(\gama) = \br f,\gama\kt
              = d \,\sum_{x \in G} f(x)\,\ovl\gama(x)\,,\\
\|f\|_{_p} = \biggl(d \,\sum_{x \in G} |f(x)|^p\biggr)^{1/p}\,,\\
(f * g)(x) = d \,\sum_{a \in G} f(x-a)\,g(a)
           = d \,\sum_{a \in G} f_a(x)\,g(a)\,,
\endgather$$
for $f,g \in \sbC^G$, $\gama \in \dG$, $1 \le p < \infty$, $x \in G$.

In the proof of Triplet Duality Theorem~2.1.5 we will need the
following preliminary qualitative version of the
{\it Smoothness-and-Decay Principle}.

\proclaim{2.3.1.~Proposition}
Let $f\: G \to \sbC$ be an $S$-continuous internal function such
that \,$\supp f \sbs \Gf$. Then \,$\wh f(\gama) \apr 0$ for each
$\gama \in \dG \sms \Go^{\ort}$.
\endproclaim

\demo{Proof}
We will proceed in a similar way as in the proof of
Proposition~2.2.5.

By the $S$-continuity of $f$, $f_a(x) \apr f(x)$ for each $x \in G$,
whenever $a \in \Go$; moreover $f(x) = f_a(x) = 0$ once
$x \in G \sms D$ for some (in fact for each) internal set $D$ such
that $\supp f  + \Go \sbs D \sbs \Gf$. Then
$$
\bigl\|\wh{f_{a\,}} - \wh f\,\bigr\|_{_\infty} \le \|f_a - f\|_{_1}
= d \,\sum_{x \in D} |f_a(x) - f(x)|
\le d\lv D\rv \max_{x \in D} |f_a(x) - f(x)|
\apr 0\,,
$$
as $d\lv D \rv  < \infty$.  Now, for any
$\gama \in \dG \sms \Go^{\ort}$ there is some $a \in \Go$ such that
$\gama(a) \not\apr 1$. Then
$$
0 \apr \bigl\|\wh{f_{a\,}} - \wh f\,\bigr\|_{_\infty} =
\max_{\chi \in \dG} \lv\,\ovl\chi(a) - 1\rv\,\bigl|\wh f(\chi)\bigr|
\ge |\gama(a) - 1|\,\bigl|\wh f(\gama)\bigr|\,.
$$
As $\gama(a) - 1 \not\apr 0$, the conclusion $\wh f(\gama) \apr 0$
follows immediately.
\enddemo

\proclaim{2.3.2.~Corollary}
Let \,$f\: G \to \sbC$ be an $S$-continuous internal function such
that \,$\supp(f) \sbs \Gf$ and \,$\|f\|_{_1} \not\apr 0$. Assume
that \,$t \in \sbR$ and \,$0 < t \le 1$, $t\not\apr 0$. Then
\,$\Spec_t(f) \sbs \Go^{\ort}$.
\endproclaim

\remark{Remark}
Both Proposition~2.3.1 and its Corollary~2.3.2 are ``soft''
statements in the sense of Tao's blogs \cite{Ta1}, \cite{Ta2}.
Their ``hard'' version is Proposition~2.2.5, providing us with the
additional information how one has to choose $\eps > 0$ and the
internal set $C$, such that $\Go \sbs C \sbs \Gf$, in order to get
$\Spec_t(f) \sbs \Bohr_\alfa(C) \sbs \Go^{\ort}$. On the other hand,
the last Corollary is fully sufficient for our purpose. We hope that
at least some readers will appreciate the advantage of nonstandard
arguments, allowing one to dispense with lots of meticulous
estimates and ``epsilontics''.
\endremark
\medskip

As made clear in Section 2.2, in order to finish the proof of
the Triplet Duality Theorem~2.1.5, as well as of Gordon's
Conjecture~1 (Theorem~2.1.4), it is enough to prove the
following:

\proclaim{2.3.3.~Proposition}
Let \,$(G,\Go,\Gf)$ be a condensing IMG group triplet with
hyperfinite abelian ambient group~$G$. Then
$$
\Go^{{\ort}{\ort}} = \Go\,.
$$
\endproclaim

\demo{Proof}
Let $\VV$ be a set of internal valuations on $G$ such that
$$
\Go = \bigl\{x \in G; (\all \ro \in \VV)(\ro(x) \apr 0)\bigr\}\,.
$$
For any $\ro \in \VV$ and a positive $r \in \sbR$ we denote the
internal closed ball of radius $r$
$$
B_\ro(r) = \{x \in G\: \ro(x) \le r\}\,,
$$
and define the internal function $h_{\ro r}\:G \to \sbR$ by
$$
h_{\ro r}(x) = \max\biggl\{1 - {\ro(x) \over r},\,0\biggr\}\,.
$$
Obviously, $h_{\ro r}$ is even, nonnegative, and
$\supp h_{\ro r} \sbs B_\ro(r)$. For $0 \not\apr r < \infty$,
both the sets $\supp h_{\ro r}$ and $B_\ro(r)$ are between $\Go$
and $\Gf$. Moreover, for $r$ noninfinitesimal, $h_{\ro r}$
is also $S$-continuous and
$\|h_{\ro r}\|_{_1} \ge d\,|B_\ro(r/2)|/2 \not\apr 0$.
By Corollary~2.3.2,
$$
\Spec_t(h_{\ro r}) \sbs \Go^{\ort}\,,
$$
for $0 < t \le 1$ once $t$ is noninfinitesimal, too. Further on,
$\supp(h_{\ro r} * h_{\ro r}) \sbs B_\ro(2r)$. In order to apply
Proposition~2.2.3 we need the estimate
$$
{\|h_{\ro r}\|_{_1} \over \|h_{\ro r}\|_{_2} \,\|1_{B_\ro(2r)}\|_{_2}}
\ge {{1 \over 2}\,|B_\ro(r/2)| \over \sqrt{|B_\ro(r)|\,|B_\ro(2r)|}}
\ge {|B_\ro(r/2)| \over 2\,|B_\ro(2r)|}
\not\apr 0\,.
$$
Thus picking some standard positive $\alfa < \pi/2$, there is a standard
positive
$$
t \le
{\|h_{\ro r}\|_{_1} \over \|h_{\ro r}\|_{_2} \,\|1_{B_\ro(2r)}\|_{_2}}
   \cd \sqrt{\cos\alfa \over 1+\cos\alfa}\,.
$$
For such $\alfa$ and $t$ we have
$$
\Bohr_{\alfa}\bigl(\Spec_t(h_{\ro r})\bigr)
\sbs B_\ro(2r) - B_\ro(2r) \sbs B_\ro(4r)\,.
$$
\vbox{\flushpar
Consequently,
$$
\Go^{{\ort}{\ort}} = \Bohr_{\alfa}\bigl(\Go^{\ort}\bigr) \sbs
\Bohr_{\alfa}\bigl(\Spec_t(h_{\ro r})\bigr) \sbs B_\ro(4r)\,.
$$
As $\ro$ and $r$ were arbitrary, this is enough to establish the
statement.
}
\enddemo

\remark{Remark}
In fact the subgroup $\Gf$ plays no role in the above statement nor
in its proof. The only thing we need to assume is the existence of
an internal set $A$ subject to $\Go \sbs A \sbs G$ such that
$|A|/|B| < \infty$ for each internal $B$ such that
$\Go \sbs B \sbs A$.
\endremark
\medskip

Having $d$ as a normalizing multiplier for the triplet $(G,\Go,\Gf)$,
then, in order to get the Fourier inversion formula and Plancherel
identity, we need to normalize the inner product, Fourier transform,
etc., on the space $\sbC^{\dG}$  of internal functions $\dG \to \sbC$,
defined on the internal dual group $\dG$, by means of the scaling
coefficient
$$
\hat d =  {1 \over d\lv G\rv}\,.
$$
In view of the canonical isomorphism of the observable trace
$\dG^\be = \Go^{\ort}\!/\Gf^{\ort}$ and the dual group $\wh{G^\be}$,
there naturally arises the question whether such a $\hat d$ is a
normalizing multiplier for the dual triplet
$\bigl(\dG,\Gf^{\ort},\Go^{\ort}\bigr)$. In that case (and only in
that case) the measure $\mb_{\hat d}$, obtained by pushing down the
Loeb measure $\lam_{\hat d}$ from $\dG$ to the observable trace
$\dG^\be \iso \wh{G^\be}$, will be a Haar measure on $\dG^\be$.
The second of Gordon's conjectures states that the response to the
above question is affirmative.

\proclaim{2.3.4.~Theorem}{\bf [Gordon's Conjecture~2]}
Let \,$(G,\Go,\Gf)$ be a condensing IMG group triplet with
hyperfinite abelian ambient group $G$. Then for any internal
set $D$ such that \,$\Go \sbs D \sbs \Gf$, \,$d = 1/|D|$ is a
normalizing multiplier for $(G,\Go,\Gf)$ and \,$\hat d = |D|/|G|$
\,is a normalizing multiplier for the dual triplet
$\bigl(\dG,\Gf^{\ort},\Go^{\ort}\bigr)$. More generally, if \,$d$
is any normalizing multiplier for \,$(G,\Go,\Gf)$, then
\,$\hat d = (d\lv G\rv)^{-1}$ \,is a normalizing multiplier for
\,$\bigl(\dG,\Gf^{\ort},\Go^{\ort}\bigr)$.
\endproclaim

\demo{Proof}
Since, for any internal set $D$ such that $\Go \sbs D \sbs \Gf$,
\,$d = 1/\lv D\rv$ is a normalizing multiplier for the triplet
$(G,\Go,\Gf)$ and $d_1/d_2 \in \FR \sms \IR$ for any its normalizing
multipliers $d_1$, $d_2$, it suffices to find a single couple
of internal sets $D$ between $\Go$, $\Gf$ and $\Gama$ between
$\Gf^{\ort}$, $\Go^{\ort}$, such that the quotient
${\lv D\rv \lv \Gama\rv/|G|}$ is neither infinite nor infinitesimal.

Starting with any internal valuation $\ro \in \VV$, we have
$\Go \sbs B_\ro(r) \sbs \Gf$ whenever $0 < r \in \bR$.
Also, we can take arbitrary standard $\alfa$ and $t$ subject to
$0 < \alfa < \pi/2$, $0 < t \le \cos\alfa$. Then the spectral set
$\Gama = \Spec_t(h_{\ro r})$ is internal and, by Lemma~2.2.2 and
Corollary~2.3.2, it satisfies the inclusions
$$
\Gf^{\ort} \sbs \Bohr_\alfa\bigl(B_\ro(r)\bigr)
   \sbs \Spec_t(h_{\ro r}) \sbs \Go^{\ort}\,.
$$
Denoting $f = h_{\ro r}$, $D = B_\ro(2r)$, we have
$\supp(f * f) \sbs D \sbs \Gf$, and multiplying the
inequalities in Proposition~2.2.4 by the factor
$\|1_D\|_{_2}^2 = |D|/|G|$ we get
$$
1 - t^2\,{\|f\|_{_2}^2\,\|1_D\|_{_2}^2  \over \|f\|_{_1}^2}
\le {\lv D\rv \lv\Spec_t(f)\rv \over |G|}
\le {1 \over t^2}\,{\|f\|_{_2}^2 \,\|1_D\|_{_2}^2
       \over \|f\|_{_1}^2}\,.
$$
As
$$
{\|f\|_{_2}^2 \,\|1_D\|_{_2}^2 \over \|f\|_{_1}^2} \le
{|B_\ro(r)|\,|B_\ro(2r)| \over \bigl({1 \over 2}|B_\ro(r/2)|\bigr)^2}
< \infty\,,
$$
the upper bound is finite for any standard $t > 0$. For the same
reason, it is possible to find a standard $t$, such that
$0 < t \le \cos\alfa$, and making the lower bound positive and
noninfinitesimal.
\enddemo

\remark{Remark}
Gordon's Conjectures 1 and 2 were proved in \cite{Go1} for
condensing IMG group triplets $(G,\Go,\Gf)$ with hyperfinite
abelian ambient group $G$ in case that $\Go$ is an intersection
and $\Gf$ is a union of countably many internal sets, and there
is an internal subgroup $K \sbs G$ such that $\Go \sbs K \sbs \Gf$.
These assumptions mean that the observable trace $\Gb = \Gf/\Go$ is
metrizable and $\sig$-compact and contains a compact open subgroup.
They particularly cover the cases of internal $\Go$ or $\Gf$,
corresponding to discrete countable or metrizable compact $\Gb$,
respectively. Using these results the {\it existence\/} of a triplet
$(G,\Go,\Gf)$ satisfying both the Conjectures with observable trace
$\Gf/\Go \iso \Gb$ was proved in \cite{Go2} for any metrizable
$\sig$-compact abelian group $\Gb$.
\endremark

\specialhead
2.4. Some (mainly) standard equivalents: Hrushovski style theorems
\endspecialhead
\flushpar
In this section we list some direct standard consequences
(in fact, equivalents) of Theorem~2.1.4 and Proposition~2.3.3.
The formulations of 2.4.2 and 2.4.3 below, to some extent,
remind of the formulation of Hrushovski's structure theorem
\cite{Hr, Theorem~1.1, Corollaries~4.13 and 4.15}
(see also \cite{BGT, Theorem~6.18}). As both the implications
2.1.4$\imp$2.4.2 and 2.3.3$\imp$2.4.3 can be proved using rather
a similar way of argumentation, we give just the proof of the
former (which is technically more complicated), introducing
an additional nonstandard equivalent 2.4.1.

In what follows the expression $[A:B]$ denotes any of the
indices $\lfloor A : B\rfloor$, $\lfloor A : B\rfloor_\ii$,
$\lceil A : B\rceil$, introduced in Section~1.2, or the quotient
$|A|/|B]$ (see also Proposition~1.5.5).

As we will show within short, Gordon's Conjecture~1 (Theorem~2.1.4)
is equivalent to a kind of {\it stability principle\/} for characters
of hyperfinite abelian groups. The following notions are needed for the
nonstandard formulation of this principle.

Let $G$ be an internal abelian group and $X$ be a (not necessarily
internal) subset of \,$G$. Two internal mappings $f\:C \to \sbT$,
$g\:D \to \sbT$, defined on sets $C,D \sbs G$, are said to be
{\it infinitesimally close on the set} \,$X$ if \,$X \sbs C \cap D$
and \,$f(x) \apr g(x)$ for all $x \in X$. We say that $g$ is
{\it almost homomorphic on the set} \,$X$ if \,$X \sbs D$ and, for
all $x,y \in X$, $x+y \in X$ implies
$$
g(x+y) \apr g(x)\,g(y)\,.
$$

\proclaim{2.4.1.~Theorem}
Let $(G,\Go,\Gf)$ be a condensing IMG group triplet with
hyperfinite abelian ambient group $G$, and \,$g\:D \to \sbT$ be
an internal mapping, where $\Gf \sbs D \sbs G$. If \,$g$ is
$S$-continuous and almost homomorphic on \,$\Gf$, then there exists
an internal $S$-continuous character $\gama\:G \to \sbT$, i.e.,
$\gama \in \Go^{\ort}$, such that $g$ and $\gama$ are
infinitesimally close on $\Gf$.
\endproclaim

Preliminarily, the reader should just notice that, due to
{\it saturation}, if the internal mapping $g\: D \to \sbT$, where
$\Gf \sbs D \sbs G$, is almost homomorphic on $\Gf$, then it must
be almost homomorphic on some symmetric internal set $C$ such that
$\Gf \sbs C \sbs D$; similarly, if \,$g$ and $\gama$ are infinitesimally
close on $\Gf$, then they must be infinitesimally close on some
symmetric internal set between $\Gf$ and $D$. At the same time,
internal mappings $g\: D \to \sbT$ almost homomorphic on the galaxy
$\Gf$ share the following property with internal homomorphisms
$\gama \in \dG$: such a $g$ is $S$-continuous on $\Gf$ if and only if,
given any (standard) $\alfa \in (0,2\pi/3)$, there is some internal
set $A$ between $\Go$ and $\Gf$ such that $\lv\arg g(x)\rv \le \alfa$
for all $x \in A$ (cf. Lemma~2.1.1 and the beginning of the proof of
Proposition~2.1.2).

The reformulation of 2.4.1 in standard terms is a highly uniform,
but rather cumbersome, {\it stability principle\/} for characters of
finite abelian groups. For the sake of its formulation, as well as
of the proof of the equivalence of 2.1.4, 2.4.1 and 2.4.2 (that way
proving Theorems~2.4.1 and 2.4.2, as well), we have to introduce some
further notions.

Let $G$ be an abelian group, and $\eps > 0$ be a (standard) real.
Two mappings $f\:C \to \bT$, \,$g\:D \to \bT$, defined on subsets
$C, D \sbs G$, are said to be {\it $\eps$-close on the set}
\,$X \sbs C \cap D$ if
$$
\lv\arg{f(x)\over g(x)}\rv \le \eps\,,
$$
for all $x \in X$. We say that $g$ is {\it $\eps$-homomorphic on the
set} \,$X \sbs D$ if for all $x,y \in X$ the condition $x+y \in X$
implies
$$
\lv\arg{g(x)\,g(y) \over g(x+y)}\rv \le \eps\,.
$$
Finally, $g$ \,is called a {\it partial $\eps$-homomorphism\/} if
it is $\eps$-homomorphic on its domain~$D$.

\proclaim{2.4.2.~Theorem}
Let $\alfa,\eps \in (0,2\pi/3)$, $k \ge 1$ and $(q_j)_{j=1}^\infty$
be any sequence of reals $q_j \ge 1$. Then there exist $m \ge 1$,
$n \ge k$ and $\delta > 0$, all depending just on $\alfa,\eps$,
$k$ and the sequence $(q_j)$, such that the following holds:

Let \,$G$ be a finite abelian group and
\,$0 \in A_n \sbs \ldots \sbs A_1 \sbs A_0 \sbs G$ be symmetric sets
such that
$$
A_j + A_j  \sbs A_{j-1}\,,
\qquad\text{and}\qquad
[A_{j-1} : A_j] \le q_j\,,
$$
for $1 \le j \le n$. Then, for every partial $\delta$-homomorphism
\,$g\:mA_0 \to \bT$ such that \,$\lv\arg g(x)\rv \le \alfa$ for
$x \in A_k$, there exists a homomorphism \,$\gama\: G \to \bT$
such that $g$ and $\gama$ are $\eps$-close on $A_0$.
\endproclaim

\demo{Proof of \,{\rm 2.1.4$\imp$2.4.2}}
Assume that 2.1.4 holds and 2.4.2 is not true. Then we can fix
some $\alfa$, $\eps$, $k$ and a sequence $(q_j)$ witnessing a
counterexample. Let $(\delta_n)$ be any strictly decreasing sequence
such that $\delta_n \to 0$. Then for all $m = n \ge k$ there exist a
finite abelian group $G_n$ and symmetric sets
$0 \in A_{nn} \sbs \ldots \sbs A_{n1} \sbs A_{n0} \sbs G_n$
such that
$$
A_{nj} + A_{nj} \sbs A_{n{\hskip.5pt}j-1}\,,
\qquad\text{and}\qquad
[A_{n{\hskip.5pt}j-1} : A_{nj}] \le q_j\,,
$$
for $1 \le j \le n$, as well as a partial $\delta_n$-homomorphism
\,$g_n\:nA_{n0} \to \bT$ \,such that $\lv\arg g_n(x)\rv \le \alfa$ for
$x \in A_{nk}$. On the other hand, for every genuine homomorphism
\,$\gama_n\:G_n \to \bT$, there is an $x_n \in A_{n0}$ such that
$$
\lv\arg{\gama_n(x_n) \over g_n(x_n)}\rv > \eps\,.
$$
Let $\DD$ be any nontrivial ultrafilter on the set
$I = \{n \in \bN\: n \ge k\}$. Then the ultraproduct
$G = \prod_{n \in I} G_n/\DD$ of the finite groups $G_n$ is a
hyperfinite abelian group (within the $\aleph_1$-saturated nonstandard
universe obtained via the ultrapoduct construction modulo~$\DD$).
For each $n \in I$, $n < j$, we put $A_{nj} = \{0\}$ and take the
ultraproduct $A_j = \prod_{n\in I} A_{nj}/\DD$ considered as an
internal subset of \,$G$. Finally we put
$$
\Go = \bigcap_{j \in \bN} A_j\,,
\qquad\qquad
\Gf = \bigcup_{j \in \bN} jA_0\,.
$$
Then $\Go$, as an intersection of countably many internal sets, is
a monadic subgroup of \,$G$ and $\Gf$, as a union of countably many
internal sets, is a galactic subgroup of \,$G$. Obviously,
$\Go \sbs \Gf$. By {\L}os' theorem, $[A_{j-1} : A_j] \le q_j$,
as well as $[jA_0 : A_1] \le q_1^j$ for each $j \ge 1$. It follows
easily that $[V : U] < \infty$ for any internal sets $U$, $V$ such
that $\Go \sbs U \sbs V \sbs \Gf$, hence, by Proposition~1.5.5,
$(G,\Go,\Gf)$ is a condensing IMG triplet with hyperfinite abelian
ambient group $G$.

Let us form the internal mapping \,$g = (g_n)_{n\in I}/\DD$. Then,
for each $n$, $g$ is a partial $\delta_n$-homomorphisms from the
internal set $\prod_{n\in I} nA_{n0}/\DD \sps \Gf$ to the ultrapower
group $\sbT = \bT^I\!/\DD$. Thus $g$ is almost homomorphic on $\Gf$
and it maps the internal set $A_k = \prod_n A_{nk}/\DD$,
satisfying $\Go \sbs A_k \sbs \Gf$, into the arc
$\{c \in \sbT\: \lv\arg c\rv \le \alfa\}$. It follows that
$g$ is $S$-continuous on $\Gf$. Then the observable trace
$g^\be = {}^\co(g\rst\Gf)$ is a continuous character of the LCA group
$G^\be = \Gf/\Go$. By Theorem~2.1.4, there is an internal character
$\gama = (\gama_n)/\DD \in \Go^{\ort}$ such that $g(x) \apr \gama(x)$
for each $x \in \Gf$. On the other hand, there is a set $J \in \DD$
such that $\gama_n \in \dG_n$ for all $n \in J$. By our assumptions,
for each $n \in J$, there is an $x_n \in A_{n0}$ such that
$\lv\arg\bigl(\gama_n(x_n)/g_n(x_n)\bigr)\rv > \eps$. Let $x_n = 0$
for $n \in I \sms J$. Then, for $x = (x_n)/\DD \in A_0 \sbs \Gf$, we
have $g(x) \not\apr \gama(x)$\,---\,a~contradiction.
\enddemo

\remark{Remark}
Concerning $m$, $n$ the result is purely existential, giving no
upper bound for them (though we can always have $m=n$). On the
other hand, it seems rather probable that even a more uniform
version of Theorem~2.4.2 is true. We conjecture that, similarly
as in \cite{MZ, Theorem~4.1}, one can take any $\delta > 0$, such
that $\delta < \min\{\eps, \pi/2, 2\pi/3-\alfa\}$, and choose
$m$, $n$ depending additionally on $\delta$.
\endremark

\demo{Proof of \,{\rm 2.4.2$\imp$2.4.1}}
Let $(G,\Go,\Gf)$ be a condensing IMG group triplet with
hyperfinite abelian ambient group $G$, $D$ be an internal set such
that $\Gf \sbs D \sbs G$ and $g\:D \to \sbT$ be an internal mapping,
$S$-continuous and almost homomorphic on $\Gf$. We can also assume
that $g(0) = 1$. Thus fixing a standard $\alpha \in (0,2\pi/3)$, there
is a symmetric internal set $V$ between $\Go$ and $\Gf$ such that
$\lv\arg g(x)\rv \le \alfa$ for $x \in V$. Let us choose,
additionally, a $\beta \in (\alfa,2\pi/3)$.

Now, it is enough to show that, for each symmetric internal set $A$
such that $V+V \sbs A \sbs \Gf$, and each standard $\eps > 0$,
such that $\alfa+\eps \le \beta$, there is a
$\gama \in \Bohr_\beta(V)$ such that
$$
\biggl|\arg \frac{g(x)}{\gama(x)}\biggr| \le \eps
$$
for all $x \in A$. Due to {\it saturation}, there is then
a $\gama \in \Bohr_\beta(V) \sbs \Go^{\ort}$ such that
$\gama(x) \apr g(x)$ for all $x \in \Gf$, i.e., $\gama$
and $g$ will be infinitesimally close on $\Gf$.

So let us fix some $A$, $V$, $\alfa$, $\beta$ and $\eps$, satisfying
the above assumptions. Then there is a sequence of symmetric
internal sets $(A_j)_{j\in\bN}$ such that $A_0 = A$, $A_1 = V$
and $A_{j+1} + A_{j+1} \sbs A_j \sps \Go$ for each $j$.
Put $q_j = [A_{j-1} : A_j]$ for $1 \le j \in \bN$. Let $m$, $n$ and
$\delta$ be the numbers guaranteed to $\alfa$, $\eps$, $k=1$ and the
sequence $(q_j)$ by Theorem~2.4.2. By the {\it transfer principle\/}
they have to work for the hyperfinite abelian group $G$ and the
internal sets $A_0,\dots,A_n$, as well. Since $g$ is almost
homomorphic on $\Gf$, the more it is $\delta$-homomorphic
on $mA_0 \sbs \Gf$. It follows that there is an internal character
$\gama \in \dG$ such that $\lv\arg g(x)\gama(x)^{-1}\rv \le \eps$
for $x \in A_0 = A$. Obviously,
$\lv\arg\gama(x)\rv \le \alfa+\eps \le \beta$ for
$x \in A_1 = V$, hence $\gama \in \Bohr_\beta(V)$.
\enddemo

Let us close the circle by proving that Theorem~2.4.1 implies
Gordon's Conjecture~1 (Theorem 2.1.4), i.e., the surjectivity
of the canonical embedding $\dG^\be \to \wh{G^\be}$. In view of
Proposition~1.2.3, this implication is plain.

\demo{Proof of \,{\rm 2.4.1$\imp$2.1.4}}
Let $(G,\Go,\Gf)$ be a condensing IMG group triplet with
hyperfinite abelian ambient group $G$; let $\Gb = \Gf/\Go$ denote
its observable trace. Assume that $\ggam\:\Gb \to \bT$ is a
continuous character of \,$\Gb$. By Proposition~1.2.3, there is an
internal mapping $g\:D \to \sbT$ such that $\Gf \sbs D \sbs G$,
$g$ is $S$-contiuous on $\Gf$, and
$$
\ggam\bigl(x^\be\bigr) = {}^\co g(x)\,,
$$
for all $x \in \Gf$. Then
$$
g(x+y) \apr \ggam\bigl((x+y)^\be\bigr) =
\ggam\bigl(x^\be + y^\be\bigr) =
\ggam\bigl(x^\be\bigr)\,\ggam\bigl(y^\be\bigr)
\apr g(x)\,g(y)\,,
$$
for all $x,y \in \Gf$, i.e., $g$ is almost homomorphic on $\Gf$.
By 2.4.1, there is a $\gama \in \Go^{\ort}$ such that
$$
\gama(x) \apr g(x) \apr \ggam\bigl(x^\be\bigr)\,,
$$
for each $x \in \Gf$, i.e., $\gama^\be = \ggam$.
\enddemo

A brief inspection of the above proofs is instructive. In order
to derive 2.4.2 from 2.1.4 it suffices to assume that Gordon's
Conjecture~1 is true just for condensing IMG group triplets
$(G,\Go,\Gf)$ with hyperfinite $G$ (in some $\aleph_1$-saturated
nonstandard universe), such that $\Go$ is an intersection and $\Gf$
is a union of {\it countably many\/} internal sets. This was indeed
Gordon's original formulation of his conjecture in \cite{Go1}.
Conversely, in order to prove 2.4.1 (from which 2.1.4 easily
follows by the virtue of 1.2.3) as a consequence of 2.4.2, it is
enough to suppose that the latter is true for a single fixed
$\alfa \in (0,2\pi/3)$ and $k=1$, only.

Proposition~2.3.3, forming the essential part of the Triplet Duality
Theorem~2.1.5, can be restated in standard terms as follows:

\proclaim{2.4.3.~Theorem}
Let $\alfa,\beta \in (0,2\pi/3)$ and $(q_j)_{j=1}^\infty$ be any
sequence of reals $q_j \ge 1$. Then there exist an $n \in \bN$,
depending just on $\alfa$, $\beta$ and the sequence $(q_j)$,
such that the following holds:

Let \,$G$ be a finite abelian group and
\,$0 \in A_n \sbs \ldots \sbs A_1 \sbs A_0 \sbs G$ be symmetric
sets such that
$$
A_j + A_j  \sbs A_{j-1}\,,
\qquad\text{and}\qquad
[A_{j-1} : A_j] \le q_j\,,
$$
for $1 \le j \le n$. Then \,$\Bohr_\beta(\Bohr_\alfa(A_n)) \sbs A_0$.
\endproclaim

Theorem~2.4.3 can be derived from Proposition~2.3.3 using the
ultraproduct construction, similarly as (and more easily than)
the implication 2.1.4$\imp$2.4.2. To this end it would be enough to
assume that 2.3.3 is true just in case when $\Go$ is an intersection
of a countable family of internal subsets of \,$G$, again.
On the other hand, in order to derive 2.3.3 from 2.4.3, it would
suffice to suppose that the latter is true for one fixed couple
$\alfa, \beta \in (0,2\pi/3)$, only.

\remark{Remark}
Theorem 2.4.1 still holds for condensing IMG group triplets
$(G,\Go,\Gf)$ with an arbitrary {\it internal\/} abelian ambient
group $G$, and not just a hyperfinite one. Moreover, most of the
results of Sections 2.1 and 2.3 admit analogous generalizations.
Similarly, both Theorems~2.4.2 and 2.4.3 (and maybe even the
strengthening of 2.4.2 mentioned in the Remark following the
proof of 2.4.1$\imp$2.4.2) remain true when replacing\,
\,{\it``Let \,$G$ be a finite abelian group and
\,$0 \in A_n \sbs \dots \sbs A_1 \sbs A_0 \sbs G$ be symmetric
sets\dots''}
\,by\,
{\it``Let \,$G$ be a Hausdorff LCA group and
\,$A_n \sbs \dots \sbs A_1 \sbs A_0$ be compact symmetric
neighborhoods of \,$0$ in $G$\dots''}
\,(without changing the rest) in their formulation.

In order to prove all that, it would be necessary (and sufficient)
to deal with condensing IMG triplets of the form $(G,\Go,\Gf)$,
where $G$ is an arbitrary internal Hausdorff LCA group (i.e., a
\,${}^*$(Hausdorff LCA) group), and define their dual triplets
$\bigl(\dG,\Gf^{\ort},\Go^{\ort}\bigr)$, with $\dG$ denoting the
internal dual group of \,$G$, as well as to generalize the results
of Section 2.2 from finite abelian groups to arbitrary Hausdorff
LCA groups. That, however, though possible and\,---\,as we shall see
in the next section\,---\,unavoidable to some extent, would be at
odds with the leading intentions of the present paper, namely to study
the LCA groups and the Fourier transform on functional spaces over them
by means of their approximations by (hyper)finite abelian groups and
the discrete Fourier transform on them.

That way generalized Theorem~2.4.3 can be schematically written
as follows:
$$
(\all \alfa,\beta)(\all q_1,\dots, q_j,\dots)(\exs n)
(\all G)(\all A_0,A_1\dots,A_n)
\bigl[(\all j \le n)(\dots) \imp \dots\bigr]\,.
$$
However, from the Pontryagin-van\,Kampen Duality Theorem it follows,
only,
$$
(\all \alfa,\beta)(\all q_1,\dots, q_j,\dots)(\all G)
(\all A_0,A_1\dots,A_j,\dots)
\bigl[(\all j)(\dots) \imp (\exs n)(\dots)\bigr]\,.
$$
In other words, such a generalization of 2.4.3 adds a considerable
uniformness to the above consequence of the Pontryagin-van\,Kampen
Duality Theorem.

Analogous remarks could be made on behalf of the indicated
generalizations of Theorems 2.4.1 and 2.4.2 and their relation to
some standardly formulated stability results for characters of LCA
groups with respect to the compact-open topology from \cite{MZ} and
\cite{Zl2}. However, the big amount of technicalities needed for
their formulation and comparison would take us too far from the
main lines and goals of the present paper.
\endremark

\specialhead
2.5.~Simultaneous approximation of an LCA group and its dual
\endspecialhead
\flushpar
Now, we are going to apply the results of the previous sections
to triplets arising from HFI approximations of LCA groups. One can
expect that an HFI approximation $\eta\:G \to \Gb$ of an LCA group
$\Gb$ gives rise to an HFI approximation $\phi\:\dG \to \dGb$ of
the dual group $\dGb$. This is true in some sense, however, in
general, the connection between these two approximations is not
as straightforward as one wished. Some of the material of this
section is partly covered in the Introduction and \S2.4 of
Gordon's book \cite{Go2}.

As we will be dealing with plenty of more or less canonical
isomorphisms of topological groups, first we have to make clear
which of them we intend to exploit for identification of the
isomorphic objects, and which of them we still view as isomorphisms
of different objects.

Let $\Gb$ be a Hausdorff LCA group, viewed as a subgroup of its
nonstandard extension $\sGb$. Let us denote
$$
\IG = \Mon(0)
\qquad\text{and}\qquad
\FG = \Ns(\sGb)\,.
$$
Then $\bigl(\sGb,\IG,\FG\bigr)$ is a condensing IMG group triplet
with internal abelian ambient group $\sGb$. We identify the
canonically isomorphic LCA groups $\Gb$ and the observable trace
$\FG/\IG$ of this triplet, as well as the standard part map
$\xb \mto {}^\co\xb\:\FG \to \Gb$ and the observable trace map
$\xb \mto \xb^\be\:\FG \to \FG/\IG$.

Similar but less straightforward identification applies also to
the dual group $\dGb$. More specifically, due to the {\it transfer
principle}, the nonstandard extension $\sdGb$ of the dual group
$\dGb$ coincides with the internal dual group $\wh{\sGb}$ of
\,$\sGb$, consisting of all internal ${}^*$continuous characters
$\ggam\:\sGb \to \sbT$. However, the elementary embedding
$\dGb \to \sdGb$ sends a character $\ggam\:\Gb \to \bT$ not to
itself but to its nonstandard extension ${}^*\ggam\:\sGb \to \sbT$.
Putting
$$
\IdG = \Mon(1_{\sGb})\,,
\qquad\text{and}\qquad
\FdG = \Ns\bigl(\sdGb\bigr)\,,
$$
one can easily realize that
$$\align
\IdG = \FG^{\ort} &=
\bigl\{\ggam \in \sdGb\:
   (\all \xb \in \FG)(\ggam(\xb) \apr 1)\bigr\}\,,\\
\FdG = \IG^{\ort} &=
\bigl\{\ggam \in \sdGb\:
(\all \xb \in \IG)(\ggam(\xb) \apr 1)\bigr\}\,,
\endalign$$
so that the condensing IMG group triplet
$\bigl(\sdGb,\IdG,\FdG\bigr)$ coincides with the dual triplet
$\bigr(\wh{\sGb},\FG^{\ort},\IG^{\ort}\bigr)$ of $(\sGb,\IG,\FG)$.
We identify the isomorphic groups $\dGb$ and $\FdG/\IdG$,
as well as the standard part map
$\ggam \mto {}^\co\ggam\:\FdG \to \dGb$ and the observable trace map
$\ggam \mto \ggam^\be\:\FdG \to \FdG/\IdG$, again.

Now assume that $\eta\: G \to \sGb$ is an HFI approximation of
\,$\Gb$ by a hyperfinite abelian group $G$ and $(G,\Go,\Gf)$ is the
condensing IMG group triplet arising from this approximation, i.e.,
$$
\Go = \eta^{-1}[\IG]
\qquad\text{and}\qquad
\Gf = \eta^{-1}[\FG]\,.
$$
Then $\eta\:(G,\Go,\Gf) \to \bigl(\sGb,\IG,\FG\bigr)$ is
a triplet  isomorphism, and its observable trace
$\eeta = \eta^\be\:\Gf/\Go \to \FG/\IG$ becomes an isomorphism
between the two representations of the original LCA group $\Gb$ as
the observable trace $\Gb \iso \Gf/\Go$ and the nonstandard hull
$\Gb \iso \FG/\IG$.

Let us form also the dual triplet
$\bigl(\dG,\Gf^{\ort},\Go^{\ort}\bigr)$. Making use of the
canonical isomorphism $\Go^{\ort}\!/\Gf^{\ort} \to \wh{\Gf/\Go}$
from Theorem~2.1.4 we identify the dual group $\wh{\Gf/\Go}$ of
the observable trace $G^\be = \Gf/\Go$ and the observable
trace $\dG^\be = \Go^{\ort}\!/\Gf^{\ort}$ of the dual triplet
$\bigl(\dG,\Gf^{\ort},\Go^{\ort}\bigr)$.

Keeping in mind the just introduced identifications and applying
the duality functor to the isomorphism $\eeta^{-1}\:\Gb \to \Gf/\Go$
of LCA groups, we obtain the isomorphism
$\pphi\:\Go^{\ort}\!/\Gf^{\ort} \to \dGb$ between their duals,
given by
$$
\pphi(\ggam) = \ggam \co \eeta^{-1}\,,
$$
for $\ggam \in \wh{\Gf/\Go} = \Go^{\ort}\!/\Gf^{\ort}$.

By Corollary~1.2.4, any isomorphism
$\pphi\:\Go^{\ort}\!/\Gf^{\ort} \to \FdG/\IdG$ between the
two representations of the dual group
$\dGb \iso \Go^{\ort}\!/\Gf^{\ort}$ and $\dGb \iso \FdG/\IdG$ is
the observable trace of some internal mapping $\phi\:\dG \to \sdGb$.
This is to say that
$\phi\:\bigl(\dG,\Gf^{\ort},\Go^{\ort}\bigr) \to
\bigl(\sdGb,\IdG,\FdG\bigr)$ is a triplet isomorphism and
$$
\pphi\bigl(\gama^\be\bigr) = \gama^\be \co \eeta^{-1}
= \phi(\gama)^\be\,,
$$
for $\gama \in \Go^{\ort}$. Then $\phi$ necessarily is almost
homomorphic on $\Go^{\ort}$, and we have
$$
\Gf^{\ort} = \phi^{-1}\bigl[\IdG\bigr]
\qquad\text{and}\qquad
\Go^{\ort} = \phi^{-1}\bigl[\FdG\bigr]\,,
$$
so that $\phi\:\dG \to \sdGb$ is an HFI approximation of \,$\dGb$
and $\bigl(\dG,\Gf^{\ort},\Go^{\ort}\bigr)$ is the IMG group triplet
arising from this approximation. In the particular case of $\pphi$
given by the assignment $\ggam \mto \ggam \co \eeta^{-1}$, as above,
we finally have
$$
(\phi\,\gama)(\eta\,x) \apr
\bigl(\pphi\,\gama^\be\bigr)\bigl(\eeta\,x^\be\bigr) =
\bigl(\gama^\be \co \eeta^{-1}\bigr)\bigl(\eeta\,x^\be\bigr)
= \gama^\be\bigl(x^\be\bigr) \apr \gama(x)\,,
$$
for any $x \in \Gf$, $\gama \in \Go^{\ort}$.

\remark{Remark}
Though this is not an essential point, in view of Corollary 1.5.12
we can assume, for convenience' sake, that the approximation $\eta$
is injective and preserves 0 and inverses. Then there is some
internal, necessarily surjective mapping $\eta'\:\sGb \to G$ such
that $\eta' \co \eta = \Id_G$, \,$\xb \apr (\eta \co \eta')(\xb)$,
for $\xb \in \FG$, and $\eta'$ preserves 0 and inverses, as well.
Now, it is natural to define the mapping $\phi$ by the assignment
$\gama \mto \gama \co \eta'$, for $\gama \in \dG$. Such a $\phi$
would be injective and strictly preserving the pointwise
multiplication of functions (hence the trivial character $1_G$ and
pointwise inverses, too). Unfortunately, this natural attempt does
not work. The reason is that $\eta$, $\eta'$, in spite of being
almost homomorphic on $\Gf$, $\FG$, respectively, are not genuine
homomorphisms, in general. Hence the mapping $\gama \co \eta'$
(though preserving 0 and inverses) would be just almost
homomorphic on $\FG$, again, and one cannot assure that
$\gama \co \eta' \in \sdGb$, for $\gama \in \dG$. Thus there
seems to be no canonical way how to determine the HFI approximation
$\phi\:\dG \to \sdGb$ of \,$\dGb$ right away from the HFI
approximation $\eta\:G \to \sGb$ of \,$\Gb$.

On the other hand, for each $\ggam \in \dGb$, the composition
${}^*\ggam \co \eta\:G \to \sbT$ is almost homomorphic on $\Gf$
(in fact, it is an $S$-continuous lifting of $\ggam$).
By Theorem~2.1.4,  there is a genuine homomorphism
$\gama \in \Go^{\ort} \sbs \dG$ such that
$\ggam = ({}^*\ggam \co \eta)^\be = \gama^\be$, more precisely,
$$
({}^*\ggam \co \eta)(x) \apr \gama(x)
$$
for each $x \in \Gf$. Though ${}^*\ggam \co \eta \nin \dG$,
in general, we shall see in Section~3.3 that it can be used
directly instead of the genuine homomorphism
$\gama \in \Go^{\ort}$ in approximation of the Fourier transform
on $\Gb$ by means of the inner product on the hyperfinite
dimensional unitary space $\sbC^G$.
\endremark
\medskip

Now assume that $G$ and  $H$ are hyperfinite abelian groups,
$\eta\:G \to \sGb$ is an HFI approximation of the Hausdorff LCA
group $\Gb$ and $\phi\:H \to \sdGb$ is an HFI approximation of
its dual group $\dGb$. Let $(G,\Go,\Gf)$ and $(H,\Ho,\Hf)$
be the condensing IMG triplets arising from these approximations,
i.e.,
$$\xalignat{2}
\Go &= \eta^{-1}[\IG]\,, & \Gf &= \eta^{-1}[\FG]\,,\\
\Ho &= \phi^{-1}\bigl[\IdG\bigr]\,, &
   \Hf &= \phi^{-1}\bigl[\FdG\bigr]\,.
\endxalignat$$
Following Gordon \cite{Go2, p.~148} we say that the approximations
$\eta$ and $\phi$ are {\it dual\/} to each other or that they form
an {\it adjoint pair\/} if \,$H = \dG$, $\Ho = \Gf^{\ort}$,
$\Hf = \Go^{\ort}$, i.e., if the triplet $(H,\Ho,\Hf)$
coincides with the dual triplet
$\bigl(\dG,\Gf^{\ort},\Go^{\ort}\bigr)$ of $(G,\Go,\Gf)$, and
$$
(\phi\,\gama)(\eta\,x) \apr \gama(x)
$$
holds for all $x \in \Gf$, $\gama \in \Hf$.
Inspecting our accounts and making use of the notions just introduced,
we find that, in any sufficiently saturated nonstandard universe, we
have proved the following results.

\proclaim{2.5.1.~Theorem}
Let $(G,\Go,\Gf)$ be an IMG group triplet with hyperfinite
abelian ambient group $G$, arising from an HFI approximation
$\eta\:G \to \sGb$ of the Hausdorff LCA group $\Gb$, and the
isomorphism $\eeta:\Gf/\Go \to \Gb$ of LCA groups be the observable
trace of \,$\eta$. Let further $\phi\:\dG \to \sdGb$ be an HFI
approximation of the dual LCA group $\dGb$ such that its observable
trace $\pphi = \phi^\be\:\Go^{\ort}\!/\Gf^{\ort} \to \dGb$ is the
dual isomorphism corresponding to $\eeta^{-1}\:\Gb \to \Gf/\Go$.
Then the HFI approximations $\eta$ and $\phi$ are dual to each
other.
\endproclaim

\proclaim{2.5.2.~Corollary}%
{\bf [Adjoint Hyperfnite LCA Group Approximation Theorem]}
To each HFI approximation $\eta\:G \to \sGb$ of a Hausdorff LCA
group $\Gb$ by a hyperfinite abelian group $G$ there exists a dual
HFI approximation $\phi\:\dG \to \sdGb$ of the dual LCA group $\dGb$
by the dual hyperfinite abelian group $\dG$.
\endproclaim

On the other hand, even if \,$\eta$ is injective, our accounts,
so far, are not sufficient to establish the analogous property
for $\phi$.

It is worthwhile to realize that, in a sufficiently saturated nonstandard
universe, the conditions defining the notion  of an adjoint pair of HFI
approximations are redundant to some extent. Here is a more detailed
account of their relation.

\proclaim{2.5.3.~Lemma}
Let \,$G$ be a hyperfinite abelian group, $\eta\:G \to \sGb$,
$\phi\:\dG \to \sdGb$ be HFI approximations of the Hausdorff LCA
groups \,$\Gb$, $\dGb$, respectively, and \,$(G,\Go,\Gf)$,
$\bigl(\dG,\Ho,\Hf\bigr)$ be the condensing IMG triplets arising
from these approximations, such that
$$
(\phi\,\gama)(\eta\,x) \apr \gama(x)
$$
for all \,$x \in \Gf$, $\gama \in \Hf$. Then, as consequence,
{\parindent 21pt
\item{\sl(a)}
the inclusions \,$\Hf \sbs \Go^{\ort}$, \,$\Gf \sbs \Ho^{\ort}$,
\,$\Go \sbs \Hf^{\ort}$, and \,$\Ho \sbs \Gf^{\ort}$ are satisfied;
\item{\sl(b)}
any two of the reverse inclusions
\,$\Go^{\ort}\sbs \Hf$, \,$\Ho^{\ort} \sbs \Gf$,
\,$\Hf^{\ort} \sbs \Go$, and \,$\Gf^{\ort} \sbs \Ho$
are equivalent.
\item{}}
\vskip-12pt
\endproclaim

\demo{Proof}
(a) In order to prove the first inclusion, take any $\gama \in \Hf$.
As \,$\Hf = \phi^{-1}\bigl[\FdG\bigr]$ and \,$\FdG = \IG^{\ort}$, we
have $(\phi\,\gama)(\xb) \apr 1$ for each $\xb \in \IG = \Mon(0)$.
Since $\Go = \eta^{-1}[\IG]$, this implies
$(\phi\,\gama)(\eta\,x) \apr 1$ for each $x \in \Go$. From
$\Go \sbs \Gf$ we obtain
$$
\gama(x) \apr (\phi\,\gama)(\eta\,x) \apr 1\,,
$$
hence \,$\gama \in \Go^{\ort}$. From $\Hf \sbs \Go^{\ort}$ we
readily  get $\Go \sbs \Go^{\ort\ort} \sbs \Hf^{\ort}$. The third
and the forth inclusion now follow by the symmetry of the situation.

(b) For brevity's sake let us number the inclusions consecutively
from (i) to (iv).
First we prove (i)$\imp$(iv). Assume that $\Go^{\ort}\sbs \Hf$ and
take any $\gama \in \Gf^{\ort} \sbs \Go^{\ort} \sbs \Hf$. Then
$\phi(\gama) \in \FdG = \IG^{\ort}$, i.e., $\phi(\gama)$ is
$S$-continuous. As \,$\Ho = \phi^{-1}\bigl[\IdG\bigr]$ and
\,$\IdG = \FG^{\ort}$, in order to establish that $\gama \in \Ho$,
we are to show that \,$(\phi\,\gama)(\xb) \apr 1$ for each
$\xb \in \FG = \Ns(\sGb)$. Since there is an $x \in \Gf$ such
that $\xb \apr \eta(x)$, we have indeed
$$
(\phi\,\gama)(\xb) \apr (\phi\,\gama)(\eta\,x)
   \apr \gama(x) \apr 1\,.
$$
(ii)$\imp$(iii) can be proved analogously by a symmetry argument.

Next we show (iv)$\imp$(ii). From $\Gf^{\ort} \sbs \Ho$ we get
$\Ho^{\ort} \sbs \Gf^{\ort\ort} = \Gf$ by the Triplet Duality
Theorem~2.1.5 (see also Lemma~2.1.7). By symmetry we have
(iii)$\imp$(i), as well, closing the circle of implications.
\enddemo

In other words, under the assumptions of the Lemma, the HFI
approximations $\eta$ and $\phi$ form an adjoint pair if and only
if anyone (hence all) of the inclusions from (b) is true. We are
inclining to consider the first inclusion $\Go^{\ort} \sbs \Hf$
from (b) as intuitively the most appealing and  fundamental one.
In ``unzipped form'' it states that, for all $\gama \in \dG$,
$$
(\all x \in \Go)\bigl(\gama(x) \apr 1\bigr) \imp
   (\all \xb \in \IG)\bigl((\phi\,\gama)(\xb) \apr 1\bigr)\,,
$$
and, as $\gama(x) \apr (\phi\,\gama)(\eta\,x)$, it is equivalent to
$$
\bigl(\all \xb \in \eta[\Go]\bigr)
   \bigl((\phi\,\gama)(\xb) \apr 1\bigr) \imp
     (\all \xb \in \IG)\bigl((\phi\,\gama)(\xb) \apr 1\bigr)\,.
$$
The remaining inclusions in (b) can be ``unzipped'' in a similar way.

The standard meaning of the above accounts and results can be
formulated in terms of approximating systems. The equivalence of
both formulations could be established by referring to Nelson's
translation algorithm \cite{Ne}. However, this would obscure some
connections\,---\,that's why we offer a more detailed exposition,
revealing some insights.

Let $\Gb$ be a Hausdorff LCA group, $\Kb \sbs \Gb$, $\Gamb \sbs \dGb$
be compact sets, and $\Ub \sbs \Gb$, $\Omb \sbs \dGb$ be neighborhoods
of the neutral elements in $\Gb$, $\dGb$, respectively. Let further
$0 < \alfa < 2\pi/3$. Then the pairs $(\Kb,\Ub)$, $(\Gamb,\Omb)$
are called {\it $\alfa$-adjoint\/} if
$$\align
&\Ub \sbs \Bohr_\alfa(\Gamb) \sbs \Bohr_\alfa(\Omb) \sbs \Kb\,, \\
&\Omb \sbs \Bohr_\alfa(\Kb) \sbs \Bohr_\alfa(\Ub) \sbs \Gamb\,.
\endalign$$
The reader should notice that in such a case we necessarily have
$\Gamb = \Bohr_\alfa(\Ub)$ and $\Kb = \Bohr_\alfa(\Omb)$, so that
the $\alfa$-adjoint pairs  $(\Kb,\Ub)$, $(\Gamb,\Omb)$ are uniquely
determined by their second items $\Ub$, $\Omb$. They give rise to
a couple of \,$\alpha$-adjoint pairs if and only if
\,$\Ub \sbs \Bohr_\alfa(\Omb)$, or, equivalently,
$\Omb \sbs \Bohr_\alfa(\Ub)$; this is to say that
$$
\lv \arg\ggam(\xb)\rv \le \alfa
$$
for all $\xb \in \Ub$, $\ggam \in \Omb$.

Two directed double bases
$\bigl((\Kb_i,\Ub_i)\bigr)_{i \in I}$ in $\Gb$ and
$\bigl(\Gamb_i,\Omb_i)_{i \in I}$ in $\dGb$ over a directed poset
$(I,\le)$ are called {\it $\alfa$-adjoint\/} if, for each $i \in I$,
the pairs $(\Kb_i,\Ub_i)$, $(\Gamb_i,\Omb_i)$ are $\alfa$-adjoint.

Let us start with the following easy but useful observation.

\proclaim{2.5.4.~Lemma}
Let \,$\Gb$ be a Hausdorff LCA group and \,$0 < \alfa < 2\pi/3$.
Then there exist $\alfa$-adjoint directed double bases
$\bigl((\Kb_i,\Ub_i)\bigr)_{i \in I}$ in \,$\Gb$ and
\,$\bigl(\Gamb_i,\Omb_i)_{i \in I}$ in \,$\dGb$ over some directed
poset $(I,\le)$. Besides $\Gamb_i = \Bohr_\alfa(\Ub_i)$,
$\Kb_i = \Bohr_\alfa(\Omb_i)$ which are automatically satisfied,
one can additionally require that $\Ub_i = \Bohr_\alfa(\Gamb_i)$,
$\Omb_i = \Bohr_\alfa(\Kb_i)$ for each $i \in I$, as well.
\endproclaim

\demo{Proof}
In view of the Pontryagin-van\,Kampen Duality Theorem,
the statement becomes obvious once we realize that, for
each $\alfa \in (0,2\pi/3)$ and any DD base
$\bigl((\Kb_i,\Ub_i)\bigr)_{i \in I}$ in $\Gb$, the sets
$$
\Gamb_i = \Bohr_\alfa(\Ub_i)\,,
\qquad\qquad
\Omb_i = \Bohr_\alfa(\Kb_i)
$$
form a DD base $\bigl((\Gamb_i,\Omb_i)\bigr)_{i \in I}$
in $\dGb$.
\enddemo

Let $\Gb$ be a Hausdorff LCA group and $(I,\le)$ be a directed poset.
Assume that, for each $i \in I$, there are finite abelian groups
$G_i$, $H_i$, endowed with mappings $\eta_i\:G_i \to \Gb$,
$\phi_i\:H_i \to \dGb$, such that the families
$(\eta_i\:G_i \to \Gb)_{i \in I}$ and $(\phi_i\:H_i \to \dGb)_{i \in I}$
form approximating systems of the group $\Gb$ and of its dual group $\dGb$,
respectively. The two approximating systems are said to be {\it dual\/}
to each other or to form an {\it adjoint pair\/} if \,$H_i = \dG_i$, for
each $i \in I$, and the following conditions hold:
\smallskip
{\parindent 21pt
\item{1.}
for each $\eps > 0$ and for any compact sets $\Kb \sbs \Gb$,
$\Gamb \sbs \dGb$ there exists an $i \in I$ such that
$$
\biggl|\arg\frac{(\phi_j\,\gama)(\eta_j\,x)}{\gama(x)}\biggr|
   \le \eps\,,
$$
for all $j \ge i$ in $I$ and $x \in \eta_j^{-1}[\Kb]$,
$\gama \in \phi_j^{-1}[\Gamb]$.
\smallskip
\item{2.}
for some (or, equivalently, for each) $\alfa \in (0,2\pi/3)$ and for
every neighborhood  $\Ub$ of \,0 in $\Gb$ there exist a compact set
$\Gamb \sbs \dGb$ and an $i \in I$ such that, for all $j \ge i$ in $I$,
$$
\Bohr_\alfa\bigl(\eta_j^{-1}[\Ub]\bigr) \sbs \phi_j^{-1}[\Gamb]\,.
$$
\item{}}

A closer connection between adjoint pairs of HFI approximations
and adjoint pairs of approximating systems can be established by
means of the ultraproduct construction. We use essentially the
same notation as that fixed prior to Proposition~1.5.8, extending
it also to the dual group $\dGb$ and its HFI approximation
and approximating system. In particular, the front star ${}^*$
denotes the utrapower of the corresponding object and $G$, $\dG$
denote the ultraproducts of the finite groups $G_i$, $\dG_i$,
respectively. The cartesian product $I \cx \bN$ of the directed
poset $(I,\le)$ with the (totally) ordered set $(\bN,\le)$
is ordered coordinatewise, i.e., $(i,m) \le (j,m)$ if and only
if \,$i \le j$ and $m \le n$. Obviously, $(I \cx \bN, \le)$ is
a directed poset, again. Also it is clear that if
$\bigl((\Kb_i,\Ub_i)\bigr)_{i \in I}$ and
$\bigl((\Gamb_i,\Omb_i)\bigr)_{i \in I}$ are $\alfa$-adjoint DD bases
over $(I,\le)$ of \,$\Gb$ and \,$\dGb$, respectively, then putting
$\Kb_{in} = \Kb_i$, $\Ub_{in} = \Ub_i$, $\Gamb_{in} = \Gamb_i$,
$\Omb_{in} = \Omb_i$, for $i \in I$, $n \in \bN$, we obtain
$\alfa$-adjoint DD bases
$\bigl((\Kb_{in},\Ub_{in})\bigr)_{(i,n) \in I \cx \bN}$ and
$\bigl((\Gamb_{in},\Omb_{in})\bigr)_{(i,n) \in I \cx \bN}$
over $(I \cx \bN, \le)$ of \,$\Gb$ and $\dGb$, respectively, again.
These conventions are used throughout the following Proposition
and its proof, as well as in Theorems~2.5.6 and 2.5.9 to follow.

\proclaim{2.5.5.~Proposition}
Let \,$\Gb$ be a Hausdorff LCA group, $0 < \alfa < 2\pi/3$, and $\DD$
be a directed ultrafilter over a directed poset $(I,\le)$.
Let further $\bigl((\Kb_i,\Ub_i)\bigr)_{i \in I}$ and
$\bigl((\Gamb_i,\Omb_i)\bigr)_{i \in I}$ be $\alfa$-adjoint DD bases
of \,$\Gb$ and \,$\dGb$, respectively. Assume that $(G_i)_{i \in I}$ is
a system of finite abelian groups endowed with maps
$\eta_i\:G_i \to \Gb$, forming a well based approximating
system of the group $\Gb$ with respect to the DD base
$\bigl((\Kb_i,\Ub_i)\bigr)_{i \in I}$, and with maps
$\phi_i\:\dG_i \to \dGb$, forming a well based approximating
system of the dual group $\dGb$ with respect to the DD base
$\bigl((\Gamb_i,\Omb_i)\bigr)_{i \in I}$. Then the following
conditions hold true:
{\parindent 21pt
\item{\sl(a)}
If the approximating systems
$\bigl(\eta_i\:G_i \to \Gb\bigr)_{i\in I}$,
\,$\bigl(\phi_i\:\dG_i \to \dGb\bigr)_{i \in I}$ form an adjoint
pair, then the ultraproduct HFI approximations $\eta\:G \to \sGb$,
\,$\phi\:\dG \to \sdGb$ form an adjoint pair.
\item{\sl(b)}
If the ultraproduct HFI approximations $\eta\:G \to \sGb$,
\,$\phi\:\dG \to \sdGb$ form an adjoint pair, then there is
a function $\tau\: I \cx \bN \to I$ such that
$i \le \tau(i,n) \le \tau(i,n+1)$ for all $i \in I$, $n \in \bN$,
the approximating systems
$\bigl(\eta_{\tau(i,n)}\:G_{\tau(i,n)} \to \Gb\bigr)_{(i,n) \in I \cx \bN}$,
$\bigl(\phi_{\tau(i,n)}\:\dG_{\tau(i,n)} \to \dGb\bigr)_{(i,n) \in I \cx \bN}$
are well based with respect to the $\alfa$-adjoint DD bases
$\bigl((\Kb_{in},\Ub_{in})\bigr)_{(i,n)}$,
$\bigl((\Gamb_{in},\Omb_{in})\bigr)_{(i,n)}$,
respectively, and they form an adjoint pair.
\item{}}
\vskip-12pt
\endproclaim

\demo{Proof}
Let us form the corresponding IMG group triplets $(G,\Go,\Gf)$,
$\bigl(\dG,\Ho,\Hf\bigr)$, corresponding to the HFI approximations
$\eta\:G \to \sGb$, \,$\phi\:\dG \to \sdGb$.

(a) Assume that the approximating systems
$(\eta_i)_{i\in I}$, $(\phi_i)_{i \in I}$ form an adjoint pair.
Then, as easily seen, condition~1 from the definition of adjoint
pair of approximating systems implies the relation
$(\phi\,\gama)(\eta\,x) \apr \gama(x)$, for $x \in \Gf$,
$\gama \in \Hf$. Condition~2 implies the inclusion
$\Go^{\ort} \sbs \Hf$. The argument can be completed by
referring to Lemma~2.5.3.

(b) As first we realize that, for each $i \in I$, $\eta$ is a
$(\sKb_i,\sUb_i)$ approximation of \,$\sGb$ and $\phi$ is
a $({}^*\Gamb_i,{}^*\Omb_i)$ approximation of \,$\sdGb$. This is
to say that there is a set $J^i_1 \in \DD$ such that $\eta_j$ is
a $(\Kb_i,\Ub_i)$ approximation of \,$\Gb$ and $\phi_j$ is a
$(\Gamb_i,\Omb_i)$ approximation of \,$\dGb$ for each $j \in J^i_1$.

Second, we have $\Ho = \Gf^{\ort}$ and $\Hf = \Go^{\ort}$, thus in
particular,
$$
\bigcup_{i \in I} \Bohr_\alfa\bigl(\eta^{-1}[\sUb_i]\bigr)
   = \Go^{\ort} \sbs \Hf =
        \bigcup_{k \in I} \phi^{-1}[{}^*\Gamb_k]\,.
$$
Hence there a function $\sig\:I \to I$ such that
$$
\Bohr_\alfa\bigl(\eta^{-1}[\sUb_i]\bigr)
   \sbs \phi^{-1}[{}^*\Gamb_k]\,,
$$
for each $i \in I$ and $k \ge \sig(i)$. This means that
there is a set $J^i_2 \in \DD$ such that, for any $i \in I$,
$k \ge \sig(i)$ and $j \in J^i_2$,
$$
\Bohr_\alfa\bigl(\eta_j^{-1}[\Ub_i]\bigr)
   \sbs \phi_j^{-1}[\Gamb_k]\,.
$$

Finally, $(\phi\,\gama)(\eta\,x) \apr \gama(x)$, for $x \in \Gf$,
$\gama \in \Hf$. Let us fix any sequence $(\eps_n)_{n\in\bN}$, such
that $0 < \eps_{n+1} < \eps_n < 2\pi/3$ for each $n$, converging
to~$0$. Now, choosing any $i \in I$, $n \in \bN$, both the sets
$\Kb_i \sbs \Gb$, $\Gamb_i \sbs \dGb$ are compact, hence
$\eta^{-1}[\sKb] \sbs \Gf$ and $\phi^{-1}[{}^*\Gamb] \sbs \Hf$,
therefore,
$$
\lv \arg\frac{(\phi \gama)(\eta x)}{\gama(x)}\rv
   \le \eps_n\,,
$$
for all $x \in \eta^{-1}[\sKb_i]$, $\gama \in \phi^{-1}[{}^*\Gamb_i]$.
Thus there is a set $J^{in}_3 \in \DD$ such that
$$
\lv \arg\frac{(\phi_j \gama)(\eta_j x)}{\gama(x)}\rv
   \le \eps_n\,,
$$
for any $x \in \eta_j^{-1}[\Kb_i]$, $\gama \in \phi_j^{-1}[\Gamb_i]$,
whenever $j \in J^{in}_3$.

Then the set $J^{in} = J^i_1 \cap J^i_2 \cap J^{in}_3$ belongs to
$\DD$, hence it is cofinal in $(I,\le)$. Thus there is a function
$\tau\:I \cx \bN \to I$ such that $\tau(i) \in J^{in}$,
$i \le \tau(i,n) \le $ and $\sig(i) \le \tau(i,n)$, for $i \in I$,
$n \in \bN$. We can arrange that $\tau$ satisfies the additional
inequality $\tau(i,n) \le \tau(i,n+1)$ by an induction argument.

Now, it is clear that the approximating systems
$\bigl(\eta_{\tau(i,n)}\:G_{\tau(i,n)} \to \Gb\bigr)_{(i,n) \in I \cx \bN}$,
\,$\bigl(\phi_{\tau(i,n)}\:\dG_{\tau(i,n)} \to \dGb\bigr)_{(i,n) \in I \cx \bN}$
of the LCA groups $\Gb$, $\dGb$ are well based with respect to their
$\alfa$-adjoint DD~bases $\bigl((\Kb_{in},\Ub_{in})\bigr)_{(i,n)}$,
$\bigl((\Gamb_{in},\Omb_{in})\bigr)_{(i,n)}$, respectively, and it is
routine to check that they form an adjoint pair.
\enddemo

\remark{Remark}
(a) It is clear from the proof that Proposition~2.5.5 would remain
true if we changed the definition of dual approximating systems
by replacing its second condition by an analogous condition
corresponding to any of the three inclusions
\,$\Ho^{\ort} \sbs \Gf$, \,$\Hf^{\ort} \sbs \Go$,
\,$\Gf^{\ort} \sbs \Ho$, equivalent to \,$\Go^{\ort}\sbs \Hf$
by Lemma~2.5.3. The reader is invited to formulate them as an
exercise.

(b) If there is a net $(\eps_i)_{i \in I}$ over $(I,\le)$ of positive
reals $\eps_i < 2\pi/3$, converging to~0, in particular, if \,$I = \bN$,
then the function $\tau\: I \cx \bN \to I$ in (b) can be replaced by a
function $\tau\: I \to I$, giving rise to an adjoint pair of approximating
systems \hbox{$\bigl(\eta_{\tau(i)}\:G_i \to \Gb\bigr)_{i \in I}$,}
$\bigl(\phi_{\tau(i)}\:\dG_i \to \dGb\bigr)_{i \in I}$, well based
with respect to the original $\alfa$-adjoint pair of DD bases
$\bigl((\Kb_i,\Ub_i)\bigr)_{i \in I}$,
$\bigl((\Gamb_i,\Omb_i)\bigr)_{i \in I}$. Analogous remarks apply
to (the proof of) Theorem~2.5.6 and to Theorem~2.5.9 to follow.
\endremark
\medskip

Proposition~1.5.7, Theorem~1.5.11, Corollary~2.5.2 and Proposition~2.5.5
yield the following consequence, giving more precision to Theorem~3
stated without proof in the Introduction to \cite{Go2}. However, in order
to apply the above mentioned results we need to guarantee that the nonstandard
universe obtained by means of the ultraproduct construction is sufficiently
saturated. To this end we need to impose some additional conditions on the
directed index poset $(I,\le)$ and the directed ultrafilter $\DD$
on $I$.

A directed poset $(I,\le)$ is called {\it regular\/} if it is infinite
and all its subsets $[i)$ have the same cardinality as $I$ itself. $(I,\le)$
is called {\it countably incomplete\/} if there is some strictly increasing
sequence $(k_n)_{n \in \bN}$ of elements of \,$I$ such that
$\bigcap_n [k_n) = \emptyset$. The term {\it regular countably incomplete
directed poset\/} will be abbreviated to {\it RCID poset}.
A thorough inspection of the proof of Theorem~6.1.4 in Chang-Keisler's book
\cite{CK} shows that there exist $(\card I)^+$-good directed ultrafilters
on every RCID poset $(I,\le)$. Any such an ultrafilter is automatically
countably incomplete and the nonstandard universe obtained by the
ultraproduct construction with respect to it is $(\card I)^+$-saturated.

Given any DD base $\bigl((\Kb_i,\Ub_i)\bigr)_{i \in I}$ over some directed
poset $(I,\le)$, we can assume, without loss of generality, that
$(I,\le)$ is an RCID poset. Otherwise we can put $J = I \cx I \cx \bN$
and order $J$ componentwise, i.e., $(i_1,j_1,n_1) \le (i_2,j_2,n_2)$
if and only if \,$i_1 \le i_2$, $j_1 \le j_2$ and $n_1 \le n_2$. Further
we can define $\Kb_{ijn} = \Kb_i$ and $\Ub_{ijn} = \Ub_i$ for $i,j \in I$,
$n \in \bN$. It is clear that $(J,\le)$ is an RCID poset and
$\bigl((\Kb_{ijn},\Ub_{ijn})\bigr)_{(i,j,n) \in J}$ is a DD base over
$(J,\le)$, again.

\proclaim{2.5.6.~Theorem}
Let \,$\Gb$ be a Hausdorff LCA group and \,$0 < \alfa < 2\pi/3$.
Then $\Gb$, $\dGb$ admit a pair of \,$\alfa$-adjoint DD bases over
some directed poset $(J,\le)$ and an adjoint pair of approximating
systems $\bigl(\zeta_j\:G_j \to \Gb\bigr)_{j \in J}$,
\,$\bigl(\psi_j\:\dG_j \to \dGb\bigr)_{j \in J}$ by finite abelian groups,
well based with respect to these DD~bases. Moreover, it is possible
to arrange that at least one of these approximating systems is injective.
\endproclaim

\demo{Proof}
Let $\bigl((\Kb_i,\Ub_i)\bigr)_{i \in I}$,
$\bigl((\Gamb_i,\Omb_i)\bigr)_{i \in I}$ be any pair of
$\alfa$-adjoint DD bases of \,$\Gb$ and \,$\dGb$, respectively,
over an RCID poset $(I,\le)$. Let
$(\eta_i\:G_i \to \Gb)_{i \in I}$ be an approximating
system of \,$\Gb$, well based with respect to the DD base
$\bigl((\Kb_i,\Ub_i)\bigr)_{i \in I}$, and $\DD$ be a $\kappa$-good
directed ultrafilter on $(I,\le)$, where $\kappa = (\card I)^+$. Let us
form the ultraproduct HFI approximation $\eta = (\eta_i)/\DD\:G \to \sGb$.
within the $\kappa$-saturated nonstandard universe obtained via the
ultraproduct construction modulo~$\DD$. Then, by Corollary~2.5.2,
$\eta$ has a dual approximation $\phi = (\phi_i)/\DD\:\dG \to \sdGb$.
Then the DD bases $\bigl((\Kb_{in},\Ub_{in})\bigr)_{(i,n)}$,
$\bigl((\Gamb_{in},\Omb_{in})\bigr)_{(i,n)}$, described prior to
Proposition~2.5.5, and the approximating systems
$\bigl(\eta_{\tau(i,n)}\:G_{\tau(i,n)} \to \Gb\bigr)_{(i,n)}$,
$\bigl(\phi_{\tau(i,n)}\:\dG_{\tau(i,n)} \to \dGb\bigr)_{(i,n)}$
over the directed poset $(I \cx \bN, \le)$, constructed in
its proof, have all the properties required.
\enddemo

So far, so good. However, it is neither the HFI approximations nor
the approximating systems but always a single ``sufficiently good''
pair of finite approximations which is decisive for applications.
Let us examine ``how good'' adjoint pairs of finite approximations
can be guaranteed to exist.

Let $\Gb$ be a Hausdorff LCA group and $0 < \eps < \alfa \le \pi/3$.%
\footnote{The only reason for this restriction is to guarantee that
$\alfa - \eps > 0$ and $\alfa + \eps < 2\pi/3$.}
Let further $\Kb \sbs \Gb$, $\Gamb \sbs \dGb$ be compact sets and
$\Ub \sbs \Gb$, $\Omb \sbs \dGb$ be neighborhoods of the neutral
elements in $\Gb$, $\dGb$, respectively, such that the pairs
$(\Kb,\Ub)$,  $(\Gamb,\Omb)$ are $\alfa$-adjoint. Given a finite
abelian group $G$ and mappings $\eta\:G \to \Gb$,
$\phi\:\dG \to \dGb$, we say that $\eta$, $\phi$ form a
{\it strongly $(\alfa,\eps)$-adjoint pair approximations\/} of
\,$\Gb$, $\dGb$, respectively, with respect to $(\Kb,\Ub)$,
$(\Gamb,\Omb)$ if
$$
\lv \arg\frac{(\phi\,\gama)(\eta\,x)}{\gama(x)}\rv \le \eps
$$
for all $x \in \eta^{-1}[\Kb]$, $\gama \in \phi^{-1}[\Gamb]$,
and there exist neighborhoods $\Vb \sbs \Ub$, $\Ypb \sbs \Omb$
of the neutral elements in $0 \in \Gb$, $1 \in \dGb$, respectively,
satisfying the inclusions
$$\align
\Bohr_\alfa\bigl(\eta^{-1}[\Ub]\bigr) &\sbs
\phi^{-1}\bigl[\Bohr_\eps(\Vb)\bigr]\,, \\
\Bohr_\alfa\bigl(\phi^{-1}[\Omb]\bigr) &\sbs
\eta^{-1}\bigl[\Bohr_\eps(\Ypb)\bigr]\,,
\endalign$$
such that $\eta$ is a $(\Kb,\Vb)$ approximation of \,$\Gb$ and
$\phi$ is a $(\Gamb,\Ypb)$ approximation of \,$\dGb$. We also say
that the strong $(\alfa,\eps)$-adjointness of $\eta$, $\phi$ is
{\it witnessed\/} by the sets $\Vb$, $\Ypb$.

Let us remark that the left hand side expressions denote Bohr sets
in the finite abelian groups $\dG$, $G$, while the right hand ones
are Bohr sets in the LCA groups $\dGb$, $\Gb$, respectively.

\proclaim{2.5.7.~Lemma}
Let \,$\Gb$, $\alfa$, $\eps$, $(\Kb,\Ub)$, $(\Gamb,\Omb)$ and
\,$G$ be as above. Let the mappings \,$\eta\:G \to \Gb$,
$\phi\:\dG \to \dGb$ form a strongly $(\alfa,\eps)$-adjoint pair of
approximations of \,$\Gb$, $\dGb$, respectively, with  respect to
$(\Kb,\Ub)$, $(\Gamb,\Omb)$, witnessed by the sets $\Vb \sbs \Ub$,
$\Ypb \sbs \Omb$. Then for any sets $\Xb \sbs \Gb$,
$\Delb \sbs \dGb$, such that \,$\Ub \sbs \Xb \sbs \Kb$,
\,$\Omb \sbs \Delb \sbs \Gamb$, we have
$$\gather
\phi^{-1}\bigl[\Bohr_{\alfa-\eps}(\Xb)\bigr] \sbs
\Bohr_\alfa\bigl(\eta^{-1}[\Xb]\bigr)\,,\\
\eta^{-1}\bigl[\Bohr_{\alfa-\eps}(\Delb)\bigr] \sbs
\Bohr_\alfa\bigl(\phi^{-1}[\Delb]\bigr)\,.
\endgather$$
If, additionally, $\Xb + \Vb \sbs \Kb$, \,$\Delb\,\Ypb \sbs \Gamb$,
then also
$$\gather
\Bohr_{\alfa-\eps}\bigl(\eta^{-1}[\Xb + \Vb]\bigr) \sbs
\phi^{-1}\bigl[\Bohr_{\alfa+\eps}(\Xb)\bigr]\,,\\
\Bohr_{\alfa-\eps}\bigl(\phi^{-1}[\Delb\,\Ypb]\bigr) \sbs
\eta^{-1}\bigl[\Bohr_{\alfa+\eps}(\Delb)\bigr]\,.
\endgather$$
\endproclaim

\demo{Proof}
Let us prove the first inclusion; then the second will follow by
a symmetry argument. From $\Ub \sbs \Xb \sbs \Kb$ we deduce that
$$
\Bohr_{\alfa-\eps}(\Xb) \sbs \Bohr_{\alfa-\eps}(\Ub)
   \sbs \Bohr_\alfa(\Ub) = \Gamb\,,
$$
hence
$\lv\arg\bigl((\phi\,\gama)(\eta\,x)/\gama(x)\bigr)\rv \le \eps$
for all $x \in \eta^{-1}[\Xb]$,
$\gama \in \phi^{-1}\bigl[\Bohr_{\alfa-\eps}(\Xb)\bigr]$, and
$$
\lv\arg\gama(x)\rv \le
   \lv\arg(\phi\,\gama)(\eta\,x)\rv +
      \lv\arg\frac{\gama(x)}{(\phi\,\gama)(\eta\,x)}\rv
         \le (\alfa - \eps) + \eps = \alfa\,.
$$

Let us assume that $\Xb + \Vb \sbs \Kb$ and prove the first
inclusion in the ``additional'' part of the Lemma. The inclusions
$$
\Bohr_{\alfa-\eps}\bigl(\eta^{-1}[\Xb + \Vb]\bigr)
\sbs \Bohr_\alfa\bigl(\eta^{-1}[\Ub]) \sbs
\phi^{-1}\bigl[\Bohr_\eps(\Vb)\bigr]
$$
show that $\lv\arg(\phi\,\gama)(\vb)\rv \le \eps$ for any
$\gama \in \Bohr_{\alfa-\eps}\bigl(\eta^{-1}[\Xb + \Vb]\bigr)$,
$\vb \in \Vb$. For such a $\gama$, we are to show that
$\lv\arg(\phi\,\gama)(\xb)\rv \le \alfa + \eps$, for any
$\xb \in \Xb$. As $\eta$ is a $(\Kb,\Vb)$ approximation of \,$\Gb$,
there is an $x \in G$ such that $\eta(x) = \xb + \vb$ for some
$\vb \in \Vb$; then $x \in \eta^{-1}[\Xb + \Vb]$. It follows that
$(\phi\,\gama)(\xb) = (\phi\,\gama)(\eta\,x)/(\phi\,\gama)(\vb)$,
hence
$$\align
\lv\arg(\phi\,\gama)(\xb)\rv
&\le \lv\arg\frac{(\phi\,\gama)(\eta\,x)}{\gama(x)}\rv
+ \lv\arg\gama(x)\rv + \lv\arg(\phi\,\gama)(\vb)\rv\\
&\le \eps + (\alfa - \eps)  + \eps = \alfa + \eps\,.
\endalign$$
The second ``additional'' implication follows from the first by
symmetry of the situation.
\enddemo

\proclaim%
{2.5.8.~Strongly Adjoint Finite LCA Group Approximation Theorem}
Let \,$\Gb$ be a Hausdorff LCA group and
\,$0 < \eps < \alfa \le \pi/3$. Let further $\Kb \sbs \Gb$,
$\Gamb \sbs \dGb$ be compact sets, $\Ub \sbs \Gb$,
$\Omb \sbs \dGb$ be neighborhoods of the neutral elements in $\Gb$,
$\dGb$, respectively, such that the pairs $(\Kb,\Ub)$,
$(\Gamb,\Omb)$ are $\alfa$-adjoint. Then there exist a finite
abelian group $G$ and mappings $\eta\:G \to \Gb$,
$\phi\:\dG \to \dGb$ forming a strongly $(\alfa,\eps)$-adjoint
pair of approximations of \,$\Gb$, $\dGb$ with respect to
$(\Kb,\Ub)$, \,$(\Gamb,\Omb)$, respectively. One can arrange
additionally that at least one of the approximations $\eta$,
$\phi$ is injective.
\endproclaim

\demo{Proof}
Let $G$ be a hyperfinite abelian group in a sufficiently saturated
nonstandard universe, $\zeta\:G \to \sGb$, $\psi\:\dG \to \sdGb$ be
an adjoint pair of HFI approximations of \,$\Gb$, $\dGb$, and
$(G,\Go,\Gf)$, $\bigl(\dG,\Ho,\Hf\bigr)$ be the corresponding IMG
group triplets, arising from $\zeta$, $\psi$, respectively. Then,
for any internal neighborhoods $\Wb \sbs \IG$, \,$\Thb \sbs \IdG$
of the neutral elements in $\sGb$, \,$\sdGb$, respectively, we have
$$\gather
\Bohr_\alfa\bigl(\zeta^{-1}[\sUb]\bigr) \sbs \Go^{\ort} = \Hf
= \psi^{-1}\bigl[\FdG\bigr] = \psi^{-1}\bigl[\IG^{\ort}\bigr]
\sbs \psi^{-1}\bigl[\Bohr_\eps(\Wb)\bigr] \,,\\
\Bohr_\alfa\bigl(\psi^{-1}[{}^*\Omb]\bigr) \sbs \Ho^{\ort} = \Gf
= \zeta^{-1}[\FG] = \zeta^{-1}\bigl[\IdG{}^{\ort}\bigr]
\sbs \zeta^{-1}\bigl[\Bohr_\eps(\Thb)\bigr]\,.
\endgather$$
Then the pairs $(\sKb,\sUb)$, $({}^*\Gamb,{}^*\Omb)$ are
$\alfa$-adjoint, and, by saturation, there is some pair of internal
neighborhoods $\Wb \sbs \IG$, $\Thb \sbs \IdG$ of \,$0 \in \sGb$
and $1 \in \sdGb$ such that $\zeta$ is a $(\sKb,\Wb)$ approximation of
\,$\sGb$ and $\psi$ is a $({}^*\Gamb,\Thb)$ approximation of \,$\sdGb$.
Now it is clear, that $\zeta$, $\psi$ form a strongly
$(\alfa,\eps)$-adjoint pair of internal approximations, witnessed by
$\Wb$, $\Thb$. By the {\it transfer principle\/} there exists
a strongly $(\alfa,\eps)$-adjoint pair of finite approximations
of \,$\Gb$, $\dGb$, respectively, with respect to $(\Kb,\Ub)$,
$(\Gamb,\Omb)$, witnessed by some neighborhoods $\Vb \sbs \Ub$,
$\Ypb \sbs \Omb$ of \,$0 \in \Gb$ and $1 \in \dGb$, respectively.
If at least one of the starting HFI approximations $\zeta$, $\psi$
were injective, then one can guarantee the same property for the
corresponding finite approximation $\eta$, $\phi$, respectively,
as well.
\enddemo

Theorem~2.5.8 implies the following strengthening of Theorem~2.5.6.
In its formulation we use the last convention formulated prior to
Proposition~2.5.5.

\proclaim{2.5.9.~Theorem}
Let \,$\Gb$ be a Hausdorff LCA group, $0 < \alfa \le \pi/3$,
$\bigl((\Kb_i,\Ub_i)\bigr)_{i \in I}$,
$\bigl((\Gamb_i,\Omb_i)\bigr)_{i \in I}$ be an $\alfa$-adjoint
pair of DD bases of \,$\Gb$ and its dual \,$\dGb$, respectively,
over a directed poset $(I,\le)$, and $(\eps_n)_{n \in \bN}$ be a
decreasing sequence of positive reals $< \alfa$, converging to~$0$.
Then there exists an adjoint pair of approximating systems
\hbox{$\bigl(\eta_{in}\:G_{in} \to \Gb\bigr)_{(i,n) \in I \cx \bN}$,}
\,$\bigl(\phi_{in}\:\dG_{in} \to \dGb\bigr)_{(i,n) \in I \cx \bN}$
of \,$\Gb$ and $\dGb$, respectively, by finite abelian groups, well
based with respect to their DD bases
$\bigl((\Kb_{in},\Ub_{in})\bigr)_{(i,n) \in I \cx \bN}$,
$\bigl((\Gamb_{in},\Omb_{in})\bigr)_{(i,n) \in I \cx \bN}$,
such that each particular pair $\eta_{in}\:G_{in} \to \Gb$,
\,$\phi_{in}\:\dG_{in} \to \dGb$ is strongly
$(\alfa,\eps_n)$-adjoint with respect to $(\Kb_i,\Ub_i)$,
$(\Gamb_i,\Omb_i)$. Moreover, it is possible to arrange that
at least one of these approximating systems is injective.
\endproclaim

\demo{Proof}
Let us pick any $i \in I$, $n \in \bN$. According to Theorem~2.5.8,
there are a $(\Kb_i,\Ub_i)$ approximation $\eta_{in}\: G_{in} \to \Gb$
and a $(\Gamb_i,\Omb_i)$ approximation $\phi_{in}\:\dG_{in} \to \dGb$,
forming a strongly $(\alfa,\eps_n)$-adjoint pair witnessed
by some neighborhoods $\Vb_{in} \sbs \Ub_i$, $\Ypb_{in} \sbs \Omb_i$
of \,$0 \in \Gb$, $1 \in \dGb$, respectively. Then the systems of
mappings $(\eta_{in})_{(i,n)}$, $(\phi_{in})_{(i,n)}$ satisfy all
the requirements of the Theorem. The supplement on injectivity is
already notorious.
\enddemo

The ``HFI parts'' of the following three examples are due to Gordon
\cite{Go2}; we are adding  the standard counterparts mainly with
the aim to illustrate the strong adjointness phenomenon. As we shall
see in these fairly important cases, some conditions of the strong
$(\alfa,\eps)$-adjointness are satisfied, so to say, automatically
and in a stricter form than required by the definition.

In the first example, building on items (a), (b) of Example~1.5.9,
we construct certain pairs of adjoint approximations for the
group $\bT$ and its dual group $\wh\bT \iso \bZ$, where
$\gama(\xb) = \xb^\gama$ for $\xb \in \bT$, $\gama \in \bZ$.
Similarly, the (hyper)finite cyclic group $\bZ_n$ is identified
with its dual group $\wh\bZ_n$ via the pairing
$\gama(a) = \ee^{2\pi\ii a\gama/n}$, for $a,\gama \in \bZ_n$.
Let us recall that $\bZ_n$ is still represented as the group
of absolutely smallest remainders modulo~$n$.

\example{2.5.10.~Example}
If \,$n \in \sbNi$, then the internal homomorphism
$\eta\:\bZ_n \to \sbT$, $\eta(a) = \ee^{2\pi\ii a/n}$, and
the inclusion map $\phi\:\bZ_n \to \sbZ$ give rise to the
mutually dual IMG group triplets $(\bZ_n,\Go,\bZ_n)$,
$\bigl(\bZ_n,\{0\},\bZ\bigr)$, with normalizing multipliers
$d = n^{-1}$, $\hat d = 1$, respectively, where
$$
\Go = \Bigl\{a \in \bZ_n\: \frac{a}{n} \apr 0\Bigr\}
    = \bZ^{\ort}\,.
$$
For $a,\gama \in \bZ_n$ we even have
$$
(\phi\,\gama)(\eta\,a) =
   \Bigl(\ee^{\frac{2\pi\ii a}{n}}\Bigr)^{\gama}
      = \ee^{\frac{2\pi\ii a\gama}{n}} = \gama(a)\,,
$$
while the infinitesimal nearness \,$\apr$ \,in place of the
second equality and $\gama \in \bZ$ would suffice to establish
the adjointness of the HFI approximations $\eta$, $\phi$ of
the mutually dual groups $\bT$, $\bZ$, respectively.

Passing to the standard situation, let $n, k \in \bN$ and
$0 < r \le \alfa \le \pi/3$ be such that $n > \pi/r$ and
$1 \le k < n/4$. We put
$$
\Ub = \bigl\{\xb \in \bT\: \lv\arg\xb\rv \le r\bigr\}\,,
\qquad\qquad
\Gama = \bigl\{\gama \in \bZ\: \lv\gama\rv \le k\bigr\}\,.
$$
Then we have
$$
\Bohr_\alfa(\Ub) =
\Bigl\{\gama \in \bZ\: \lv\gama\rv \le \frac{\alfa}{r}\Bigr\}\,,
\qquad\qquad
\Bohr_\alfa(\Gama) =
\Bigl\{\xb \in \bT\: \lv\arg\xb\rv \le \frac{\alfa}{k}\Bigr\}\,,
$$
hence the pair $(\bT,\Ub)$ in $\bT$ and the pair
$\bigl(\Gama,\{0\}\bigr)$ in $\bZ$ are $\alfa$-adjoint
if and only if (as the remaining conditions are trivial)
$\Gama = \Bohr_\alfa(\Ub)$, which is equivalent to
$$
k = \Bigl\lfloor\frac{\alfa}{r}\Bigr\rfloor\,.
$$
If \,$r = \alfa/k$, then also $\Bohr_\alfa(\Gama) = \Ub$.
The $(\bT,\Ub)$ approximation $\eta\:\bZ_n \to \bT$ of \,$\bT$ and
the $\bigl(\Gama,\{0\}\bigr)$ approximation $\phi\:\bZ_n \to \bZ$
of \,$\bZ$ (both given as in the HFI case) still satisfy
$(\phi\,\gama)(\eta\,a) = \gama(a)$ for $a, \gama \in \bZ_n$.
In order to show that they are strongly $(\alfa,\eps)$-adjoint for
an $\eps \in (0,\alfa)$ with respect to the pairs $(\bT,\Ub)$,
$\bigl(\Gama,\{0\}\bigr)$ we need to point out some witnessing sets.
$\Yps = \{0\} \sbs \bZ$ is plain. A~simple computation shows
that
$$
\Bohr_\alfa\bigl(\eta^{-1}[\Ub]\bigr) =
\Bigr\{\gama \in \bZ_n\:
   \lv\gama\rv \le \frac{\alfa}{r}\Bigl\}\,,
$$
and for $\Vb = \{\xb \in \bT\: \lv\arg\xb\rv \le s\}$, where
$0 < s \le r$, we have
$$
\phi^{-1}\bigl[\Bohr_\eps(\Vb)\bigr] =
\Bigl\{\gama \in \bZ_n\: \lv\gama\rv \le \frac{\eps}{s}\Bigr\}\,.
$$
Thus the inclusion $\Bohr_\alfa\bigl(\eta^{-1}[\Ub]\bigr) \sbs
\phi^{-1}\bigl[\Bohr_\eps(\Vb)\bigr]$ can be guaranteed by choosing
any positive
$$
s \le \frac{r\,\eps}{\alfa}\,.
$$
In order $\eta$ to be a $(\bT,\Vb)$ approximation of \,$\bT$ we need
that $n > \pi/s$.
\endexample

In the next example we describe a pair of adjoint approximations for
the self-dual group $\bR$, building on item (c) of Example~1.5.9.
More precisely, $\bR$ is identified with its dual group $\wh\bR$
via the pairing $\ggam(\xb) = \ee^{\ii\xb\ggam}$, for
$\xb,\ggam \in \bR$. The passage to differently scaled pairings
$\ggam(\xb) = \ee^{2\pi\ii\xb\ggam/T}$, with any $T > 0$, is
straightforward.

\example{2.5.11.~Example}
Let \,$n \in \sbNi$, and $d,d'$ be positive infinitesimals such
that both the hyperreal numbers $n d$, $n d'$ are infinite. Then
the internal mappings $\eta\:\bZ_n \to \sbR$, $\eta(a) = ad$, and
$\phi\:\bZ_n \to \sbR$, $\phi(\gama) = \gama d'$, are HFI
approximations of the group $\bR$, inducing IMG group triplets
$(\bZ_n,\Go,\Gf)$, $(\bZ_n,\Ho,\Hf)$, where
$$\xalignat{2}
\Go &= \{a \in \bZ_n\: ad \apr 0\}\,,&
\Gf &= \{a \in \bZ_n\: \lv a\rv d < \infty\}\,, \\
\Ho &= \{\gama \in \bZ_n\: \gama d' \apr 0\}\,,&
\Hf &= \{\gama \in \bZ_n\: \lv\gama\rv d' < \infty\}\,,
\endxalignat$$
with normalizing multipliers $d$, $d'$, respectively. One can
easily verify that these triplets are mutually dual if and only if
\,$n d d' \in \FR \sms \IR$. For any $a,\gama \in \bZ_n$ we have
$$
(\phi\,\gama)(\eta\,a) = \ee^{\ii a\gama dd'}\,,
\qquad\qquad
\gama(a) = \ee^{\frac{2\pi \ii a\gama}{n}}\,,
$$
and the two expressions are infinitesimally close for all
$a \in \Gf$, $\gama \in \Hf$ if and only if
$$
n dd' \apr 2\pi\,.
$$
Hence this condition is equivalent to the adjointness of the
HFI approximations $\eta$, $\phi$. Then the scaling coefficient
$\hat d = 1/nd$, dual to $d$, is a normalizing multiplier for the
triplet $(\bZ_n,\Ho,\Hf)$, as well. Under the particular choice
$$
d' = 2\pi\hat d = \frac{2\pi}{nd}
$$
we even have $(\phi\,\gama)(\eta\,a) = \gama(a)$ for all
$a,\gama \in \bZ_n$.

The corresponding standard situation is framed by some
$n,k,m \in \bN$ and positive $d,d',r,\rho,\alfa \in \bR$, such
that $1 \le k,m < n/4$, $r,\rho \le \alfa \le \pi/3$, and
$d/2 < r \le kd$, $d'2 < \rho \le md'$.
We put
$$\xalignat{2}
\Ub &= [-r,r]\,, & \Kb &= [-kd,kd]\,, \\
\Omb &= [-\rho,\rho]\,, & \Gamb &= [-md',md']\,.
\endxalignat$$
Then
$$\xalignat{2}
\Bohr_\alfa(\Ub) &=
\Bigl[-\frac{\alfa}{r},\frac{\alfa}{r}\Bigr]\,, &
\Bohr_\alfa(\Kb) &=
\Bigl[-\frac{\alfa}{kd},\frac{\alfa}{kd}\Bigr]\,, \\
\Bohr_\alfa(\Omb) &=
\Bigl[-\frac{\alfa}{\rho},\frac{\alfa}{\rho}\Bigr]\,, &
\Bohr_\alfa(\Gamb) &=
\Bigl[-\frac{\alfa}{md'},\frac{\alfa}{md'}\Bigr]\,,
\endxalignat$$
thus the pairs $(\Kb,\Ub)$, $(\Gamb,\Omb)$ are $\alfa$-adjoint if
and only if
$$
\rho \cd kd = r \cd md' = \alfa\,.
$$
The reverse proportionality of the lengths $kd$ and $\rho$,
as well as that of $md'$ and $r$ are worthwhile to notice.
The $(\Kb,\Ub)$ approximation $\eta\:\bZ_n \to \bR$ and the
$(\Gamb,\Omb)$ approximation $\phi\:\bZ_n \to \bR$ of \,$\bR$
are given by the same formulas as in the HFI case. Then
$$
\lv\arg\frac{(\phi\,\gama)(\eta\,a)}{\gama(a)}\rv =
\lv a\gama\rv dd'\lv 1 - \frac{2\pi}{ndd'}\rv\,.
$$
This expression is \,$\le \eps$ for an $\eps \in (0,\alfa)$ and all
$a \in \eta^{-1}[\Kb]$, $\gama \in \phi^{-1}[\Gamb]$ if and only if
$$
kd\,md' \biggl(1 - \frac{2\pi}{ndd'}\biggr) \le \eps\,.
$$
Putting $d' = 2\pi/nd$ we can even achieve that
$(\phi\,\gama)(\eta\,a) = \gama(a)$ for all $a,\gama \in \bZ_n$.
It remains to describe some intervals
$\Vb = [-s,s]$, $\Ypb = [-\sig,\sig]$, where $0 < s \le r$,
$0 < \sig \le \rho$, witnessing the strong
$(\alfa,\eps)$-adjointness of the approximations $\eta$, $\phi$
with respect to the pairs $(\Kb,\Ub)$, $(\Gamb,\Omb)$. Since
$$\xalignat{2}
\Bohr_\alfa\bigl(\eta^{-1}[\Ub]\bigr) &=
\Bigr\{\gama \in \bZ_n\:
   \lv\gama\rv \le \frac{nd\alfa}{2\pi r}\Bigl\}\,,&
\phi^{-1}\bigl[\Bohr_\eps(\Vb)\bigr] &=
\Bigl\{\gama \in \bZ_n\: \lv\gama\rv \le \frac{\eps}{sd'}\Bigr\}\,,\\
\Bohr_\alfa\bigl(\eta^{-1}[\Omb]\bigr) &=
\Bigr\{a \in \bZ_n\:
   \lv a\rv \le \frac{nd'\alfa}{2\pi\rho}\Bigl\}\,,&
\phi^{-1}\bigl[\Bohr_\eps(\Ypb)\bigr] &=
\Bigl\{a \in \bZ_n\: \lv a\rv \le \frac{\eps}{\sig d}\Bigr\}\,,
\endxalignat$$
the inclusions
\,$\Bohr_\alfa\bigl(\eta^{-1}[\Ub]\bigr) \sbs
\phi^{-1}\bigl[\Bohr_\eps(\Vb)\bigr]$,
\,$\Bohr_\alfa\bigl(\eta^{-1}[\Omb]\bigr) \sbs
\phi^{-1}\bigl[\Bohr_\eps(\Ypb)\bigr]$ are equivalent to the
inequalities
$$
s \le \frac{2\pi\,r\,\eps}{ndd'\alfa}\,,
\qquad\qquad
\sig \le \frac{2\pi\,\rho\,\eps}{ndd'\alfa}\,,
$$
respectively. If \,$ndd' = 2\pi$, then they take the simple form
$$
s \le \frac{r\,\eps}{\alfa}\,,
\qquad\qquad
\sig \le \frac{\rho\,\eps}{\alfa}\,,
$$
analogous to that in Example~2.5.10. In order $\eta$ to be
a $(\Kb,\Vb)$ approximation and $\phi$ to be a $(\Gamb,\Ypb)$
approximation of \,$\bR$ we need that $d < 2s$ and
$d' < 2\sig$. As $\rho \cd kd = r \cd md' = \alfa$, this implies
$$
k > \frac{\alfa}{2\rho s}\,,
\qquad\text{and}\qquad
m > \frac{\alfa}{2r\sig}\,,
$$
which is possible only if
$$
n > \frac{2\alfa}{\min\{\rho s, r\sig\}}\,.
$$
\endexample

In the last of our examples we construct adjoint approximations
to those from item (d) of Example~1.5.9.

\example{2.5.12.~Example}
Let $\Gb$ be a Hausdorff LCA group with a DD base
$\bigl((\Kb_i,\Ub_i)\bigr)_{i \in I}$ consisting of compact
open subgroups $\Ub_i \sbs \Kb_i$ of \,$\Gb$. As
$\Bohr_\alfa(\Xb)$ coincides with the annihilator
$\Xb^\perp = \Bohr_0(\Xb)$ for $\alfa \in (0,2\pi/3)$ and
any subgroup $\Xb \sbs \Gb$, the system
$\bigl(\bigl(\Ub_i^\perp,\Kb_i^\perp\bigr)\bigr)_{i \in I}$
forms a DD base of the dual group $\dGb$ consisting of
compact open subgroups, again. Picking an $i \in I$ and
putting $\Ub = \Ub_i$, $\Kb = \Kb_i$, and $\Omb = \Kb^\perp$,
$\Gamb = \Ub^\perp$, it is obvious that the pairs
$(\Kb,\Ub)$, $(\Gamb,\Omb)$ are $\alfa$-adjoint for each
$\alfa \in (0,2\pi/3)$. The quotients
$G = \Kb/\Ub$ and $\Gamb/\Omb$ are finite abelian groups
and the latter can be canonically identified with the dual
group $\dG$. Any right inverse map $\eta\:G \to \Kb \sbs \Gb$
to the canonical projection $\zeta\:\Kb \to \Kb/\Ub$ is a
$(\Kb,\Ub)$ approximation of \,$\Gb$, and similarly, any right
inverse map $\phi\:\dG \to \Gamb \sbs \dGb$ to the canonical
projection $\xi\:\Gamb \to \Gamb/\Omb$ is a $(\Gamb,\Omb)$
approximation of \,$\dG$. Then one can easily verify by a
straightforward computation that
$$
(\phi\,\gama)(\eta\,a) = \gama(a)
$$
for any $a \in G$, $\gama \in \dG$. Finally, the equalities
$$\align
\dG &= \Bohr_\alfa\bigl(\eta^{-1}[\Ub]\bigr) =
\phi^{-1}\bigl[\Bohr_\eps(\Ub)\bigr]\,,\\
G &= \Bohr_\alfa\bigl(\phi^{-1}[\Omb]\bigr) =
\eta^{-1}\bigl[\Bohr_\eps(\Omb)\bigr]\,
\endalign$$
show that, for any $0 < \eps < \alfa \le \pi/3$, the strong
$(\alfa,\eps)$-adjointness of the approximations $\eta$, $\phi$
with respect to the pairs $(\Kb,\Ub)$, $(\Gamb,\Omb)$ is witnessed
by the very sets $\Vb = \Ub$, $\Ypb = \Omb$.

If \,$\Ub \sbs \Kb$ are $^*$compact $^*$open subgroups of \,$\sGb$,
such that  $\Ub \sbs \IG$, $\FG \sbs \Kb$, then
$\Omb = \Kb^{\ort} = \Kb^\perp$, $\Gamb = \Ub^{\ort} = \Ub^\perp$
are $^*$compact $^*$open subgroups of \,$\sdGb$, such that
$\Omb \sbs \IdG$, $\FdG \sbs \Gamb$. The quotients $G = \Kb/\Ub$,
$\dG = \Gamb/\Omb$ are mutually dual hyperfinite abelian groups.
The HFI approximations $\eta\:G \to \sGb$, $\phi\:\dG \to \sdGb$
can be constructed in essentially the same way as in the standard
situation above. The corresponding IMG group triplets $(G,\Go,\Gf)$,
$\bigl(\dG,\Ho,\Hf\bigr)$, where
$$\xalignat{2}
\Go &= \zeta[\IG]\,, & \Gf &= \zeta[\FG]\,,\\
\Ho &= \xi\bigl[\IdG\bigr]\,, & \Hf &= \xi\bigl[\FdG\bigr]\,,
\endxalignat$$
are, clearly, mutually dual, as well. For any internal subgroup
$\Xb$ of \,$\Gb$, such that $\IG \sbs \Xb \sbs \FG$,
$d = [\Xb:\Ub]^{-1}$ and $\hat d = (\lv G\rv d)^{-1} = [\Kb:\Xb]^{-1}$
can serve as adjoint normalizing multipliers for the triplets
$(G,\Go,\Gf)$ and $\bigl(\dG,\Ho,\Hf\bigr)$, respectively.
The equality $(\phi\,\gama)(\eta\,a) = \gama(a)$ holds even for
all $a \in G$, $\gama \in \dG$, while
$(\phi\,\gama)(\eta\,a) \apr \gama(a)$ for $a \in \Gf$,
$\gama \in \Hf$ would be enough to establish the adjointness of
\,$\eta$, $\phi$.
\endexample

In \cite{Go2} also pairs of adjoint HFI approximations for
$\tau$-adic solenoids $\Sigma_\tau$ and their dual groups of
$\tau$-adic rationals
$$
\bQ^{(\tau)} = \bigcup_{n \in \bN}
   \Bigl(\frac{1}{\tau_n}\Bigr)\,\bZ =
   \Bigl\{\frac{a}{\tau_n}\: a \in \bZ \et n \in \bN\Bigr\}
$$
are described. These, as well as the corresponding finite
approximations, can be constructed combining some ideas from
Examples~2.5.10 and 2.5.12.

In view of the supplements on injectivity after 2.5.2 and in 2.5.8,
as well as of the last three examples it is natural to formulate
the following

\proclaim{Conjecture}
Let \,$\Gb$ be a Hausdorff LCA group, and $\alfa$, $\eps$, $\Kb$,
$\Ub$, $\Gamb$, $\Omb$ be as above. Then there exist a finite
abelian group $G$ and \,{\rm injective\/} mappings
$\eta\:G \to \Gb$, $\phi\:\dG \to \dGb$ forming a strongly
$(\alfa,\eps)$-adjoint pair of approximations of \,$\Gb$, $\dGb$,
respectively, with respect to $(\Kb,\Ub)$, $(\Gamb,\Omb)$.
\endproclaim

Let us close this chapter by a brief discussion of the roles of
standard and nonstandard methods in the proof of the Strongly
Adjoint Finite LCA Group Approximation Theorem~2.5.8. Our starting
point was the proof of existence of arbitrarily good (standard)
finite approximations of any single LCA group $\Gb$ (Theorem~1.5.11)
from which we derived the existence of (nonstandard) hyperfinite
approximations of \,$\Gb$ by means of the {\it transfer principle\/}
(Corollary~1.5.12). Next we proved another nonstandard result,
namely the existence of an adjoint approximation to any HFI
approximation of \,$\Gb$ (Theorem~2.5.1 and Corollary~2.5.2).
Both in their formulation and proof the (inherently nonstandard)
Gordon's Conjecture~1 (Theorem~2.1.4) and the Triplet Duality
Theorem~2.1.5 were crucial. Finally we turned back, deriving the
existence of (standard) arbitrarily good strongly adjoint pairs
of finite approximations of \,$\Gb$ and its dual group $\dGb$
(Theorem~2.5.8) from the nonstandard Corollary~2.5.2, using the
{\it transfer principle\/} in the ``opposite direction''.
\newpage

{}\phantom{x}
{}\bigskip

\head
3. The Fourier Transform in Hyperfinite Dimensional Ambience
\endhead
{}\bigskip\bigskip

\flushpar
Our last Chapter deals with the second of the main topics of this paper,
namely the analysis of the discrete Fourier transform on some subspaces
of the hyperfinite dimensional linear space $\sbC^G$, arising from a
condensing IMG group triplet $(G,\Go,\Gf)$ with hyperfinite abelian
ambient group $G$, and its application to the Fourier transform on
various spaces of functions $\fb\:\Gb \to \bC$ defined on its observable
trace, the Hausdorff LCA group $\Gb = \Gf/\Go$. In particular, we will
formulate and prove a generalization of the third of Gordon's Conjectures
to approximations of the Fourier transforms
$\Lb^1(\Gb) \to \Cbo\bigl(\dGb\bigr)$,
$\Mb(\Gb) \to \Cbu\bigl(\dGb\bigr)$ and
$\Lbp(\Gb) \to \Lbq\bigl(\dGb\bigr)$, for adjoint exponents
$1 < p \le 2 \le q < \infty$, by the discrete Fourier transform
$\sbC^G \to \sbC^{\dG}$.

Throughout the first three sections of Chapter~3, $(G,\Go,\Gf)$ denotes
a condensing IMG group triplet with hyperfinite abelian ambient group~$G$
in a sufficiently saturated nonstandard universe. Its observable trace is
denoted by $\Gb = G^\be = \Gf/\Go$, and it is a Hausdorff LCA group.
Then $d$ denotes a normalizing coefficient for the triplet and all the
norms $\LV\cd\RV_{_p}$, for $1 \le p < \infty$, on the linear space
$\sbC^G$ are defined using $d$; similarly, $\mb = \mbd$ denotes the
Haar measure on $\Gb$, obtained by pushing down the Loeb measure
$\lam_d$ on $G$ and the norms $\LV\cd\RV_{_p}$ on the Lebesgue spaces
$\Lb^p(\Gb)$ are defined via $\mb$. Analogous convention is adopted for
the dual group $\dGb$ and its Haar measure $\mb_{\hat d}$ obtained from
the normalizing multiplier $\hat d = \big(d\lv G\rv\bigr)^{-1}$  for
the dual triplet $\bigl(\dG,\Gf^{\ort},\Go^{\ort}\bigr)$.

\specialhead
3.1.~A characterization of liftings
\endspecialhead
\flushpar
In Section~1.4 we just remarked that for a locally compact
Hausdorff space $\Xb$, represented as the observable trace
$\Xb \iso X^\be = \Xf/E$ of an IMG triplet $(X,E,\Xf)$ with
hyperfinite $X$, not every $S$-integrable function $f \in \sbC^X$
is lifting of some function $\fb \in \Lb^1(\Xb)$, but were not able
to describe these liftings more closely. For Hausdorff LCA groups,
however, we can give an intuitively appealing characterization of
such liftings in terms of certain continuity condition.

Let $\Nn$ be an arbitrary internal norm on the vector space
$\sbC^G$. An internal function $f\:G \to \sbC$ is called
{\it $S$-continuous with respect to the norm \,$\Nn$}, or,
briefly, {\it $S^{\nn}$-continuous} if
$$
\Nn(f_a - f) \apr 0
$$
for each $a \in \Go$. In case of the $p$-norms we speak
about $S^p$-continuous functions. In particular,
$S^{\infty}$-continuity coincides with the usual notion
of $S$-continuity.

The following Lemma is obvious, once we realize that
$$
(f * g)_a - f * g  = (f_a - f) * g\,,
\qquad\text{and}\qquad
\Nn(f * g) \le \Nn(f)\LV g\RV_{_1}
$$
for any functions $f,g \in \sbC^G$, $a \in G$ and an internal
translation invariant norm $\Nn$.

\proclaim{3.1.1.~Lemma}
Let $\Nn$ be an internal translation invariant norm on $\sbC^G$
and $f,g \in \sbC^G$. If \,$f$ is $S^{\nn}$-continuous and
\,$\|g\|_1 < \infty$ then $f * g$ is $S^{\nn}$-continuous, as well.
\endproclaim

In analyzing the structure of $S^{\nn}$-continuous functions we
will make use of a family of internal functions akin to the family
$h_{\ro r}$ defined in the proof of Proposition~2.3.3. For any
valuation $\ro \in \VV$ (see the text preceding Proposition~1.5.1)
and $r > 0$ we put
$$
\th_{\ro r} = \|h_{\ro r}\|_{_1}^{-1} h_{\ro r}\,.
$$
Then each of the functions $\th_{\ro r}$ is $S$-continuous, even,
nonnegative, and satisfies both $\|\th_{\ro r}\|_{_1} = 1$
and $\|\th_{\ro r}\|_{_\infty} < \infty$. Moreover,
$\Go \sbs \supp\th_{\ro r} \sbs B_\ro(r)$, thus, in particular,
$\th_{\ro r} \in \CCc(G,\Go,\Gf)$.

The family of internal functions $\th_{\ro r}$ behaves like an
approximate unit for the operation of convolution on the set of
all $S^{\nn}$-continuous functions in the sense of \cite{HR1},
\cite{HR2}. The precise formulation follows.

\proclaim{3.1.2.~Lemma}
Let $\Nn$ be any internal norm on $\sbC^G$. Then for every
$S^{\nn}$-continuous function $f \in \sbC^G$ the system of functions
$\th_{\ro r} * f$, where $\ro \in \VV$, $0 < r \in \bR$, converges
to the function $f$ with respect to the norm $\Nn$ in the following
sense: for each (standard) $\eps > 0$ there is an internal set $Q$
such that $\Go \sbs Q \sbs G$ and for any $\ro$, $r$ the inclusion
$B_\ro(r) \sbs Q$ implies $\Nn(f - \th_{\ro r} * f) \le \eps$.
Consequently, if \,$\ro \in {}^*\VV$, $0< r \in \sbR$ are such that
$B_\ro(r) \sbs \Go$, then $\Nn(f - \th_{\ro r} * f) \apr 0$.
\endproclaim

Let us remark, that for each internal set $Q \sps \Go$ there are
indeed a $\ro \in \VV$ and an $r > 0$ such that $B_\ro(r) \sbs Q$,
so that the situation described in the Lemma is not merely
hypothetical.

\demo{Proof}
The $S^{\nn}$-continuity of $f$ means, in standard terms, that
for each $\eps > 0$ there is an internal set $Q \sps \Go$ such
that $\Nn(f_a - f) \le \eps$, for $a \in Q$. Now, assume that
$B_\ro(r) \sbs Q$. As $\th_{\ro r}$ is nonnegative,
$\|\th_{\ro r}\|_{_1} = 1$, and
$\supp\th_{\ro r} \sbs B_\ro(r)$, we have
$$
f - \th_{\ro r} * f =
\LV \th_{\ro r}\RV_{_1} f - d \,\sum_{a \in G} \th_{\ro r}(a)\,f_a
= d \,\sum_{a \in B_\ro(r)} \th_{\ro r}(a)(f - f_a)\,,
$$
hence
$$
\Nn(f - \th_{\ro r} * f) \le
d \,\sum_{a \in B_\ro(r)} |\th_{\ro r}(a)|\,\Nn(f - f_a) \le
\|\th_{\ro r}\|_{_1} \max_{a \in B_\ro(r)} \Nn(f - f_a) \le \eps\,,
$$
since $\Nn(f - f_a) \le \eps$ for $a \in B_\ro(r) \sbs Q$.

The infinitesimal supplement is an immediate consequence of the
standard statement just proved.
\enddemo

The last of our lemmas deals with a density condition for certain
$S^{\nn}$-continuous functions. To this end denote
$\CCc^{\nn,1}(G,\Go,\Gf)$ the $\FC$-linear subspace of the internal
space $\sbC^G$, consisting of all $S^{\nn}$-continuous functions
$f \in \sbC^G$ satisfying $\Nn(f) < \infty$,
$\LV f\RV_{_1} < \infty$ and $\supp f \sbs \Gf$. If $\Nn$ is the
$p$-norm $\LV\cd\RV_{_p}$ we write $\CCc^{p,1}(G,\Go,\Gf)$.

\proclaim{3.1.3.~Lemma}
Let \,$\Nn$ be any internal norm on $\sbC^G$. If the subspace
$\CCc(G,\Go,\Gf)$ is contained in the subspace
\,$\CCc^{\nn,1}(G,\Go,\Gf)$ then \,$\CCc(G,\Go,\Gf)$ is dense in
\,$\CCc^{\nn,1}(G,\Go,\Gf)$ with respect to the norm $\Nn$.
\endproclaim

\demo{Proof}
Taking any function $f \in \CCc^{\nn,1}(G,\Go,\Gf)$, we know that
the system of functions $\th_{\ro r} * f$, where $\ro \in \VV$,
$r > 0$, converges to $f$ with respect to $\Nn$ in the sense of
Lemma~3.1.2. It remains to show that
$\th * f \in \CCc(G,\Go,\Gf)$ for each such a function
$\th = \th_{\ro r}$.

Clearly,
$$
\supp(\th * f) \sbs \supp\th + \supp f \sbs \Gf\,,
$$
and, according to the fact that the max-norm $\|\cd\|_{_\infty}$ is
translation invariant,
$$
\|\th * f\|_{_\infty} \le
   \|\th\|_{_\infty} \|f\|_{_1} < \infty\,.
$$
For the same reason, as the function $\th$ is $S$-continuous and
$\|f\|_1 < \infty$, the $S$-continuity of the function $\th * f$
follows from Lemma~3.1.1 applied to the norm $\|\cd\|_{_\infty}$.
\enddemo

As all the internal norms $\LV\cd\RV_{_p}$, $1 \le p \le \infty$,
on $\sbC^G$ are translation invariant and satisfy the inclusion
$\CCc(G,\Go,\Gf) \sbs \CCc^{p,1}(G,\Go,\Gf)$, all the
Lemmas~3.1.1--3 apply particularly to them.

Further, let us denote
$$
\MM(G,\Go,\Gf) = \bigl\{f \in \sbC^G\: \LV f\RV_{_1} < \infty
\et \bigl(\alli Z \sbs G \sms \Gf\bigr)
\bigl(\LV f \cd 1_Z\RV_{_1} \apr 0\bigr)\bigr\}
$$
the $\FC$-linear subspace of \,$\sbC^G$, formerly denoted as
$\MM(G,\Gf,d)$. Notice that displaying $\Go$ and hiding $d$ is fully
justified and unambiguous, as $\MM(G,\Gf,d) = \MM(G,\Gf,d')$ for
any normalizing multipliers $d$, $d'$ of the triplet $(G,\Go,\Gf)$.
The relation of the subspace $\MM(G,\Go,\Gf) = \MM(G,\Gf,d)$ to the
Banach space $\Mb(\Gb)$ of all complex regular Borel measures with
finite variation on $\Gb$ via weak liftings is described in
Proposition~1.4.3. Similarly as in Section~1.4 we put
$$
\MM^p(G,\Go,\Gf) =
   \bigl\{f \in \sbC^G\: \lv f\rv^p \in \MM(G,\Go,\Gf)\bigr\}\,,
$$
for $1 \le p < \infty$.

We are going to characterize the subspaces $\LL^p(G,\Go,\Gf)$
of the internal linear space $\sbC^G$ formed by liftings
$f \in \MM^p(G,\Go,\Gf)$ of functions $\fb \in \Lbp(\Gb)$, for
$1 \le p < \infty$. The following Theorem resembles an early
theorem by Rudin \cite{Rd1}, characterizing measures
$\mmu \in \Mb(\Gb)$ arising from functions $\fb \in \Lb^p(\Gb)$
(i.e., $\dif\mmu = \fb\dif\mb$) as those for which the shift
$\ab \mto \mmu_\ab(\Bb) = \mmu(\Bb - \ab)$ is a continuous
function $\Gb \to \bC$, for any Borel set $\Bb \sbs \Gb$.

\proclaim{3.1.4~Theorem}
Let \,$1 \le p < \infty$ and $f \in \MM^p(G,\Go,\Gf)$ be an
internal function. Then $f \in \LL^p(G,\Go,\Gf)$ if and only if
$f$ is $S^p$-continuous.
\endproclaim

\demo{Proof}
For brevity's sake, let us denote
$$
\CMp(G,\Go,\Gf) =
\{f \in \MM^p(G,\Go,\Gf)\: \text{\rm $f$ is $S^p$-continuous}\}\,.
$$
Then we are to prove that $\LL^p(G,\Go,\Gf) = \CMp(G,\Go,\Gf)$.
Clearly, they both are subspaces of the internal vector space
$\sbC^G$ and contain the subspace $\CCc(G,\Go,\Gf)$. We divide
the proof into three simpler Claims. Putting them together, the
Theorem easily follows.
\smallskip
\flushpar
{\it Claim\/}~1. $\LL^p(G,\Go,\Gf)$ is closed in $\sbC^G$ with
respect to the norm $\LV\cd\RV_{p}$.
\smallskip
This, however, is almost obvious. If $(f_n)_{n\in\bN}$ is a sequence
in $\LL^p(G,\Go,\Gf)$, converging to a function $f \in \sbC^G$, and
each $f_n$ is a lifting of an $\fb_n \in \Lbp(\Gb)$, then the
sequence $(\fb_n)_{n\in\bN}$ satisfies the Bolzano-Cauchy condition,
hence it converges to a function $\fb \in \Lbp(\Gb)$. It is routine
to check that $f$ is a lifting of $\fb$, i.e.,
$f \in \LL^p(G,\Go,\Gf)$.
\smallskip
\flushpar
{\it Claim\/}~2. $\LL^p(G,\Go,\Gf) \sbs \CMp(G,\Go,\Gf)$\,.
\smallskip
It suffices to show that each lifting $f$ of a function
$\fb \in \Lbp(\Gb)$ is $S^p$-continuous. It is known that the
shift $\ab \mto \fb_{\ab}$ is a uniformly continuous mapping
$\Gb \to \Lb^p(\Gb)$ (see \cite{Pd} or \cite{Rd2}). Translating
this condition into the language of infinitesimals, one readily
obtains the $S^p$-continuity of $f$.
\smallskip
\flushpar
{\it Claim\/}~3. $\CCc(G,\Go,\Gf)$ is dense in $\CMp(G,\Go,\Gf)$
with respect to the norm $\LV\cd\RV_{_p}$.
\smallskip
According to Lemma~3.1.3 it is enough to show that the subspace
$\CCc^{p,1}(G,\Go,\Gf)$ is dense in $\CMp(G,\Go,\Gf)$ with respect
to $\LV\cd\RV_{_p}$. Let $f \in \CMp(G,\Go,\Gf)$. Then
\hbox{$\|f \cd 1_{G\sms Z}\|_{_p} \apr 0$} for any internal set
$Z \sbs G \sms \Gf$. Due to {\it saturation}, there exists a
sequence of internal sets $(A_n)_{n\in\bN}$ such that
$\Go \sbs A_n \sbs \Gf$ and
$$
\|f \cd 1_{G \sms A_n}\|_{_p} \to 0
\quad\text{for $n \to \infty$,}
$$
and we can additionally assume that $A_n + A_n \sbs A_{n+1}$.
Then, for each $n$, there is an $S$-continuous function
$g_n \in \sbC^G$ such that $g_n(x) = 1$ for $x \in A_n$,
$g_n(x) = 0$ for $x \in G \sms A_{n+1}$ and
$0 \le g_n(x) \le 1$ for $x \in A_{n+1} \sms A_n$. We put
$f_n = f \cd g_n$. From
$$
\|f - f_n\|_{_p} \le \|f \cd 1_{G \sms A_n}\|_{_p}
$$
it follows that $\|f - f_n\|_{_p} \to 0$.

Let us show that $f_n \in \CCc^{p,1}(G,\Go,\Gf)$
for each $n$. Clearly, $\supp f_n \sbs A_{n+1} \sbs \Gf$
and $\|f_n\|_{_p} \le \|f\|_{_p} < \infty$. According
to H\"older's inequality,
$$
\|f_n\|_{_1} = \|f\,g_n\|_{_1} \le
   \LV f\RV_{_p} \LV g_n\RV_{_q} < \infty\,,
$$
where ${1 \over p} + {1 \over q} = 1$. Taking any $a \in \Go$
we have
$$\align
\LV(f_n)_a - f_n\RV_{_p} &\le
\LV f_a((g_n)_a - g_n)\RV_{_p} + \LV(f_a - f)g_n\RV_{_p} \\
&\le \LV f_a\RV_{_p} \LV(g_n)_a - g\RV_{_\infty}
+ \LV f_a - f\RV_{_p} \LV g_n\RV_{_\infty} \apr 0\,,
\endalign$$
showing that $f_n$ is $S^p$-continuous. Hence
$f_n \in \CCc^{p,1}(G,\Go,\Gf)$.
\enddemo

Theorem~3.1.4 and Proposition~1.4.4 yield the following

\proclaim{3.1.5.~Corollary}
Let $1 \le p < \infty$. Then for any measurable function
$\fb\:\Gb \to \bC$ the following conditions are equivalent:
{\parindent 21pt
\item{\sl (i)}\,
$\fb \in \Lbp(\Gb)$;
\item{\sl (ii)}\,
$\fb$ has an S$^p$-integrable lifting;
\item{\sl (iii)}\,
$\fb$ has an $S^p$-continuous lifting $f \in \MM^p(G,\Go,\Gf)$.
}
\endproclaim

\remark{Remark\/ \rm 1}
Let $1 \le p < \infty$. We adopt a similar and equally justified
convention $\SS^p(G,\Go,\Gf) = \SS^p(G,\Gf,d)$ like that for
$\MM^p(G,\Go,\Gf)$. Then, according to Theorem~3.1.4, we have
$$\align
\LL^p(G,\Go,\Gf) &= \bigl\{f \in \MM^p(G,\Go,\Gf)\: (\all a \in G)
\bigl(a \apr 0 \imp \LV f_a - f\RV_{_p} \apr 0\bigr)\bigr\}\,,\\
\SS^p(G,\Go,\Gf) &=
\bigl\{f \in \MM^p(G,\Go,\Gf)\: \bigl(\alli A \sbs G\bigr)
\bigl(d\lv A\rv \apr 0 \imp
\LV f \cd 1_A\RV_{_p} \apr 0\bigr)\bigr\}\,,
\endalign$$
so that the characterization of $\LL^p(G,\Go,\Gf)$ differs from
the definition of $\SS^p(G,\Go,\Gf)$ just in replacing the condition
of absolute continuity by that of $S^p$-continuity. As it follows
from 3.1.4 and 3.1.5, $\LL^p(G,\Go,\Gf) \sbs \SS^p(G,\Go,\Gf)$,
i.e., for a function $f \in \MM^p(G,\Go,\Gf)$, $S^p$-continuity
implies absolute continuity (but not vice versa). However, one
would like to have a more direct proof of this inclusion.
\endremark
\smallskip

\remark{Remark\/ \rm 2}
Given any condensing IMG triplet $(X,E,\Xf)$ with hyperfinite
ambient set $X$ and a nonnegative internal function $d\:X \to \sbR$,
the observable trace $\Xb = \Xf/E$ is a Hausdorff locally compact
space, so that it still makes sense to ask which internal functions
$f\:X \to \sbC$ are liftings of functions $\fb \in \Lbp(\Xb,\mb)$,
where $\mb = \mbd$ is the Lebesgue measure on $\Xb$ obtained by
pushing down the Loeb measure $\lam_d$. However, as long as no
group structure on $X$ is involved, $\LL^p(X,E,\Xf)$ cannot be
characterized in terms of \,$S^p$-continuity. It would be
nice to have some reasonable intrinsic characterization of
\,$\LL^p(X,E,\Xf)$ within such a more general setting, at least
for constant $d(x) = d$ such that $d\lv A\rv \not\apr 0$ for some
and $d\lv A\rv < \infty$ for each internal set $A \sbs \Xf$.
\endremark
\smallskip

\remark{Remark\/ \rm 3}
The characterizing conditions of $\LL^p(G,\Go,\Gf)$ make
sense also for $p = \infty$. More precisely, the conjunction
of \,$\LV f \RV_{_\infty} < \infty$,
\,$\LV f(x) \cd 1_Z\RV_{_\infty} \apr 0$ for each internal
set $Z \sbs G \sms \Gf$ and $S^\infty$-continuity defines the
subspace $\CCo(G,\Go,\Gf)$ of $S$-continuous internal functions
$f\:G \to \sbC$ which are finite on the whole $G$ and
infinitesimal outside of \,$\Gf$. Thus we could formally write
$\LL^\infty(G,\Go,\Gf) = \CCo(G,\Go,\Gf)$. This, however, would
interfere with the adopted standard notation, as such an
$\LL^\infty(G,\Go,\Gf)$ would be formed just by the liftings of
functions in $\Cbo(\Gb)$ which is a proper closed subspace of
the Banach space $\Lb^\infty(\Gb)$ (cf\. Proposition~1.3.1).
\endremark

\specialhead
3.2.~The Smoothness-and-Decay Principle
\endspecialhead
\flushpar
The more smooth is a function $f\:\bR^n \to \bC$, the more rapidly
its Fourier transform $\wh f\:\bR^n \to \bC$ decays, and vice
versa, the more rapidly a function  $f\:\bR^n \to \bC$ decays,
the smoother is its Fourier transform $\wh f\:\bR^n \to \bC$. This
vague informal statement is known as the {\it Smoothness-and-Decay
Principle\/} and\,---\,jointly with the {\it Uncertainty
Principle\/} to which it is closely related\,---\,belongs to
fundamental heuristic principles of Fourier or time-frequency
analysis. It can take the form of various precise mathematical
statements some of which generalize from $\bR$ or $\bR^n$ to
arbitrary LCA groups (see, e.g., \cite{Gc}, \cite{Ta3} for
discussion).

The view through the lense of an IMG group triplet $(G,\Go,\Gf)$
with hyperfinite abelian ambient group $G$ and its dual triplet
$\bigl(\dG,\Gf^{\ort},\Go^{\ort}\bigr)$ offers an intuitively
appealing explanation of this principle for internal functions
$f\: G \to \sbC$, based on the Fourier inversion formula
$$
f(x) = \hat d \,\sum_{\gama \in \dG} \wh f(\gama)\,\gama(x)
$$
in which both $S$-continuous characters $\gama \in \Go^{\ort}$
as well as non-$S$-continuous characters
$\gama \in \dG \sms \Go^{\ort}$ occur. If $f$ is smooth or continuous
(in some intuitive meaning of these words), then the contribution of
the non-$S$-continuous characters to the above expansion of \,$f$
must be negligible in some sense. This condition causes a kind of
quick decay of \,$\wh f$. The other way round, viewing the elements
$x \in G$ as characters of the dual group $\dG$, the Fourier
transform of \,$f$ can be expressed as their linear combination:
$$
\wh f(\gama) = d \,\sum_{x \in G} f(x)\,\ovl\gama(x)
= d \,\sum_{x \in G} f(-x)\,x(\gama)\,.
$$
If \,$f$ decays quickly, i.e., if the values of \,$f$ on the
infinite elements $x \in G \sms \Gf$ are somehow negligible, then
the values of its Fourier transform are essentially determined by
the values of \,$f$ on the finite elements $x \in \Gf$, which
happen to coincide with the $S$-continuous characters of \,$\dG$
by the Triplet Duality Theorem~2.1.5. If, additionally, neither
the coefficients $f(x)$, for $x \in \Gf$, are too big, then we can
reasonably expect $\wh f$ to be smooth or continuous in some sense.

The next theorem is a fairly general precise statement of this kind
of the {\it Smoothness-and-Decay Principle}. Both its formulation as
well as its proof borrow some ideas from a paper by Pego \cite{Pg}.

A~pair of internal norms $\Nn$ on $\sbC^G$ and $\Mm$ on $\sbC^{\dG}$
is called {\it Fourier compatible\/} if the Fourier transform
$\FF\: \sbC^G \to \sbC^{\dG}$ is a bounded linear operator with
respect to the norms $\Nn$, $\Mm$, i.e.,
$$
\Nn(f) < \infty \imp \Mm\bigl(\wh f\,\bigr) < \infty
$$
for each $f \in \sbC^G$. This is equivalent to the $S$-continuity
of \,$\FF$, i.e.,
$$
\Nn(f) \apr 0 \imp \Mm\bigl(\wh f\,\bigr) \apr 0\,.
$$
Also recall from Section~2.2 that a norm $\Nn$ on $\sbC^G$ is called
{\it absolute\/} if $\Nn(f) \le \Nn(g)$ for any $f,g \in \sbC^G$
such that $|f(x)| \le |g(x)|$ for all $x \in G$.

\proclaim{~3.2.1.~Theorem}{\bf [Smoothness-and-Decay Principle]}
Let \,$\Nn$, $\Mm$ be Fourier compatible internal norms on the
linear spaces $\sbC^G$, $\sbC^{\dG}$, respectively. Then for every
function $f \in \sbC^G$ the following implications hold:
{\parindent 21pt
\item{\sl (a)}
If \,$\Mm$ is absolute and \,$f$ is $S^{\nn}$-continuous,
then \,$\Mm\bigl(\wh f \cd 1_\Gama\bigr) \apr 0$ for every internal
set \,$\Gama \sbs \dG \sms \Go^{\ort}$.
\item{\sl (b)}
If \,$\Nn$ is absolute, $\Nn(f) < \infty$, and
\,$\Nn(f \cd 1_X) \apr 0$ for every internal set
\,$X \sbs G \sms \Gf$, then $\wh f$ is $S^{\mm}$-continuous.
\item{}
\vskip-12pt}
\endproclaim

\demo{Proof}
(a) In this part of proof we will once more make use of the
families of internal functions $h_{\ro r}$ and $\th_{\ro r}$
(cf\. Lemma~3.1.2 and its proof).

Assume that \,$\Nn(f_a - f) \apr 0$ for any $a \in \Go$.
We will show that \,$\Mm\bigl(\wh f \cd 1_\Gama\bigr) \apr 0$
for every internal set
$\Gama \sbs \dG \sms \Go^{\ort}$. Let us fix any (standard)
$t \in (0,1)$. By Corollary~2.3.2,
$\Gama \cap \Spec_t(h_{\ro r}) = \emptyset$
for every $\ro \in \VV$ and standard $r > 0$. As $\VV$ is upward
directed, by {\it saturation\/} there are a $\tau \in {}^*\VV$,
satisfying $\ro \le \tau$ for all $\ro \in \VV$, and a positive
$s \apr 0$, such that
\,$\Gama \cap \Spec_t(h_{\tau s}) = \emptyset$ still holds.
Let us denote $\Dela = \dG \sms \Spec_t(h_{\tau s})$ and recall
that \,$\th_{\tau s} = \LV h_{\tau s}\RV_{_1}^{-1} h_{\tau s}$.
As $h_{\tau s}$ is even and nonnegative, so is $\th_{\tau s}$,
hence
$$
\wh\th_{\tau s}(1_G) = d \,\sum_{a \in G} \th_{\tau s}(a)
= \LV\th_{\tau s}\RV_{_1} = 1\,,
$$
and for $\gama \in \Dela$ we have
$$
\bigl|\wh\th_{\tau s}(\gama)\bigr| =
\LV h_{\tau s}\RV_{_1}^{-1}\bigl|\wh h_{\tau s}(\gama)\bigr| < t\,,
$$
hence,
$$
1 - t < 1 - \bigl|\wh\th_{\tau s}(\gama)\bigr|
\le \bigl|1 - \wh\th_{\tau s}(\gama)\bigr|\,.
$$
Further on, $\Gama \sbs \Dela$, and as the norm $\Mm$ is absolute,
$$
(1 - t)\,\Mm\bigl(\wh f \cd 1_\Gama\bigr)
   \le (1 - t)\,\Mm\bigl(\wh f \cd 1_\Dela\bigr)
   \le \Mm\bigl(\bigl(1_{\dG} - \wh\th_{\tau s}\bigr)\,\wh f\,\bigr)
     = \Mm\bigl((f - \th_{\tau s} * f)\sphat\,\bigr)\,.
$$
Due to our choice of \,$\tau$ and $s$ we have $B_\tau(s) \sbs \Go$,
consequently, $\Nn(f - \th_{\tau s} * f) \apr 0$ by the virtue of
Lemma~3.1.2. As the norms $\Nn$, $\Mm$ are Fourier compatible, this
implies that \,$\Mm\bigl((f - \th_{\tau s} * f)\sphat\,\bigr) \apr 0$,
as well. Since $t \not\apr 1$, we can conclude that
\,$\Mm(\wh f \cd 1_\Gama) \apr 0$.

(b) Assume that $\Nn(f)$ is finite and $\Nn(f \cd 1_X) \apr 0$ for
each internal set $X \sbs G \sms \Gf$. We are to show
that $\Mm\bigl(\wh f_{\,\gama} - \wh f\,\bigr) \apr 0$ for any
$\gama \in \Gf^{\ort}$. First notice that
$$
\wh f_{\,\gama} - \wh f = \bigl((\gama - 1_G)\,f\bigr)\sphat\,.
$$
As $\gama \in \Gf^{\ort}$, \,$\gama(x) \apr 1$ for each $x \in \Gf$.
Due to {\it saturation}, there is an internal set $Y$ such that
$\Gf \sbs Y \sbs G$ and $\gama(y) \apr 1$ for each $y \in Y$;
then $X = G \sms Y \sbs G \sms \Gf$. Let us denote
$$
\eps = \LV(\gama - 1_G)\cd 1_Y\RV_{_\infty}
     = \max_{y \in Y} \lv\gama(y)- 1\rv\,;
$$
obviously, $\eps \apr 0$. As $\Nn$ is absolute,
$$\align
\Nn\bigl((\gama - 1_G)\,f) &\le
\Nn\bigl((\gama - 1_G)\,f \cd 1_Y\bigr)
   + \Nn\bigl((\gama - 1_G)\,f \cd 1_X\bigr) \\
&\le \eps\,\Nn(f) + 2\,\Nn(f \cd 1_X) \apr 0\,.
\endalign$$
Therefore, $\Mm\bigl(\wh f_{\,\gama} - \wh f\,\bigr) \apr 0$,
as well.
\enddemo

The last Theorem applies to any pair of norms $\LV\cd\RV_{_p}$ on
$\sbC^G$ and $\LV\cd\RV_{_q}$ on $\sbC^{\dG}$ for $1 \le p \le 2$
and $q = p/(p-1)$, including $p=1$, $q=\infty$, in which case we
have:

\proclaim{3.2.2.~Corollary}
For every function $f \in \sbC^G$ the following conditions hold:
{\parindent 21pt
\item{\sl (a)}
If \,$f$ is $S^1$-continuous, then
\,$\wh f(\gama) \apr 0$ for all \,$\gama \in \dG \sms \Go^{\ort}$.
\item{\sl (b)}
If \,$\LV f\RV_{_1} < \infty$ and \,$\LV f \cd 1_X\RV_{_1} \apr 0$ for
every internal set \,$X \sbs G \sms \Gf$, then $\wh f$ is S-continuous,
i.e., $\wh f(\gama) \apr \wh f(\chi)$ for all \,$\gama,\chi \in \dG$
such that \,$\gama(x) \apr \chi(x)$ for each $x \in \Gf$.
\item{\sl (c)}
\,$\FF\bigl[\LL^1(G,\Go,\Gf)\bigr] \sbs
\CCo\bigl(\dG,\Gf^{\ort},\Go^{\ort}\bigr)$.
}
\endproclaim

Notice that (c) is a hyperfinite dimensional version of the
Riemann-Lebesgue lemma, and (b) can be written as a similar
inclusion
$$\FF\bigl[\MM(G,\Go,\Gf)\bigr] \sbs
   \CCbu\bigl(\dG,\Gf^{\ort}\bigr)\,.
$$

\proclaim{3.2.3.~Corollary}
Let \,$1 < p \le 2$ and \,$q = p/(p-1)$ be its dual exponent. Then for
every function $f \in \sbC^G$ the following conditions hold:
{\parindent 21pt
\item{\sl (a)}
If \,$f$ is $S^p$-continuous, then
$\bigl\|\wh f \cd 1_\Gama\bigr\|_{_q} \apr 0$ for every internal
set \,$\Gama \sbs \dG \sms \Go^{\ort}$.
\item{\sl (b)}
If \,$\LV f\RV_{_p} < \infty$ and \,$\LV f \cd 1_X\RV_{_p} \apr 0$
for every internal set \,$X \sbs G \sms \Gf$, then $\wh f$ is
$S^q$-continuous, i.e.,
\,$\bigl\|\wh f_{\,\gama} - \wh f\,\bigl\|_{_q} \apr 0$ for all
\,$\gama \in \Gf^{\ort}$.
\item{\sl (c)}
\,$\FF\bigl[\LL^p(G,\Go,\Gf)\bigr] \sbs
\LL^q\bigl(\dG,\Gf^{\ort},\Go^{\ort}\bigr)$.
\item{}
\vskip-8pt}
\endproclaim

In the Hilbert space case $p=q=2$  the last Corollary can be
slightly strengthened. Applying 3.2.3 both to the Fourier
transform $\FF\:\sbC^G \to \sbC^{\dG}$ and its inverse
$\FF^{-1}\:\sbC^{\dG} \to \sbC^G$, for functions satisfying
$\LV f\RV_{_2} < \infty$, we get equivalences in (a), (b) and
equality in (c).

\proclaim{3.2.4.~Corollary}
For every function $f \in \sbC^G$ such that $\LV f\RV_{_2} < \infty$
the following conditions hold:
{\parindent 21pt
\item{\sl (a)}
\,$f$ is $S^2$-continuous if and only if
\,$\bigl\|\wh f \cd 1_\Gama\bigr\|_{_2} \apr 0$ for every internal
set \,$\Gama \sbs \dG \sms \Go^{\ort}$.
\item{\sl (b)}
\,$\LV f \cd 1_X\RV_{_2} \apr 0$ for every internal set
\,$X \sbs G \sms \Gf$ if and only if
\,$\wh f$ is $S^2$-continuous.
\item{\sl (c)}
\,$\FF\bigl[\LL^2(G,\Go,\Gf)\bigr] =
\LL^2\bigl(\dG,\Gf^{\ort},\Go^{\ort}\bigr)$.
}
\endproclaim

The following result follows directly from Corollary~3.2.4.

\proclaim{3.2.5.~Corollary}
For every function $f \in \sbC^G$ such that
\,$\LV f\RV_{_2} < \infty$ the following conditions are equivalent:
{\parindent 21pt
\item{\sl (i)}
\,$f \in \LL^2(G,\Go,\Gf)$;
\item{\sl (ii)}
\,$\wh f \in \LL^2\bigl(\dG,\Gf^{\ort},\Go^{\ort}\bigr)$;
\item{\sl (iii)}
both \,$f$ and \,$\wh f$ are $S^2$-continuous;
\item{\sl (iv)}
$\LV f \cd 1_X\RV_{_2} \apr 0$ for every internal set
$X \sbs G \sms \Gf$ and
\,$\bigl\|\wh f \cd 1_\Gama\bigr\|_{_2} \apr 0$
for every internal set \,$\Gama \sbs \dG \sms \Go^{\ort}$.
\item{}
\vskip-8pt}
\endproclaim

Corollary~3.2.5 generalizes a result by Albeverio, Gordon and
Khrennikov \cite{AGK}, where the equivalence of conditions (i),
(ii) and (iv) in case there is an internal subgroup $K$ of $G$
such that $\Go \sbs K \sbs \Gf$ was proved. This assumption is
equivalent to the existence of a compact open subgroup of $G^\be$.
It is also mentioned there without proof that the group of reals
$\bR$, as well, can be represented as $\bR \iso G^\be=\Gf/\Go$ for
some triplet $(G,\Go,\Gf)$ satisfying that way abridged version of
Corollary~3.2.5.

\specialhead
3.3.~Hyperfinite dimensional approximation of the Fourier transform: \\
\phantom{3.3.}~Generalized Gordon's Conjecture 3
\endspecialhead
\flushpar
Various versions of the {\it Smoothness-and-Decay Principle\/}
proved in the previous section make it possible to approximate
the classical Fourier transform on various functional spaces
related to the LCA group $\Gb = \Gf/\Go$ by the discrete
Fourier transform on the hyperfinite dimensional linear space
$\sbC^G$. For the sake of clear distiction, we denote
$\Fb(\fb) = \wh\fb$ the classical Fourier transform of a function
$\fb\:\Gb \to \bC$ and $\FF(f) = \wh f$ the discrete Fourier
transform of an internal function $f\:G \to \sbC$.

The discrete hyperfinite dimensional Fourier transform
$\FF\:\sbC^G \to \sbC^{\dG}$ approximates the classical
Fourier transform $\Fb\:\Lb^1(\Gb) \to \Cbo\bigl(\dGb)$
in the following sense:

\proclaim{3.3.1.~HFD Fourier Transform Approximation Theorem}
Let the internal function $f \in \LL^1(G,\Go,\Gf)$ be a lifting
of a function $\fb \in \Lb^1\bigl(\Gb\bigr)$. Then the internal
function
$\FF(f) = \wh f \in \CCo\bigl(\dG,\Gf^{\ort},\Go^{\ort}\bigr)$
is a lifting of the function
\,$\Fb(\fb) = \wh\fb \in \Cbo\bigl(\dGb\bigr)$.
\endproclaim

\demo{Proof}
Let $f \in \LL^1(G,\Go,\Gf)$ be a lifting of
$\fb \in \Lb^1\bigl(\Gb\bigr)$. Then
$\wh f \in \CCo\bigl(\dG,\Gf^{\ort},\Go^{\ort}\bigr)$
by Corollary~3.2.2(c). Thus it suffices to prove that
$$
\wh\fb\bigl(\gama^\be\bigr) =
{\vphantom{f}}^{\co\!}\wh f(\gama)
$$
for each $\gama \in \Go^{\ort}$. However, as $\gama$ is bounded and
$S$-continuous, i.e., $\gama \in \CCbu(G,\Go)$, it is routine to
check that the internal function $f\,\ovl\gama \in \LL^1(G,\Go,\Gf)$
is a lifting of the function
$\fb\,\ovl\gama^{\,\be} \in \Lb^1\bigl(\dGb\bigr)$. Then
$$
\wh\fb\bigl(\gama^\be\bigr) =
\int \fb\,\ovl\gama^{\,\be}\,\dif\mb
= \,{\vphantom{\biggl(}}^{\co\!} \biggl(d\,\sum_{x \in G}
   f(x)\,\ovl\gama(x)\biggr) =
   {\vphantom{\wh f}}^{\co}\wh f(\gama)
$$
by \cite{ZZ, Proposition~3.5}, see also the text preceding
Proposition~1.4.3.
\enddemo

For $1 < p \le 2$ and $1/p + 1/q = 1$, the Fourier transform
$\Fb\:\Lbp(\Gb) \to \Lbq\bigl(\dGb\bigr)$ is defined as the
continuous extension (with respect to the norms $\LV\cd\RV_p$ on
$\Lbp(\Gb)$ and $\LV\cd\RV_q$ on $\Lbq\bigl(\dGb\bigr)$) of the
restriction of the Fourier transform
$\Fb\:\Lb^1(\Gb) \to \Cbo\bigl(\dGb)$ to the dense subspace
$\Lbp(\Gb) \cap \Lb^1(\Gb)$ of $\Lbp(\Gb)$ (see \cite{HR2} or
\cite{Rd2}). For functions in this subspace everything works like
in the proof above. Thus, by a continuity argument, Theorem~3.3.1
together with Corollary~3.2.3(c) give rise to HFD approximations
of the classical Fourier transforms
$\Fb\:\Lbp(\Gb) \to \Lbq\bigl(\dGb)$ in a similar way.
The case $p=q=2$ of the Fourier-Plancherel transform
$\Fb\:\Lb^2(\Gb) \to \Lb^2\bigl(\dGb\bigr)$ settles
Gordon's Conjecture~3.

\proclaim{3.3.2.~Theorem}{\bf [Generalized Gordon's Conjecture~3]}
Let \,$1 < p \le 2$ and \,$2 \le q < \infty$ be its dual exponent.
Let the internal function $f \in \LL^p(G,\Go,\Gf)$ be a lifting of
a function $\fb \in \Lbp\bigl(\Gb\bigr)$. Then the internal function
$\FF(f) = \wh f \in \LL^q\bigl(\dG,\Gf^{\ort},\Go^{\ort}\bigr)$
is a lifting of the function
\,$\Fb(\fb) = \wh\fb \in \Lbq\bigl(\dGb\bigr)$.
\endproclaim

The HFD Fourier Transform Approximation Theorem~3.3.1 extends
to the Fourier-Stieltjes transform
$\Fb\: \Mb(\Gb) \to \Cbu\bigl(\dGb)$, as well.

\proclaim{3.3.3.~HFD Fourier-Stieltjes Transform Approximation Theorem}
Let the internal function $g \in \MM(G,\Go,\Gf)$ be a weak lifting
of a complex regular Borel measure $\mmu \in \Mb(\Gb)$. Then the
internal function
$\FF(g) = \wh g \in \CCbu\bigl(\dG,\Gf^{\ort}\bigr)$
is a lifting of the function
$\Fb(\mmu) = \wh\mmu \in \Cbu\bigl(\dGb\bigr)$.
\endproclaim

\demo{Proof}
Let $g \in \MM(G,\Go,\Gf)$ be a weak lifting of
\,$\mmu \in \Mb\bigl(\Gb\bigr)$. Then
$\wh g \in \CCbu\bigl(\dG,\Gf^{\ort}\bigr)$ by Corollary~3.2.2(b).
Thus it suffices to prove that
$$
\wh\gb\bigl(\gama^\be\bigr) =
{\vphantom{\wh g}}^{\co\,}\wh g(\gama)
$$
for each $\gama \in \Go^{\ort}$. For the same reason as in the
proof of Theorem~3.3.1 we have
$$
\int f^\be \dif\mmu =
{\vphantom{\biggl(}}^\co\!\biggl(d\,
   \sum_{x \in G} f(x)\,g(x)\biggr)
$$
for every internal function $f \in \CCb(G,\Go,\Gf)$. For
$f = \ovl\gama \in \Go^{\ort} \sbs \CCbu(G,\Go)$ this gives
$$
\wh\mmu\bigl(\gama^\be\bigr) = \int \ovl\gama^{\,\be}\,\dif\mmu
= {\vphantom{\biggl(}}^\co\!\biggl(d\,
   \sum_{x \in G} \ovl\gama(x)\,g(x)\biggr)
= {\vphantom{\wh g}}^{\co\,}{\wh g}(\gama)\,.
$$
\enddemo

In particular, if \,$\gb \in \Lb^1(\Gb)$, $\dif\mmu = \gb \dif\mb$
and \,$g \in \LL^1(G,\Go,\Gf)$ is a lifting of \,$\gb$, then
$$
\wh\gb\bigl(\gama^\be\big) = \wh\mmu\bigl(\gama^\be\big)
= \int \ovl\gama^{\,\be}\,\gb \dif\mb
= {\vphantom{\biggl(}}^\co\!\biggl(d\,
\sum_{x \in G} \ovl\gama(x)\,g(x)\biggr)
= {\vphantom{\wh g}}^{\co\,}{\wh g}(\gama)\,,
$$
for each $\gama \in \Go^{\ort}$, reproving Theorem~3.3.1. This
account indicates that it is Theorem~3.3.3 which is crucial for
hyperfinite dimensional approximations of the Fourier transform
on LCA groups. Therefore we address the issue raised in the Remark
closing the introductory part of Section~2.5 primarily for the
Fourier-Stieltjes transform.

Assume, for the rest of this section, that $(G,\Go,\Gf)$ is an IMG
group triplet with hyperfinite abelian ambient group $G$, arising
from an HFI approximation $\eta\:G \to \sGb$ of the Hausdorff LCA
group $\Gb$. Let us denote $\FF_\eta\:\sbC^G \to \sbC^{\sdGb}$
the internal linear operator given by
$$
\FF_\eta(f)(\cchi) = \br f, \cchi \co \eta\kt_G =
d \,\sum_{x \in G} f(x)\,\ovl\cchi(\eta\,x)\,,
$$
for $f \in \sbC^G$, $\cchi \in \sdGb$. The {\it modified
discrete Fourier transform\/} $\FF_\eta$, defined by means of
the internal inner product on $\sbC^G$, can be employed for
the approximation of the classical Fourier transform on $\Gb$,
without the need to mention the adjoint HFI approximation
$\phi\:\dG \to \sdGb$ of the dual group $\dG$.

\proclaim{3.3.4.~Theorem}
Let \,$\Fb\:\Mb(\Gb) \to \Cbu\bigl(\dGb\bigr)$ be the
Fourier-Stieltjes transform on $\Gb$, $\mmu \in \Mb(\Gb)$
and \,$g \in \MM(G,\Go,\Gf)$ be a weak lifting of \,$\mmu$.
Then, for each \,$\ggam \in \dGb$,
$$
\Fb(\mmu)(\ggam) = \wh\mmu(\ggam) \apr \FF_\eta(g)({}^*\ggam)\,.
$$
\endproclaim

\demo{Proof}
As ${}^*\ggam \co \eta$ is almost homomorphic and $S$-continuous
on $\Gf$, by Theorem~2.1.4 there is a $\gama \in \Go^{\ort}$ such
that ${}^*\ggam(\eta\,x) \apr \gama(x)$ for each $x \in \Gf$.
According to Theorem~3.3.3,
$$
\wh\mmu(\ggam) \apr \wh g(\gama) = \br g,\gama\kt\,.
$$
By the virtue of {\it saturation}, there is an internal set $X$
such that $\Gf \sbs X \sbs G$ and
${}^*\ggam(\eta\,x) \apr \gama(x)$ holds for all $x \in X$.
Denoting $Y = G \sms X$ we have
$$\align
\bigl|\br g, {}^*\ggam \co \eta\kt - \br g, \gama\kt\bigr|
&\le  \bigl|\br g, ({}^*\ggam \co \eta - \gama) \cd 1_X\kt\bigr|
   + \bigl|\br g  \cd 1_Y, {}^*\ggam \co \eta - \gama\kt\bigr| \\
&\le \LV g \RV_{_1}
   \LV ({}^*\ggam \co \eta - \gama)\cd 1_X\RV_{_\infty}
+ \LV g \cd 1_Y\RV_{_1} \LV {}^*\ggam \co \eta - \gama\RV_{_\infty}
\apr 0\,,
\endalign$$
as $\LV g \RV_{_1}$ and
$\LV {}^*\ggam \co \eta - \gama\RV_{_\infty}$ are finite, and
$\LV ({}^*\ggam \co \eta - \gama)\cd 1_X\RV_{_\infty}$ and
$\LV g \cd 1_Y\RV_{_1}$ are infinitesimal. The needed conclusion
$\FF_\eta(g)({}^*\ggam) =
\br g, {}^*\ggam \co \eta\kt \apr\wh\mmu(\ggam)$  is obvious, now.
\enddemo

Theorem~3.3.4, jointly with Theorems~3.3.1 and 3.3.2, respectively,
yield the following two corollaries.

\proclaim{3.3.5.~Corollary}
Let \,$\Fb\:\Lb^1(\Gb) \to \Cbo\bigl(\dGb\bigr)$ be the Fourier
transform on $\Gb$, \hbox{$\fb \in \Lb^1(\Gb)$} and
\,$f \in \LL^1(G,\Go,\Gf)$ be a lifting of \,$\fb$. Then, for each
\,$\ggam \in \dGb$,
$$
\Fb(\fb)(\ggam) = \wh\fb(\ggam) \apr \FF_\eta(f)({}^*\ggam)\,.
$$
\endproclaim

\proclaim{3.3.6.~Corollary}
Let \,$1 < p \le 2 \le q < \infty$ be dual exponents and
\,$\Fb\:\Lbp(\Gb) \to \Lbq\bigl(\dGb\bigr)$ be the Fourier
transform on $\Gb$. Let further $\fb \in \Lbp(\Gb)$ and
\,$f \in \LL^p(G,\Go,\Gf)$ be a lifting of \,$\fb$. Then
$$
\Fb(\fb)(\ggam) \apr \FF_\eta(f)({}^*\ggam)
$$
for almost all \,$\ggam \in \dGb$ with respect to the Haar
measure on $\dGb$.
\endproclaim

\specialhead
3.4.~Standard consequences:
Simultaneous approximation of a function \\
\phantom{3.4.}~and its Fourier transform
\endspecialhead
\flushpar
In this section we are going to apply the previous nonstandard
results to approximations of functions $\fb\:\Gb \to \bC$ and
their Fourier transforms $\wh\fb\:\dGb \to \bC$ by functions
$f\:G \to \bC$ defined on some finite abelian group $G$ and their
discrete Fourier transforms $\wh f\:\dG \to \bC$. To this end we
need to introduce some approximate standard counterparts of various
types of liftings.

Assume that $\Ub \sbs \Gb$ is a neighborhood of \,$0 \in \Gb$,
$\Kb \sbs \Gb$ is a compact set of positive Haar measure and
$\delta > 0$. Let further $G$ be a finite abelian group and
$\eta\:G \to \Gb$ be any mapping such that $\eta(0) \in \Ub$.
Then
\smallskip
{\parindent 21pt
\item{(a)}
given a complex Borel measure $\mmu$ on $\Gb$, a function $f\:G \to \Gb$
is said to be a {\it weak $(\Ub,\delta)$ lifting\/} of \,$\mmu$ on $\Kb$
with respect to $\eta$ if, for each $\xb \in \Kb$, with possible exception
of a subset $\Ab \sbs \Kb$ of measure $\mb(\Ab) \le \delta$,
$$
\lv\,\frac{\mmu(\xb + \Ub)}{\mb(\Ub)} -
   \frac{1}{\lv\eta^{-1}[\Ub]\rv}\,
      \sum_{a \in \eta^{-1}[\xb + \Ub]} f(a)\,\rv
         \le \delta \,;
$$
\item{(b)}
given a measurable function $\fb\:\Gb \to \bC$, a function
$f\:G \to \bC$ is said to be a {\it $(\Ub,\delta)$ lifting\/}
of \,$\fb$ on $\Kb$ with respect to $\eta$ if, for each
$\xb \in \Kb$, with possible exception of a subset
$\Ab \sbs \Kb$ of measure $\mb(\Ab) \le \delta$,
$$
\lv\,\frac{1}{\mb(\Ub)}\,\int_{\xb+\Ub} \fb \dif\mb -
   \frac{1}{\lv\eta^{-1}[\Ub]\rv}\,
      \sum_{a \in \eta^{-1}[\xb + \Ub]} f(a)\,\rv
         \le \delta \,;
$$
\item{(c)}
given a function $\fb\:\Gb \to \bC$, a function $f\:G \to \bC$ is
said to be a {\it $(\Ub,\delta)$ approxim\-ation\/} of \,$\fb$ on
$\Kb$ with respect to $\eta$ if, for any $\xb \in \Kb$,
$a \in \eta^{-1}[\xb + \Ub]$,
$$
\lv\fb(\xb) - f(a)\rv \le \delta\,.
$$
\item{}
\vskip-12pt
}

Intuitively, if \,$\eta$ is a $(\Kb,\Ub)$ approximation of \,$\Gb$,
then $f$ is a $(\Ub,\delta)$ lifting of \,$\fb$ on $\Kb$ with
respect to $\eta$ if the mean value
$\bigl(\int_{\xb+\Ub} \fb \dif\mb\bigr)\big/\mb(\Ub)$ of \,$\fb$
on the neighborhood $\xb + \Ub$ of ``$\delta$-almost each'' point
$\xb \in \Kb$ can be approximated by the ``almost average''
$\bigl(\sum_{a \in \eta^{-1}[\xb + \Ub]} f(a)\bigr)\big/
\lv\eta^{-1}[\Ub]\rv$
with an error at most~$\delta$. Obviously, $f$ is a is a
$(\Ub,\delta)$ lifting of a function $\fb\in \Lb^1(\Gb)$ on $\Kb$
with respect to $\eta$ if and only if it is a weak
$(\Ub,\delta)$ lifting of the measure $\mmu\in \Mb(\Gb)$,
such that $\dif\mmu = \fb\dif\mb$, on $\Kb$ with respect to $\eta$.

If \,$\eta$ is a $(\Kb,\Ub)$ approximation of \,$\Gb$ and $\fb$
satisfies \,$\lv\fb(\xb) - \fb(\yb)\rv \le \delta$ for any
$\xb \in \Kb$, $\yb \in \xb+\Ub$, then, obviously, the composed
function $f =\fb \co \eta$ is a $(\Ub,\delta)$ approximation of
\,$\fb$ on $\Kb$ with respect to $\eta$. The other way round, if
$\fb$ is measurable and $f$ is any $(\Ub,\delta)$ approximation of
\,$\fb$ on $\Kb$ with respect to $\eta$, then it as a
$(\Ub,\delta_1)$ lifting of $\fb$, where
$$
\delta_1 = \sup_{\xb \in \Kb}
\left(\biggl(1 + \frac{\lv\eta^{-1}[\xb + \Ub]\rv}
   {\lv\eta^{-1}[\Ub]\rv}\biggr)\,\delta
+ \biggl|1 - \frac{\lv\eta^{-1}[\xb + \Ub]\rv}
   {\lv\eta^{-1}[\Ub]\rv}\biggr| \lv\fb(\xb)\rv\right)
$$
The straightforward verification of this fact is left to the reader.
It should be noted that if $\eta$ is a $(\Kb',\Vb)$ approximation
of \,$\Gb$, where $\Kb + \Ub \sbs \Kb'$ and $\Vb \sbs \Ub$
is small enough, then the quotients
$\lv\eta^{-1}[\xb + \Ub]\rv\!\bigl/\lv\eta^{-1}[\Ub]\rv$
can be made arbitrarily close to~1 for $\xb \in \Kb$.
\smallskip

Let us additionally introduce some sets of functions $G \to \bC$
approximatively related to the linear spaces $\MM(G,\Go,\Gf)$,
$\LL^p(G,\Go,\Gf)$, $\CCbu(G,\Go)$ and $\CCo(G,\Go,\Gf)$,
respectively, defined for IMG triplets $(G,\Go,\Gf)$.
Given a finite $(\Kb,\Ub)$ approximation $\eta\:G \to \Gb$,
a $p \in [1,\infty)$, a neighborhood $\Vb \sbs \Ub$ of \,$0 \in \Gb$
and $\delta,\eps > 0$, we denote
$$\align
\MM(G,\eta,\Ub,\Kb,\delta) &= \biggl\{f \in \bC^G\:
   \frac{\mb(\Ub)}{\lv\eta^{-1}[\Ub]\rv}\,
      \sum_{a \in \eta^{-1}[\Gb\sms\Kb]} \lv f(a)\rv
         \le \delta \biggr\}\,,\\
\LL^p(G,\eta,\Ub,\Kb,\delta,\Vb,\eps) &= \biggl\{f \in \bC^G\:
   \lv f\rv^p \in \MM\bigl(G,\eta,\Ub,\Kb,\delta^p\bigr)
\ \& \\
&\qquad \bigl(\all a \in \eta^{-1}[\Vb]\bigr)
\biggl(\frac{\mb(\Ub)}{\lv\eta^{-1}[\Ub]\rv}\,
\sum_{x \in G} \lv f_a(x) - f(x)\rv^p \le \eps^p\biggr)\biggr\}\,,\\
\CCbu(G,\eta,\Kb,\Vb,\eps) &=\ \bigr\{f \in \bC^G; \\
&\qquad \bigl(\all x,y \in \eta^{-1}[\Kb]\bigr)
   \bigl(\eta(x) - \eta(y) \in \Vb \imp
      \lv f(x) - f(y)\rv \le \eps\bigr)\bigl\}\,,\\
\CCo(G,\eta,\Kb,\delta,\Vb,\eps) &=
   \bigl\{f \in \CCbu(G,\eta,\Kb,\Vb,\eps)\:
      \bigl(\all x \in \eta^{-1}[\Gb\sms\Kb]\bigr)\bigl(\lv f(x)\rv
         \le \delta\bigr)\bigr\}\,.
\endalign$$

The following two propositions could be proved directly in
a standard way, and the corresponding nonstandard results on
existence of liftings could be obtained from them. We, in turn,
will derive them from their nonstandard counterparts. Though
this could be achieved by applying Nelson's translation algorithm,
we will provide a more detailed argumentation.

\proclaim{3.4.1.~Proposition}
Let \,$\Gb$ be an LCA group and \,$\mmu \in \Mb(\Gb)$ a complex
regular Borel measure on $\Gb$. Let further $\Ub$ be a neighborhood
of \,$0 \in \Gb$, $\Kb \sbs \Gb$ be a compact set such that
\,$\Ub \sbs \Kb$, and \,$\delta > 0$. Then there exist a finite
abelian group $G$, a $(\Kb,\Ub)$ approximation $\eta\:G \to \Gb$
and a weak $(\Ub,\delta)$ lifting $f\:G \to \bC$ of \,$\mmu$ on
$\Kb$ with respect to $\eta$. If additionally
$\lv\mmu\rv(\Gb \sms \Kb) < \delta$, then one can guarantee that
$f \in \MM(G,\eta,\Ub,\Kb,\delta)$.
\endproclaim

\demo{Proof}
Let \,$\eta\:G \to \sGb$ be any HFI approximation of \,$\Gb$ by
a hyperfinite abelian group $G$ and $(G,\Go,\Gf)$ be the IMG
group triplet arising from $\eta$. Then there is a weak lifting
$f \in \MM(G,\Go,\Gf)$ of \,$\mmu$. Obviously, $\eta$ is a
$(\sKb,\sUb)$ approximation of \,$\sGb$ and \,$f$ is a weak
$(\sUb,\delta)$ lifting of \,${}^{*\!}\mmu$ on $\sKb$ with
respect to $\eta$. By the {\it transfer principle}, there exist
a $(\Kb,\Ub)$ approximation $\eta\:G \to \Gb$ and of $\Gb$ by
a finite abelian group $G$ and a weak $(\Ub,\delta)$ lifting
$f\:G \to \bC$ of \,$\mmu$ on $\Kb$ with respect to $\eta$.

If, additionally, the variation $\lv\mmu\rv$ of $\mmu$ is
$\delta$-concentrated on $\Kb$, then it is clear that the lifting
$f \in \MM(G,\Go,\Gf)$ belongs to
$\MM(G,\eta,\sUb,\sKb,\delta)$. Then the existence of a
finite $(\Kb,\Ub)$ approximation $\eta\:G \to \Gb$ and a weak
\,$(\Ub,\delta)$ lifting $f \in \MM(G,\eta,\Ub,\Kb,\delta)$ of
\,$\mmu$ on $\Kb$ with respect to $\eta$ follows from the
{\it transfer principle}, again.
\enddemo

\proclaim{3.4.2.~Proposition}
Let \,$\Gb$ be an LCA group, $1 \le p < \infty$ and
\,$\fb \in \Lbp(\Gb)$. Let further $\Ub$ be a neighborhood of
\,$0 \in \Gb$, $\Kb \sbs \Gb$ be a compact set such that
\,$\Ub \sbs \Kb$, and \,$\delta > 0$. Then there exist a finite
abelian group $G$, a $(\Kb,\Ub)$ approximation $\eta\:G \to \Gb$
and a $(\Ub,\delta)$ lifting $f\:G \to \bC$ of \,$\fb$ on
$\Kb$ with respect to $\eta$. If additionally
$$
\LV\fb \cd 1_{\Gb \sms \Kb}\RV_{_p} < \delta
\qquad\text{and}\qquad
\LV\fb_{\ab} - \fb\RV_{_p} < \eps
$$
for all \,$\ab$ from some neighborhood \,$\Vb \sbs \Ub$ of
\,$0$ in $\Gb$ and an $\eps > 0$, then one can guarantee that
\,$\eta$ is a $(\Kb,\Vb)$ approximation of \,$\Gb$  and
\,$f \in \LL^p(G,\eta,\Ub,\Kb,\delta,\Vb,\eps)$, as well.
\endproclaim

\demo{Proof}
The first statement follows right away from Proposition~3.4.1.
Moreover, the lifting $f$ of $\fb$ with respect to the HFI
approximation $\eta$ belongs to $\LL^p(G,\Go,\Gf)$.

If, additionally, $\fb$ is $\delta$-concentrated on $\Kb$ with
respect to the norm $\LV\cd\RV_{_p}$ and $\Vb \sbs \Ub$ is such
a neighborhood of \,$0 \in \Gb$ that
$\LV\fb_{\ab} - \fb\RV_{_p} < \eps$ for $\ab \in \Vb$, then $\eta$
is even a $(\sKb,\sVb)$ approximation of $\sGb$ and the
lifting $f \in \LL^p(G,\Go,\Gf)$ obviously belongs to
$\LL^p(G,\eta,\sUb,\sKb,\delta,\sVb,\eps)$. The
existence of a finite $(\Kb,\Vb)$ approximation
$\eta\:G \to \Gb$ and a $(\Ub,\delta)$ lifting
$f \in \LL^p(G,\eta,\Ub,\Kb,\delta,\Vb,\eps)$ of \,$\fb$ on $\Kb$
with respect to $\eta$ follows from the {\it transfer principle},
once again.
\enddemo

For completeness' sake let us record also the following obvious

\proclaim{3.4.3.~Proposition}
Let \,$\Gb$ be an LCA group and \,$\fb \in \Cb(\Gb)$. Let further
$\Ub$ be a neighborhood of \,$0 \in \Gb$, $\Kb \sbs \Gb$ be a
compact set such that \,$\Ub \sbs \Kb$ and \,$\eta\:G \to \Gb$ be
a  $(\Kb,\Ub)$ approximation of \,$\Gb$ by a finite abelian group
$G$. Then, for any \,$\delta,\eps > 0$, we have:
{\parindent 21pt
\item{\sl (a)}
if \,$\lv\fb(\xb) - \fb(\yb)\rv \le \eps$ for all \,$\xb \in \Kb$,
$\yb \in \xb + \Ub$, then the composed mapping $f = \fb \co \eta$ is
a $(\Ub,\eps)$ approximation of \,$\fb$ on $\Kb$ with respect to
$\eta$ and belongs to $\CCbu(G,\eta,\Kb,\Ub,\eps)$;
\item{\sl (b)}
if additionally \,$\lv\fb(\xb)\rv \le \delta$ \,for all
\,$\xb \in \Gb \sms \Kb$, then even
\,$f \in \CCo(G,\eta,\Kb,\delta,\Ub,\eps)$.
\item{}
\vskip-9pt}
\endproclaim

Our starting point is the approximation of the Fourier-Stieltjes
transform on LCA groups by the discrete Fourier transform on finite
abelian groups.

Let $\mmu \in \Mb(\Gb)$ be a complex regular Borel measure on $\Gb$.
Then its Fourier-Stieltjes transform $\wh\mmu \in \Cbu(\dGb)$ is a
bounded uniformly continuous function on the dual group~$\dGb$. We
will approximate $\wh\mmu$ by the discrete Fourier transform
$\wh f$ of a function $f\:G \to \bC$ on some finite abelian
group~$G$. The ``parameters of the approximation'' given in advance
are: a compact set $\Gamb_0 \sbs \dGb$ of positive Haar measure and
an $\eps > 0$. Then there is a relatively compact neighborhood
$\Omb_0$ of \,$1 \in \dGb$ such that, for any $\ggam \in \Gamb_0$,
$\cchi \in \ggam\,\Omb_0$,
$$
\lv\wh\mmu(\ggam) - \wh\mmu(\cchi)\rv \le \eps\,.
$$
Our goal is achieved once we find a ``sufficiently good''
pair of adjoint finite approximations $\eta\:G \to \Gb$,
$\phi\:\dG \to \dGb$, such that $\phi$ is a
$(\Gamb,\Omb)$ approximation of \,$\dGb$ for some neighborhood
$\Omb \sbs \Omb_0$ of \,$1 \in \dGb$ and a compact set
$\Gamb \sbs \dGb$ such that $\Gamb_0\,\Omb_0 \sbs \Gamb$, and
formulate some ``reasonable'' conditions making sure that for any
function $f\:G \to \bC$ which is such a ``nice'' approximation of
the measure $\mmu$, its Fourier transform $\wh f$ is an
$(\Omb,\eps)$ approximation of \,$\wh\mmu$ on $\Gamb_0$ with
respect to $\phi$.

\proclaim{3.4.4.~Finite Fourier-Stieltjes Transform Approximation
Theorem}
Let \,$\Gb$ be a Hausdorff LCA group and \,$\mmu \in \Mb(\Gb)$ be
a complex regular Borel measure on $\Gb$. Let further $\eps > 0$,
$0 < \alfa \le \pi/3$, $\Gamb_0 \sbs \dGb$ be a compact set of
positive Haar measure and \,$\Omb_0$ be a relatively compact
neighborhood of \,$1$ in $\dGb$ such that
$\lv\wh\mmu(\ggam) - \wh\mmu(\cchi)\rv \le \eps$ whenever
$\ggam \in \Gamb_0$, $\cchi \in \ggam\,\,\Omb_0$. Then there
exist $\alfa$-adjoint pairs $(\Kb,\Ub)$, $(\Gamb,\Omb)$ of subsets
in $\Gb$ and its dual group $\dGb$, respectively, such that
$\Omb \sbs \Omb_0$, $\Gamb_0\,\Omb_0 \sbs \Gamb$, a
$\delta \in (0,\alfa)$, a finite abelian group $G$, and a strongly
$(\alfa,\delta)$-adjoint pair of approximations $\eta\:G \to \Gb$,
$\phi\:\dG \to \dGb$ of \,$\Gb$, $\dGb$, respectively, with respect
to $(\Kb,\Ub)$, $(\Gamb,\Omb)$, such that, for any function
$f \in \MM(G,\eta,\Ub,\Kb,\delta)$ which is a weak
$(\Ub,\delta)$ lifting of the measure $\mmu$ on $\Kb$ with
respect to $\eta$, its Fourier transform $\wh f$ is an
$(\Omb,\eps)$ approximation of the Fourier-Stieltjes transform
$\wh\mmu$ on $\Gamb_0$ with respect to $\phi$.
\endproclaim

\demo{Proof}
Let $\bigl((\Kb_i,\Ub_i)\bigr)_{i \in I}$,
$\bigl((\Gamb_i,\Omb_i)\bigr)_{i \in I}$ be $\alfa$-adjoint DD~bases in
$\Gb$, $\dGb$, respectively, over some directed poset $(I,\le)$,
such that $\Omb_i \sbs \Omb_0$ and $\Gamb_0\,\Omb_0 \sbs \Gamb_i$ for
each $i \in I$. Let further $(\delta_n)_{n\in\bN}$ be a decreasing
sequence of positive reals $< \alpha$ converging to $0$. We can assume,
without loss of generality, that $(I,\le)$ is regular; then the set
$J = I \cx \bN$ ordered componentwise is an RCID poset, so that
there exists a $\kappa$-good directed ultrafilter on $J$, where
$\kappa = (\card I)^+$\,---\,see the comment preceding Theorem~2.5.6.
According to Theorem~2.5.9 (and using the convention formulated prior to
Proposition~2.5.5), there is an adjoint pair of approximating systems
\hbox{$\bigl(\eta_{in}\:G_{in} \to \Gb\bigr)_{(i,n) \in J}$,}
$\bigl(\phi_{in}\:\dG_{in} \to \dGb\bigr)_{(i,n) \in J}$
of \,$\Gb$, $\dGb$, respectively, by finite abelian groups, well
based with respect to the DD bases
$\bigl((\Kb_{in},\Ub_{in})\bigr)_{(i,n) \in J}$,
$\bigl((\Gamb_{in},\Omb_{in})\bigr)_{(i,n) \in J}$,
such that each particular pair $\eta_{in}$, $\phi_{in}$
is strongly $(\alpha,\delta_n)$-adjoint with respect to
$(\Kb_i,\Ub_i)$, $(\Gamb_i,\Omb_i)$

Let us assume, for contradiction, that the conclusion of the Theorem
is not true. Then, for any $(i,n) \in J$, there is an
$f_{in} \in \MM(G_{in},\eta_{in},\Ub_i,\Kb_i,\delta_n)$ which is a weak
$(\Ub_i,\delta_n)$ lifting of the measure $\mmu$ on $\Kb_i$ with
respect to $\eta_{in}$, however, $\wh f_{in}$ \,is not an
$(\Omb_i,\eps)$ approximation of \,$\wh\mmu$ on $\Gamb_0$ with
respect to $\phi_{in}$, i.e., there are a $\ggam_{in} \in \Gamb_0$ and
a $\chi_{in} \in \dG_{in}$,  such that
$\phi_{in}(\chi_{in}) \in \ggam_{in}\,\Omb_i$, but
$$
\bigl|\wh\mmu(\ggam_{in}) - \wh f_{in}(\chi_{in})\bigr| > \eps\,.
$$
Let $\DD$ be some $\kappa$-good directed ultrafilter on $(J,\le)$.
Forming the ultraproducts $G = \prod_{(i,n) \in J} G_{in}/\DD$ and
$\sGb = \Gb^J\!/\DD$, we can identify
$\dG = \prod_{(i,n) \in J} \dG_{in}/\DD$, as well as
$\wh{\sGb} = \sdGb = \dGb^J\!/\DD$. By Proposition 2.5.5(a),
we obtain an adjoint pair $\eta = (\eta_{in})_{(i,n) \in J}/\DD$,
\hbox{$\phi = (\phi_{in})_{(i,n) \in J}/\DD$} of HFI approximations
$\eta\:G \to \sGb$, $\phi\:\dG \to \sdGb$ of \,$\Gb$, $\dGb$,
respectively, giving rise to dual IMG triplets $(G,\Go,\Gf)$,
$\bigl(\dG,\Gf^{\ort},\Go^{\ort}\bigr)$ within the $\kappa$-saturated
nonstandard universe obtained via the ultraproduct construction
modulo~$\DD$. Then, as easily seen, the internal function
$f = (f_{in})_{(i,n) \in J}/\DD\:G \to \sbC$ is a weak lifting of the
measure $\mmu$ and belongs to $\MM(G,\Go,\Gf)$. By Theorem~3.3.3,
its Fourier transform $\wh f \in \CCbu\bigl(\dG,\Gf^{\ort}\bigr)$
is an $S$-continuous lifting of the function
$\wh\mmu \in \Cbu\bigl(\dGb\bigr)$, i.e.,
$\wh f(\chi) \apr \wh\mmu(\ggam)$ whenever $\chi \in G$, $\ggam\in \dGb$
satisfy $\phi(\chi) \apr \ggam$. At the same time, for the particular
choice of \,$\chi = (\chi_{in})_{(i,n) \in J}/\DD \in \dG$ and
$\ggam = {}^\co\bigl((\ggam_{in})_{(i,n) \in J}/\DD\bigr) \in \Gamb_0$,
we have $\phi(\chi) \apr \ggam$ and
$\bigl|\wh\mmu(\ggam) - \wh f(\chi)\bigr| \ge \eps$,
which is a contradiction.
\enddemo

If \,$\fb \in \Lb^1(\Gb)$, then the last Theorem applies to the
measure $\mmu \in \Mb(\Gb)$ such that $\dif\mmu = \fb\dif\mb$.
However, as $\wh\fb \in \Cbo(\dGb)$, there is  a compact set
$\Gamb_0 \sbs \dGb$ such that $\wh\fb$ \,is $\eps$-concentrated
on $\Gamb_0$ with respect to the norm $\LV\cd\RV_{_\infty}$, i.e.,
$$
\bigl|\wh\fb(\cchi)\bigr| < \eps
$$
for $\cchi \in \dGb \sms \Gamb_0$. Thus it is possible to give a
``globalized'' version of Theorem~3.4.3 in this case.

\proclaim{3.4.5.~Finite Fourier Transform Approximation Theorem}
Let \,$\Gb$ be a Hausdorff LCA group and \,$\fb \in \Lb^1(\Gb)$.
Let further $\eps > 0$, $0 < \alfa \le \pi/3$, $\Gamb_0 \sbs \dGb$
be a compact set of positive Haar measure such that $\wh\fb$ is
$\eps$-concentrated on $\Gamb_0$ with respect to the norm
$\LV\cd\RV_{_\infty}$, and \,$\Omb_0$ be a relatively compact
neighborhood of \,$1$ in $\dGb$ such that
\,$\bigl|\wh\fb(\ggam) - \wh\fb(\cchi)\bigr| \le \eps$
whenever $\ggam \in \Gamb_0$, $\cchi \in \ggam\,\Omb_0$. Then
there exist $\alfa$-adjoint pairs $(\Kb,\Ub)$, $(\Gamb,\Omb)$ of
subsets in $\Gb$ and its dual group $\dGb$, respectively, such
that $\Omb \sbs \Omb_0$, $\Gamb_0\,\Omb_0 \sbs \Gamb$, a
$\delta \in (0,\alfa)$, a~finite abelian group $G$, and a strongly
$(\alfa,\delta)$-adjoint pair of approximations $\eta\:G \to \Gb$,
$\phi\:\dG \to \dGb$ of \,$\Gb$, $\dGb$, respectively, with respect
to $(\Kb,\Ub)$, $(\Gamb,\Omb)$, such that, for any function
$f \in \LL^1(G,\eta,\Ub,\Kb,\delta,\Ub,\delta)$ which is a
$(\Ub,\delta)$ lifting of $\fb$ on $\Kb$ with respect to $\eta$,
its Fourier transform $\wh f$ is an $(\Omb,\eps)$ approximation
of the Fourier transform $\wh\fb$ on $\Gamb_0$ with respect to
$\phi$, and $\bigl|\wh f(\chi)\bigr| \le \eps$ for
$\chi \in \phi^{-1}\bigl[\dGb \sms \Gamb_0\,\Omb\bigr]$.
\endproclaim

\demo{Proof}
Let $\bigl((\Kb_i,\Ub_i)\bigr)_{i \in I}$,
$\bigl((\Gamb_i,\Omb_i)\bigr)_{i \in I}$,
$(\delta_n)_{n\in\bN}$, $J = I \cx \bN$, $\kappa = (\card J)^+$,
as well as the approximating systems
$\bigl(\eta_{in}\:G_{in} \to \Gb\bigr)_{(i,n) \in J}$,
$\bigl(\phi_{in}\:\dG_{in} \to \dGb\bigr)_{(i,n) \in J}$
be as in the proof Theorem~3.4.4. Assume that the conclusion of
Theorem~3.4.4 fails. Then, for any $(i,n) \in J$, there is a function
$f_{in} \in \LL^1(G_{in},\eta_{in},\Ub_i,\Kb_i,\delta_n,\Ub_i,\delta_n)$
which is a $(\Ub_i,\delta_n)$ lifting of $\fb$ on $\Kb_i$ with
respect to $\eta_{in}$, however, $\wh f_{in}$ \,either is not an
$(\Omb_i,\eps)$ approximation of \,$\wh\fb$ on $\Gamb_0$ with
respect to $\phi_{in}$, i.e., there are a $\ggam_{in} \in \Gamb_0$ and
a $\chi_{in} \in \dG_i$,  such that
$\phi_{in}(\chi_i) \in \ggam_{in}\,\Omb_i$, but
$$
\bigl|\wh\fb(\ggam_i) - \wh f_i(\chi_i)\bigr| > \eps\,,
$$
or there is a
$\chi_{in} \in \phi_{in}^{-1}\bigl[\dGb \sms \Gamb_0\,\Omb_i\bigr]$
such that $\bigl|\wh f_{in}(\chi_{in})\bigr| > \eps$. Let $J_1$, $J_2$
denote the subsets of \,$J$ consisting of those pairs $(i,n)$ for which
the first or the second alternative takes place, respectively.

Let us take any $\kappa$-good directed ultrafilter $\DD$ on $(J,\le)$
and form the ultraproducts, as well as the adjoint HFI approximations
$\eta\:G \to \sGb$, $\phi\:\dG \to \sdGb$ as in the proof of
Theorem~3.4.4. Then the internal function
$f = (f_{in})_{(i,n) \in J}/\DD$ belongs to $\LL^1(G,\Go,\Gf)$ and it
is a lifting of $\fb \in \Lb^1(\Gb)$ with respect to $\eta$.
By Theorem~3.3.1, $\wh f \in \CCo\bigl(\dG,\Gf^{\ort},\Go^{\ort})$
and it is an $S$-continuous lifting of $\wh\fb$, i.e.,
$\wh f(\chi) \apr \wh\fb(\ggam)$ if
\,$\phi(\chi) \apr \ggam \in \dGb$, as well as $\wh f(\chi) \apr 0$
if \,$\phi(\chi) \nin \FdG = \Ns\bigl(\sdGb\bigr)$.

As $J_1 \cup J_2 = J$ and $\DD$ is an ultrafilter, we have
either $J_1 \in \DD$ or $J_2 \in \DD$. If \,$J_1 \in \DD$, then, for
the particular choice
$\ggam = {}^\co\bigl((\ggam_{in})_{(i,n) \in J}/\DD\bigr) \in \dGb$,
$\chi = (\chi_{in})_{(i,n) \in J}/\DD \in \dG$ (with arbitrary
$\ggam_{in}$, $\chi_{in}$ for $(i,n) \in J \sms J_1$), we have
$\phi(\chi) \apr \ggam \in \Gamb_0$ and
$\bigl|\wh\fb(\ggam) - f(\chi)\bigr| \ge \eps$. If \,$J_2 \in \DD$,
then, for $\chi = (\chi_{in})_{(i,n) \in J}/\DD$ (with arbitrary $\chi_{in}$
if \,$(i,n) \in J \sms J_2$), we have
$\phi(\chi) \in \dGb \sms \Gamb_0$ and $\bigl|\wh f(\chi)\bigr| \ge \eps$.
Both possibilities lead to contradictions.
\enddemo

If $1 < p \le 2 \le q < \infty$ are dual exponents, then similar
accounts lead us to the formulation of the following approximation
theorem. Let $\nb = \mb_{\dGb}$ denote the Haar measure on the dual
group $\dGb$ properly normalized to make the Fourier inversion
formula to hold.

\proclaim%
{3.4.6.~Finite Generalized Fourier Transform Approximation Theorem}
Let \,$\Gb$ be a Hausdorff LCA group, $1 < p \le 2 \le q < \infty$
be dual exponents and \,$\fb \in \Lb^p(\Gb)$. Let further $\eps > 0$,
$0 < \alfa \le \pi/3$, $\Gamb_0 \sbs \dGb$ be a compact set of positive
Haar measure such that $\wh\fb$ is $\eps$-concentrated on $\Gamb_0$
with respect to the norm $\LV\cd\RV_{_q}$, and \,$\Omb_0$ be a
relatively compact neighborhood of \,$1$ in $\dGb$ such that
\,$\bigl\|\wh\fb_{\,\oom} - \wh\fb\,\bigr\|_{_q} \le \eps$ whenever
$\oom \in \Omb_0$. Then there exist $\alfa$-adjoint pairs
$(\Kb,\Ub)$, $(\Gamb,\Omb)$ of subsets in $\Gb$ and its dual
group $\dGb$, respectively, such that $\Omb \sbs \Omb_0$,
$\Gamb_0\,\Omb_0 \sbs \Gamb$, a $\delta \in (0,\alfa)$,
a finite abelian group $G$, and a strongly $(\alfa,\delta)$-adjoint
pair of approximations $\eta\:G \to \Gb$, $\phi\:\dG \to \dGb$ of
\,$\Gb$, $\dGb$, respectively, with respect to $(\Kb,\Ub)$,
$(\Gamb,\Omb)$, such that, for any function
$f \in \LL^p(G,\eta,\Ub,\Kb,\delta,\Ub,\delta)$ which is a
$(\Ub,\delta)$ lifting of $\fb$ on $\Kb$ with respect to $\eta$,
its Fourier transform $\wh f$ is an $(\Omb,\eps)$ lifting
of the Fourier transform $\wh\fb$ on $\Gamb_0$ with respect to
$\phi$, and
$$
\frac{\nb(\Omb)}{\lv\phi^{-1}[\Omb]\rv}\,
   \sum_{\chi \in \phi^{-1}[\dGb \sms \Gamb_0\,\Omb]}
      \bigl|\wh f(\chi)\bigr|^q \le \eps^q\,.
$$
\endproclaim

\demo{Proof}
Let us start in the same way as in the proofs of Theorems~3.4.4,
3.4.5 and assume that the conclusion of Theorem~3.4.6 fails.
Then, for any $(i,n) \in J = I \cx \bN$, there is a function
$f_{in} \in \LL^p(G_{in},\eta_{in},\Ub_i,\Kb_i,\delta_n,\Ub_i,\delta_n)$
which is a $(\Ub_i,\delta_n)$ lifting of $\fb$ on $\Kb_i$ with
respect to $\eta_{in}$, however, $\wh f_{in}$ \,either is not an
$(\Omb_i,\eps)$ lifting of \,$\wh\fb$ on $\Gamb_0$ with respect to
$\phi_{in}$, i.e., (due to the regularity of $\nb$) there is a compact
set $\Delb_{in} \sbs \Gamb_0$ such that $\nb(\Delb_{in}) > \eps$ and
$$
\lv\,\frac{1}{\nb(\Omb_i)}\int_{\ggam\,\Omb_i} \wh\fb \dif\nb
- \frac{1}{\lv\phi_{in}^{-1}[\Omb_i]\rv}\,
   \sum_{\chi \in \phi_{in}^{-1}[\ggam\,\Omb_i]} \wh f_{in}(\chi)\,\rv
     > \eps
$$
for each $\ggam \in \Delb_{in}$, or
$$
\frac{\nb(\Omb_i)}{\lv\phi_{in}^{-1}[\Omb_i]\rv}\,
   \sum_{\chi \in \phi_{in}^{-1}[\dGb \sms \Gamb_0\,\Omb_i]}
      \bigl|\wh f_{in}(\chi)\bigr|^q > \eps^q\,.
$$
Again, let $J_1$, $J_2$ denote the subsets of \,$J$ consisting of those
pairs $(i,n)$ for which the first or the second alternative takes place,
respectively.

Let $\DD$ be a $\kappa$-good directed ultrafilter on $(J,\le)$
and $\eta\:G \to \sGb$, $\phi\:\dG \to \sdGb$ be the adjoint HFI
approximations as in the proof of Theorem~3.4.4. Then the internal
function $f = (f_{in})_{(i,n) \in J}/\DD$ belongs to $\LL^p(G,\Go,\Gf)$
and it is a lifting of $\fb \in \Lb^p(\Gb)$ with respect to $\eta$.
By Theorem~3.3.2, $\wh f \in \LL^q\bigl(\dG,\Gf^{\ort},\Go^{\ort}\bigr)$
and it is a lifting of $\wh\fb$, i.e.,
$\wh f(\chi) \apr \wh\fb(\ggam)$ whenever
$\phi(\chi) \apr\ggam \in \dGb$, for almost all $\chi \in \Gf^{\ort}$
with respect to the Loeb measure $\lam_{\hat d}$, where $\hat d$ is
any normalizing multiplier for the triplet
$\bigl(\dG,\Gf^{\ort},\Go^{\ort}\bigr)$. At the same time
$$
\hat d\,\sum_{\chi \in \phi^{-1}[\Thb]}
              \bigl|\wh f(\chi)\bigr|^q \apr 0
$$
for any internal set $\Thb \sbs \sdGb \sms \FdG$.

If \,$J_1 \in \DD$, then we form the internal sets
$$\gather
\Delb = \prod_{(i,n) \in J} \Delb_{in}/\DD \sbs {}^*\Gamb_0
= \Gamb_0^{\,J}/\DD \sbs \FdG\,, \\
\Omb = \prod_{(i,n) \in J} \Omb_i/\DD \sbs \IdG
\endgather$$
(with $\Delb_{in} \sbs \dGb$ arbitrary for $(i,n) \in J \sms J_1$),
as well as the observable trace ${}^{\co\!}\Delb \sbs \Gamb_0$.
Further we put
$\Dela_{in} = \phi_{in}^{-1}[\Delb_{in}\,\Omb_i] \sbs G_{in}$
for any $(i,n) \in J$, and form the internal set
$$
\Dela = \prod_{(i,n) \in J} \Dela_{in}/\DD = \phi^{-1}[\Delb]
        \sbs \Go^{\ort}\,.
$$
Then ${}^{\co\!}\Delb = \Dela^{\kr}$, hence ${}^{\co\!}\Delb$ is
Borel and
$$
\nb({}^{\co\!}\Delb) = \lam_{\hat d}\bigl(\Dela^{\kr}\bigr)
                     \ge \eps\,,
$$
where $\hat d$ is a particular normalizing multiplier on $\dG$
for which the last formula is valid.
For any $\ggam \in \dGb$ let us form the hypercomplex numbers
$$\gather
A(\ggam) = \left(\frac{1}{\nb(\Omb_i)}\int_{\ggam\,\Omb_i}
     \wh\fb \dif\nb\right)_{(i,n) \in J}\bigg/\DD\,, \\
\\
B(\ggam) = \biggl(\,\frac{1}{\lv\phi_{in}^{-1}[\Omb_i]\rv}\,
    \sum_{\chi \in \phi_{in}^{-1}[\ggam\,\Omb_i]}
      \wh f_{in}(\chi)\biggr)_{(i,n) \in J}\bigg/\DD
= \frac{1}{\lv\phi^{-1}[\Omb]\rv}
   \sum_{\chi \in \phi^{-1}[\ggam\,\Omb]} \wh f(\chi)\,.
\endgather$$
By {\L}os theorem, $\lv A(\ggam) - B(\ggam)\rv \ge \eps$ for
each $\ggam \in {}^\co\Delb$. On the other hand, for $\nb$-almost
all $\ggam \in \dGb$ we have $\wh\fb(\ggam) \apr A(\ggam)$.
At the same time, for $\lam_{\hat d\,}$-almost all
$\chi \in \Go^{\ort}$ we have $f(\chi) \apr \wh f(\ggam)$
whenever $\phi(\chi) \apr \ggam \in \dGb$. Since also
$$
\frac{\lv\phi^{-1}[\ggam\,\Omb]\rv}{\lv\phi^{-1}[\Omb]\rv} \apr 1
$$
for $\ggam \in \dGb$, we finally get $B(\ggam) \apr \wh\fb(\ggam)$
for almost all $\ggam \in \dGb$, which is a contradiction.

If \,$J_2 \in \DD$, then the function
$f = (f_{in})_{(i,n) \in J}/\DD \in \sbC^G$ satisfies
$$
\frac{{}^*\nb(\Omb)}{\lv\phi^{-1}[\Omb]\rv}\,
   \sum_{\chi \in \phi^{-1}[\dGb \sms \Gamb_0\,\Omb]}
      \bigl|\wh f(\chi)\bigr|^q > \eps^q\,,
$$
contradicting that $\wh f \in \LL^q\bigl(\dG,\Gf^{\ort},\Go^{\ort}\bigr)$
is a lifting of \,$\wh\fb \in \Lbq\bigl(\dGb\bigr)$ and $\wh\fb$ is
$\eps$-concentrated on $\Gamb_0$ with respect to $\LV\cd\RV_{_q}$.
\enddemo

In view of the massive formulations of Theorems~3.4.1--3 we cannot
spare some remarks.

\remark{Remark\/ \rm 1}
It seems that all of the last three Theorems contain some
``overkill.'' It is well possible that in many concrete cases
not all properties of the objects guaranteed to exist will prove
to be useful. One of the candidates to be reduced is the strong
$(\alfa,\eps)$-adjointness of the approximations $\eta$, $\phi$
with respect to the pairs $(\Kb,\Ub)$, $(\Gamb,\Omb)$.
\endremark
\medskip

\remark{Remark\/ \rm 2}
All the three Theorems are purely existential and give neither any
estimation of the objects to be found in terms of the given ones
nor a hint how to find them. This is the common disadvantage of
many ``soft'' results obtained by means of nonstandard analysis.
On the other hand, they indicate that it makes sense to look for
the corresponding ``hard'' counterparts, for instance for estimates
of \,$\delta$ in terms of $\eps$ (or vice versa), and the like.
\endremark
\medskip

In order to illustrate what we have in mind in the above Remarks,
we formulate two fairly general accompanying examples to
Theorems~3.4.4 and 3.4.5, in which the ``soft'' (i.e., existential)
parts are reduced to the standard equivalents of Theorem~3.3.4 and
its Corollary~3.3.5, and some explicit bounds are given, as well.
The second of these examples will be representative enough to partly
take care of Theorem~3.4.6, as well.

Instead of emphasizing that our constructions work for every
$\alfa \in (0,\pi/3]$ or even for $\alfa \in (0,2\pi/3)$, as we
used to do, we will take the advantage of choosing $\alfa > 0$ as
small as we please, now. In the first of our examples, instead of
using the strong $(\alfa,\eps)$-adjointness of the approximations
$\eta$, $\phi$, we will even manage with a weaker property.

Given a Hausdorff LCA group $\Gb$, an $\alfa \in (0,\pi/3]$
and $\alfa$-adjoint pairs of sets $(\Kb,\Ub)$ in $\Gb$ and
$(\Gamb,\Omb)$ in $\dGb$, we say that a $(\Kb,\Ub)$ approximation
$\eta\:G \to \Gb$ and a $(\Gamb,\Omb)$ approximation
$\phi\:\dG \to \dGb$ are {\it $\alfa$-adjoint\/} if
$$
\lv\arg\frac{(\phi\,\gama)(\eta\,a)}{\gama(a)}\rv \le \alfa
$$
for all $a \in \eta^{-1}[\Kb]$, $\gama \in \phi^{-1}[\Gamb]$.

For every finite $(\Kb,\Ub)$ approximation $\eta\:G \to \Gb$ we
denote $\FF_{\eta,d}\:\bC^G \to \bC^{\dGb}$ the linear operator,
called the {\it modified Fourier transform}, given by
$$
\FF_{\eta,d}(f)(\cchi) = \br f, \cchi \co \eta\kt_d
= d \,\sum_{a \in G} f(a)\,\ovl\cchi(\eta\,a)\,,
$$
for $f \in \bC^G$, $\cchi \in \dGb$, where $d > 0$ is some
scaling coefficient to be determined subsequently. In this
context we always assume that
$$
d = \frac{\mb(\Ub)}{\lv\eta^{-1}[\Ub]\rv}\,,
$$
where $\mb$ is a fixed Haar measure on $\Gb$.

\example{3.4.7.~Example}
Let $\Gb$ be a Hausdorff LCA group, $\mmu \in \Mb(\Gb)$ be a
complex regular Borel measure on $\Gb$ with finite variation
$\LV\mmu\RV$, \,$\Gamb_0 \sbs \dGb$ be a compact set of positive
Haar measure and $\alfa \in (0,\pi/3]$. Translating Theorem~3.3.4
into standard terms it follows that there exist a compact set
$\Kb_0 \sbs \Gb$, a neighborhood $\Ub_0$ of \,$0$ in $\Gb$ and
a $\delta > 0$ such that, for any compact set $\Kb \sbs \Gb$
and a neighborhood $\Ub$ of \,$0 \in \Gb$, the inclusions
$\Ub \sbs \Ub_0$, $\Kb_0 \sbs \Kb$ imply that, for every finite
$(\Kb,\Ub)$ approximation $\eta\:G \to \Gb$ and each weak
$(\Ub,\delta)$ lifting $f \in \MM(G,\eta,\Ub,\Kb,\delta)$ of
\,$\mmu$ on $\Kb$, we have
$$
\bigl|\wh\mmu(\cchi) - \FF_{\eta,d}(f)(\cchi)\bigr| \le \alfa
$$
for all \,$\cchi \in \Gamb_0$. Let $\Omb_0 = \Bohr_\alfa(\Kb_0)$.

Due to Lemma~2.5.6, there exist $\alfa$-adjoint pairs
$(\Kb,\Ub)$, $(\Gamb,\Omb)$ such that $\Ub \sbs \Ub_0$,
$\Kb_0 \sbs \Kb$, as well as $\Omb \sbs \Omb_0$,
$\Gamb_0\,\Omb_0 \sbs \Gamb$. By Proposition~3.4.1 there is
indeed a weak $(\Ub,\delta)$ lifting
\,$f \in \MM(G,\eta,\Ub,\Kb,\delta)$ of \,$\mmu$ on $\Kb$.
According to Proposition~1.4.3 we can assume that
$$
\LV f\RV_{_1} = d \,\sum_{a \in G} \lv f(a)\rv \le \LV\mmu\RV\,.
$$

Now, assume that $\phi\:\dG \to \dGb$ is a $(\Gamb,\Omb)$
approximation of the dual group $\dG$, $\alfa$-adjoint to $\eta$.
Let $\cchi \in \Gamb_0$, $\gama \in \dG$ be such that
$\phi(\gama) \in \cchi\,\Omb \sbs \Gamb$.
Then for each $a \in \eta^{-1}[\Kb]$ we have
$\eta(a) \in \Kb = \Bohr_\alfa(\Omb)$, consequently,
$$
\lv\arg\frac{\cchi(\eta\,a)}{\gama(a)}\rv \le
\lv\arg\frac{\cchi(\eta\,a)}{(\phi\,\gama)(\eta\,a)}\rv
+ \lv\arg\frac{(\phi\,\gama)(\eta\,a)}{\gama(a)}\rv
\le 2\alfa\,,
$$
hence
$\lv\cchi(\eta\,a) - \gama(a)\rv \le 2\sin\alfa \le 2\alfa$.
Further, we have
$$
\bigl|\wh\mmu(\cchi) - \wh f(\gama)\bigr| \le
\bigl|\wh\mmu(\cchi) - \FF_{\eta,d}(f)(\cchi)\bigr| +
\bigl|\FF_{\eta,d}(f)(\cchi) - \wh f(\gama)\bigr|\,,
$$
and the first summand is \,$\le \alfa$. Denoting
$K = \eta^{-1}[\Kb]$, we get for the second summand
$$\align
\bigl|\FF_{\eta,d}(f)(\cchi) - \wh f(\gama)\bigr|
&= \bigl|\br f,\ovl\cchi \co \eta - \ovl\gama\,\kt_d\bigr| \\
&\le \LV f\RV_{_1} \LV(\cchi\co\eta - \gama)
   \cd 1_K\RV_{_\infty}
+ \LV f \cd 1_{G \sms K}\RV_{_1}
   \LV\cchi\co\eta - \gama\RV_{_\infty} \\
&\le 2\alfa\LV\mmu\RV + 2\delta\,.
\endalign$$
Finally, as we obviously can choose $\delta \le \alfa$,
$$
\bigl|\wh\mmu(\cchi) - \wh f(\gama)\bigr| \le
\bigl(2\LV\mmu\RV + 1\bigr)\alfa + 2\delta
 \le \bigl(2\LV\mmu\RV + 3\bigr)\alfa\,.
$$
If we were given an $\eps > 0$ in advance, it is always
possible to arrange that the right-hand expression is
\,$ \le \eps$ by choosing $\alfa$ small enough.
\endexample

\example{3.4.8.~Example}
If \,$\Gb$ is a Hausdorff LCA group, $\fb \in \Lb^1(\Gb)$,
\,$\Gamb_0 \sbs \dGb$ is a compact set of positive Haar measure and
$\alfa \in (0,\pi/3]$, then the previous Example applies directly
to the measure $\mmu$ such that $\fb = \dif\mmu/\!\dif\mb$. Let us
additionally assume that $\bigl|\wh\fb(\cchi)\bigr| \le \eps$ for
 some ``small'' $\eps > 0$  given in advance and all
$\cchi \in \dGb \sms \Gamb_0$. To avoid trivialities we assume
that $\LV\fb\RV_{_1} \ne 0$, i.e., $\fb$ is not equal to $0$
almost everywhere.

However, using Corollary~3.3.5 instead of Theorem~3.3.4 and
translating it into standard terms it follows that there exist
a compact set $\Kb_0 \sbs \Gb$, a neighborhood $\Ub_0$ of \,$0$
in $\Gb$ and a $\delta > 0$ such that, for any compact set
$\Kb \sbs \Gb$ and a neighborhood $\Ub$ of \,$0 \in \Gb$, the
inclusions $\Ub \sbs \Ub_0$, $\Kb_0 \sbs \Kb$ imply that, for
every finite $(\Kb,\Ub)$ approximation $\eta\:G \to \Gb$ and
each $(\Ub,\delta)$ lifting
\,$f \in \LL^1(G,\eta,\Ub,\Kb,\delta,\Ub,\delta)$ of \,$\fb$
on $\Kb$, we have
$$
\bigl|\wh\fb(\cchi) - \FF_{\eta,d}(f)(\cchi)\bigr| \le \alfa
$$
for all \,$\cchi \in \Gamb_0$. Obviously, $\delta$ can be chosen
as small as we please.

Let $(\Kb,\Ub)$, $(\Gamb,\Omb)$ be $\alfa$-adjoint pairs obtained
in the same way as in Example~3.4.7. Using Proposition~3.4.2 we can
ensure that if \,$\Ub$ is small enough and $\Kb$ is big enough,
then there is some $(\Ub,\delta)$ lifting
\,$f \in \LL^1(G,\eta,\Ub,\Kb,\delta,\Ub,\delta)$ of \,$\fb$
on $\Kb$.

If \,$\phi\:\dG \to \dGb$ is a $(\Gamb,\Omb)$ approximation of
\,$\dGb$, $\alfa$-adjoint to $\eta$, then, in the same way as in
the previous Example, we could show that
$$
\bigl|\wh\fb(\cchi) - \wh f(\gama)\bigr| \le
\bigl(2\LV\fb\RV_{_1} + 1\bigr)\alfa + 2\delta
\le \bigl(2\LV\fb\RV_{_1} + 3\bigr)\alfa
$$
for $\cchi \in \Gamb_0$, $\gama \in \dG$, whenever
$\phi(\gama) \in \cchi\,\Omb$, and $\delta \le \alfa$. However,
in order to derive some estimate for $\wh f(\gama)$ when
$\phi(\gama) \in \dGb$ is ``remote'' in a sense, we have to
make use of some additional tools.

Let's start by fixing some $t \in (0,1)$. According to
Proposition~1.4.3 and Lemma~4.1 from \cite{ZZ}, we can choose $f$
subject to \,$t \LV\fb\RV{_1} < \LV f\RV_{_1} \le \LV\fb\RV_{_1}$.
Since
$$
\LV f_a - f\RV_{_1} \le \delta
$$
for $a \in \eta^{-1}[\Ub]$, from Corollary~2.2.6 we get the inclusion
\,$\Spec_t(f) \sbs \Bohr_\alfa\bigl(\eta^{-1}[\Ub]\bigr)$, whenever
$$
\delta \le 2t \LV f\RV_{_1}\sin\frac{\alfa}{2}\,.
$$
Due to our choice of $f$, the last inequality is guaranteed for
$$
\delta \le 2t^2 \LV\fb\RV_{_1}\sin\frac{\alfa}{2}\,.
$$

By the virtue of Theorem~2.5.8 we can assume that the approximations
$\eta$, $\phi$ form a strongly $(\alfa,\delta)$-adjoint pair,
witnessed by some neighborhoods $\Vb \sbs \Ub$, $\Ypb \sbs \Omb$.
Then
$$
\Spec_t(f) \sbs \Bohr_\alfa\bigl(\eta^{-1}[\Ub]\bigr)
\sbs \phi^{-1}\bigl[\Bohr_\delta(\Vb)\bigr]\,,
$$
thus for
$\gama \in \dG \sms \phi^{-1}\bigl[\Bohr_\delta(\Vb)\bigr]$
we have $\gama \nin \Spec_t(f)$, hence
$$
\bigl|\wh f(\gama)\bigr| < t \LV f\RV_{_1} \le t \LV\fb\RV_{_1}\,.
$$
Choosing $\alfa$, and subsequently also $\delta$ small enough,
we can guarantee that
$$
\bigl|\wh\fb(\cchi) - \wh f(\gama)\bigr| \le \eps
$$
for $\cchi \in \Gamb_0$, $\gama \in \phi^{-1}\bigl[\cchi\,\Omb\bigr]$,
and
$$
\bigl|\wh f(\gama)\bigr| \le \eps
$$
for $\gama \in \phi^{-1}\bigl[\dGb \sms \Bohr_\delta(\Vb)\bigr]$.
Taking $\Vb$ small enough, we can arrange that $\Gamb_0\,\Omb$ or
even $\Gamb$ is contained in $\Bohr_\delta(\Vb)$. In any case,
however, we cannot exclude that the set
$\Delb = \Bohr_\delta(\Vb) \sms \Gamb_0\,\Omb$ is nonempty.
Then the values \,$\wh f(\gama)$ for $\gama \in \phi^{-1}[\Delb]$
(if any) are out control of the present approach. This contrasts
the purely existential Theorem~3.4.5, assuring that, at least in
principle, it is possible to maintain full control of \,$\wh f$.
\endexample

\remark{Final remarks}
As the subspace $\Cbc(\Gb)$ is dense in $\Lb^1(\Gb)$, it would be
sufficient in some sense to deal with continuous functions
$\fb \in \Lb^1(\Gb)$ in the last Example~3.4.8. In that case the
composition $f = \fb \co \eta\:G \to \bC$ can be taken as the
approximate lifting of \,$\fb$.

If \,$1 < p \le 2 \le q < \infty$ are dual exponents, then the
continuous functions $\fb \in \Lb^1(\Gb) \cap \Lbp(\Gb)$
form a dense subspace in $\Lbp(\Gb)$ and their Fourier transforms
$\wh\fb \in \Lbq\bigl(\dGb\bigr)$ are continuous, as well. Thus
Example~3.4.8 can be used for the sake of approximation of the
generalized Fourier transform $\Lbp(\Gb) \to \Lbq\bigl(\dGb\bigr)$,
at the same time. That's why there's no need  to elaborate a similar
explicit accompanying example to Theorem~3.4.6, here. Realizing
how cumbersome and laborious it would inevitably be in such a general
setting, we believe that the reader will forgive us.

If \,$\eta\:G \to \Gb$, $\phi\:\dG \to \dGb$ are ``sufficiently well
adjoint'' approximations, then the function $\fb$ and its Fourier
transform $\wh\fb$, once both continuous, admit ``simultaneous''
approximate liftings: $\fb$ by means of the composition
$f = \fb \co \eta\:G \to \bC$ and  $\wh\fb$ by its discrete Fourier
transform $\wh f\:\dG \to \bC$ which is ``close'' to the composition
$\wh\fb \co \phi$, forming another approximate lifting of \,$\wh\fb$.

As pointed out by Gordon \cite{Go2} for the Fourier-Plancherel
transform $\Lb^2(\Gb) \to \Lb^2\bigl(\dGb\bigr)$, the class of
functions $\fb \in \Lbp(\Gb)$ liftable by the composition
$\fb \co \eta$ contains even more general functions, namely the
{\it Riemann integrable\/} ones, i.e., functions continuous almost
everywhere with respect to the Haar measure on $\Gb$.
\endremark

\newpage


\enddocument